\documentclass[a4paper,12pt,twoside]{article}
\usepackage{graphicx}\usepackage{a4wide}
\usepackage{amssymb,amsmath,amsthm}
\usepackage{subfigure,color,colordvi,graphicx,xcolor}
\usepackage[only,llbracket,rrbracket,interleave]{stmaryrd}
\usepackage[sort,compress]{cite}

\newcommand{\bm}[1]{\text{\boldmath $#1$\unboldmath}}

\newcommand{\bddot}{\operatorname{\bm{:}}}
\newcommand{\curl}{\operatorname{curl}}
\newcommand{\curlV}{\operatorname{curl}_{\texttt{V}}}
\newcommand{\vect}[1]{\mathbf{#1}}
\newcommand{\mat}[1]{\mathbf{#1}}

\newcommand{\grad}{\bm{\nabla}}
\newcommand{\RR}{\mathbb{R}}

\newcommand{\SSym}{\mathbb{S}}

\newcommand{\VhHat}{\ensuremath{\mathcal{\hat{V}}^h}}
\newcommand{\Vh}{\ensuremath{\mathcal{V}^h}}

\newcommand{\eltwo}{\ensuremath{\mathcal{L}_2}}

\newcommand{\Assem}{\mbox{\textsf{\textbf{\Large A}}}}

\newcommand{\nen}  {\ensuremath{\texttt{n}_{\texttt{en}}}}
\newcommand{\nfn}  {\ensuremath{\texttt{n}_{\texttt{fn}}}}

\newcommand{\nsd}  {\ensuremath{\texttt{n}_{\texttt{sd}}}}
\newcommand{\msd}  {\ensuremath{\texttt{m}_{\texttt{sd}}}}
\newcommand{\numel}{\ensuremath{\texttt{n}_{\texttt{el}}}}
\newcommand{\numfa}{\ensuremath{\texttt{n}_{\texttt{fa}}}}

\DeclareMathOperator{\tr}{tr}

\newcommand{\hv}{\hat{v}}
\newcommand{\hu}{\hat{u}}
\newcommand{\bu}{\bm{u}}
\renewcommand{\bf}{\bm{f}}
\newcommand{\bt}{\bm{t}}

\newcommand{\bhw}{\widehat{\bw}}
\newcommand{\bhu}{\widehat{\bu}}

\newcommand{\bv}{\bm{v}}

\newcommand{\bn}{\bm{n}}
\newcommand{\bx}{\bm{x}}
\newcommand{\bL}{\bm{L}}

\newcommand{\bTens}{\bm{\varsigma}}
\newcommand{\bTensV}{\bm{\varsigma}_{\texttt{V}}}
\newcommand{\btau}{\bm{\tau}}

\newcommand{\nDeg}{\ensuremath{k}}
\newcommand{\Pk}{\ensuremath{\mathcal{P}^{\nDeg}}}
\newcommand{\Hdiv}{H(\operatorname{div})}

\newcommand{\Insd}{\mat{I}_{\nsd}}
\newcommand{\Imsd}{\mat{I}_{\msd}}
\newcommand{\stress}{\bm{\sigma}}
\newcommand{\defo}[1]{\bm{\varepsilon}(#1)}
\newcommand{\gradS}{\bm{\nabla}_{\texttt{S}}}
\newcommand{\stressV}{\bm{\sigma}_{\texttt{V}}}
\newcommand{\strainV}{\bm{\varepsilon}_{\texttt{V}}}
\newcommand{\bD}{\mat{D}}
\newcommand{\bDHalf}{\bD^{1/2}}
\newcommand{\bN}{\mat{N}}
\newcommand{\bR}{\mat{R}}
\newcommand{\bT}{\mat{T}}
\newcommand{\nrr}  {\ensuremath{\texttt{n}_{\texttt{rr}}}}
\newcommand{\bg}{\bm{g}}

\newcommand{\jump}[1]{\llbracket #1\rrbracket}

\newcommand{\nodaluV}{\vect{u}}
\newcommand{\nodaluhV}{\hat{\vect{u}}}
\newcommand{\nodalLV}{\vect{L}}
\newcommand{\bxi}{\bm{\xi}}
\newcommand{\bet}{\bm{\eta}}

\newcommand{\nipe}{\texttt{n}_{\texttt{ip}}^{\texttt{e}}}
\newcommand{\sumge}{\sum_{\texttt{g}=1}^{\nipe}}

\newcommand{\bxige}{\bxi^\texttt{e}_\texttt{g}}
\newcommand{\wge}{w^\texttt{e}_\texttt{g}}
\newcommand{\bJ}{\mat{J}}

\newcommand{\nipf}{\texttt{n}_{\texttt{ip}}^{\texttt{f}}}
\newcommand{\sumgf}{\sum_{\texttt{g}=1}^{\nipf}}

\newcommand{\bxigf}{\bxi^\texttt{f}_\texttt{g}}
\newcommand{\wgf}{w^\texttt{f}_\texttt{g}}

\newcommand{\Nmat}{\bm{\mathcal{N}}}
\newcommand{\Mmat}{\bm{\mathcal{M}}}
\newcommand{\NmatHat}{\bm{\widehat{\mathcal{N}}}}

\newcommand{\bw}{\bm{w}}

\newenvironment{keywords}{\begin{quote}\emph{\textbf{Keywords:}}}{\end{quote}}

\newtheorem{lemma}{Lemma}
\theoremstyle{definition}

\newtheorem{remark}{Remark}

\newcommand{\etal}{et al.\ }

\graphicspath{{figs/}}

\begin{document}
\title{A super--convergent hybridisable discontinuous Galerkin method for linear elasticity}

\author{Ruben Sevilla\\[-1ex]
             \small Zienkiewicz Centre for Computational Engineering, \\[-1ex]
             \small College of Engineering, Swansea University, Wales, UK \\[1em]
             Matteo Giacomini, Alexandros Karkoulias, and Antonio Huerta\\[-1ex]
             \small Laboratori de C\`alcul Num\`eric (LaC\`aN), \\[-1ex]
             \small ETS de Ingenieros de Caminos, Canales y Puertos, \\[-1ex]
             \small Universitat Polit\`ecnica de Catalunya, Barcelona, Spain}
\date{\today}
\maketitle

\begin{abstract}
The first super--convergent hybridisable discontinuous Galerkin (HDG) method for linear elastic problems capable of using the same degree of approximation for both the primal and mixed variables is presented. The key feature of the method is the strong imposition of the symmetry of the stress tensor by means of the well--known and extensively used Voigt notation, circumventing the use of complex mathematical concepts to enforce the symmetry of the stress tensor either weakly or strongly. A novel  procedure to construct element--by--element a super--convergent post--processed displacement is proposed. Contrary to other HDG formulations,  the methodology proposed here is able to produce a super--convergent displacement field for low order approximations. The resulting method is robust and locking--free in the nearly--incompressible limit. An extensive set of numerical examples is utilised to provide evidence of the optimality of the method and its super--convergent properties in two and three dimensions and for different element types.
\end{abstract}


\section{Introduction}
\label{sc:Intro}

The numerical approximation of the linear elasticity equation presents several difficulties as highlighted by the extensive literature available on the topic (cf. e.g.~\cite{MR2373954, braess2001finite}).
In particular, locking phenomena in nearly incompressible and incompressible materials, construction of stable pairs of finite elements and strong enforcement of the symmetry of the stress tensor in mixed formulations, accurate computation of the stresses (classically recovered from the displacement field via numerical differentiation) and post--process procedures to improve the quality of the approximate displacement fields are some of the subjects that have attracted the attention of the scientific community over the last 40 years.

It is well--known that an accurate approximation of the linear elastic problem for nearly incompressible materials requires the discrete space in which the solution is sought to be rich enough to describe non--trivial divergence--free vector fields.
Within this context, the primal formulation where the displacement field is the sole unknown fails to provide a locking--free approximation using conforming Lagrangian finite element functions. In fact, Brenner and Sung~\cite{MR1140646}  proposed a possible remedy by means of the nonconforming Crouzeix--Raviart element~\cite{MR0343661}.
In order to circumvent this issue, two main approaches have been proposed in the literature. On the one hand, mixed formulations in which both the displacement field and the stress tensor act as unknowns of a saddle point problem~\cite{reissner1950on}. On the other hand, discontinuous Galerkin (DG) discretisations in which the approximate displacement field is sought in a bigger space and the variational formulation of the problem is modified to account for the jumps of the discrete displacement field across the element interfaces~\cite{MR1886000, MR1972650}.

Starting from the seminal paper by Reissner~\cite{reissner1950on}, mixed variational formulations of the linear elasticity equation have known a great success in the scientific community. 
The solution of the resulting saddle point problem provides an approximation of both the displacement field and the stress tensor that is not retrieved as a post--processed quantity (with a consequent loss of precision) as in the primal formulation~\cite{MR1077658}.
A major drawback of the mixed formulation lies in the difficulty of constructing a pair of finite element spaces that fulfil the requirements of Brezzi's theory~\cite{Brezzi1974} to guarantee the stability of the method. 
More precisely, concurrently imposing the balance of momentum by seeking a stress tensor in $\Hdiv$ (i.e.\ a square--integrable tensor with square--integrable row--wise divergence) and the balance of angular momentum by enforcing its symmetry proved to be an extremely difficult task~\cite{MR2449101}.
Stemming from the pioneering work by Fraejis de Veubeke~\cite{Veubeke1975}, a first approach discussed in the literature relies on maintaining the $\Hdiv$-conformity of the stress tensor while its symmetry is relaxed (cf. e.g.~\cite{MR553347, MR946367, MR1082821}). 
Among the most successful approaches, the so--called PEERS element by Arnold, Brezzi and Douglas Jr.~\cite{MR840802} introduced a Lagrange multiplier as extra variable to account for the symmetry constraint, see also~\cite{MR1005064, MR834332, MR954768, MR956898, MR859922, MR1464150}. 
In 2002, the first stable pair of finite element spaces for the discretisation of the mixed formulation of the linear elasticity equation with $\Hdiv$-conforming strongly--enforced symmetric stress tensor in two dimensions was proposed by Arnold and Winther~\cite{MR1930384}.
The corresponding three dimensional case is discussed in~\cite{Adams2005, 10.2307/40234556}.
Nevertheless, the construction of these finite element pairs is based on non--trivial techniques of exterior calculus~\cite{MR2269741} and results in a large number of degrees of freedom per element (the lowest--order approximation of the stress tensor features 24 degrees of freedom on a triangle and 162 on a tetrahedron) making their application to complex problems unfeasible.
More recently, an alternative mixed formulation featuring a tangential--continuous displacement field and a normal--normal continuous symmetric stress tensor has been proposed by Pechstein and Sch{\"o}berl~\cite{MR2826472}.
For a detailed discussion on mixed methods, the interested reader is referred to~\cite{brezzi1991mixed}.

An alternative approach to the discretisation of the linear elastic problem focuses on relaxing the $\Hdiv$-conformity of the stress tensor while strongly enforcing its pointwise symmetry.
This results in nonconforming discretisations (cf. e.g.~\cite{MR1977627, MR2484491, MR3163400}).
Moreover, owing to the fact that these methods use polynomial basis functions but no degrees of freedom is located in the vertices of the elements, Gopalakrishnan and Guzm\'{a}n~\cite{MR2831058} show that the resulting nonconforming approximation may be efficiently implemented via hybridisation.
Nevertheless, the convergence rate of the stress tensor is sub--optimal when using this nonconforming discretisation.
Among nonconforming discretisations, the DG method has experienced a great success in recent years.
The interest in DG methods for the linear elastic problem is motivated by their high--order convergence properties and their flexibility in performing local $h$-- and $p$--adaptivity. 
Beside the aforementioned works by Hansbo and Larson~\cite{MR1886000, MR1972650}, see also~\cite{MR2027288, MR2219020}.
The ability of DG methods to efficiently construct a locking--free approximation in nearly incompressible materials has been recently analysed in~\cite{MR2997905} also for the case of heterogeneous media.
In~\cite{MR2995177}, a discontinuous Petrov--Galerkin (DPG) formulation is proposed to simultaneously approximate the displacement field and the symmetric stress tensor and is shown to be $hp$--optimal.

More recently, novel discretisation techniques inspired by the previously discussed ones have been proposed.
The local discontinuous Galerkin (LDG) method~\cite{MR2220915} is based on a mixed discontinuous Galerkin formulation and provides an exactly incompressible approximation of the displacement field, that is a displacement field which is normal--continuous across inter--element boundaries and pointwise incompressible inside each element.
The method converges optimally for the displacement field whereas the strain tensor and the pressure results are sub--optimal by one order.
Moreover, contrary to the framework discussed by Gopalakrishnan and Guzm\'{a}n~\cite{MR2831058}, the method cannot be hybridised, thus resulting in a considerable number of degrees of freedom for high--order approximations.
Stemming from the work on LDG, Cockburn and co--workers have proposed the hybridisable discontinuous Galerkin (HDG) method whose analysis for the linear elasticity equation is available in~\cite{soon2009hybridizable, MR3340089}.
HDG is based on a mixed discontinuous Galerkin formulation with hybridisation and provides an optimally convergent displacement field with order $k+1$ whereas the strain and stress tensor converge sub--optimally with order $k+1/2$.
The optimal approximation of the stress tensor may be retrieved by adding matrix bubble functions to the discrete space as discussed in~\cite{MR2629995}.
The alternative HDG formulation by Qiu \etal\cite{MR3716189} exploits polynomials of different degrees for the approximation of the displacement field (order $k+1$), its trace (order $k$) and the strain tensor (order $k$).
By introducing a modified definition of the numerical trace, optimal convergence of order $k+1$ for all the unknowns is retrieved.
In~\cite{MR3283758},  Di Pietro and Ern discuss the hybrid high--order (HHO) method which features a nonconforming discretisation based on a pure displacement formulation of the linear elastic problem.
The method leads to a locking--free displacement field and a strongly--symmetric strain tensor, both converging with optimal order.
Similarly to the hybridisation in HDG, the definition of the unknowns on the faces reduces the computational cost associated with the solution of the problem, making HHO suitable for high--order approximations. 

Inspired by the works of Arnold and Brezzi~\cite{MR813687} and Stenberg~\cite{MR954768} on mixed methods, Cockburn and co--workers have investigated several procedures to construct a super--convergent post--processed solution by exploiting both the optimally convergent primal and mixed variables.
Nevertheless, as shown in~\cite{soon2009hybridizable}, the case of linear elasticity experiences a sub--optimal convergence of the strain tensor and consequently a loss of super--convergence for the post--processed displacement field.
A great effort within the HDG community is currently devoted to investigating techniques to remedy this issue and retrieve the super--convergence of the post--processed variable even for low--order approximations.
More precisely, Cockburn and co--workers have recently introduced the concept of $\bm{M}$-decomposition~\cite{MR3081483, arXiv2017BCGF} to construct discrete spaces suitable to retrieve the aforementioned super--convergence property.
This elegant theory guarantees that an HDG approximation for which the local space admits an $\bm{M}$-decomposition provides a locking--free approximate displacement field, an optimally convergent approximate stress tensor and a super--convergent post--processed displacement field obtained via an element--by--element procedure.
Nevertheless, the construction of such spaces is non--trivial and their implementation in existing HDG library is not straightforward.

This work proposes an extremely simple alternative to remedy the loss of optimality of the mixed variable in the HDG formulation of the linear elastic problem by strongly enforcing its symmetry with classical nodal--based discrete interpolation spaces.
First, the Voigt notation for symmetric tensors is recalled and the linear elasticity equation is rewritten strongly enforcing the symmetry of stress tensor (Section~\ref{sc:problem}).
In Section~\ref{sc:HDG}, the HDG framework discussed by Soon \etal \ in~\cite{soon2009hybridizable} is considered and discrete spaces featuring equal order interpolation for all the variables are employed.
By exploiting the retrieved optimal convergence rate of order $k+1$ of the mixed variable (i.e.\ the stress tensor) a novel procedure to derive a super--convergent post--processed displacement field is introduced (Section~\ref{sc:superconvergence}).
The displacement field being identified up to \emph{rigid motions} (three in 2D and six in 3D), an additional set of constraints is required for the displacement field to be unique and different solutions are discussed.
Extensive numerical tests in both two and three dimensions are presented in Section~\ref{sc:examples} to validate the convergence rates of the primal, mixed and post--processed variables, using different types of elements. Special attention is given to the nearly incompressible limit case in which the novel formulation confirms to be locking--free and the optimal convergence rates are preserved.
Section~\ref{sc:Conclusion} summarises the discussed results whereas the implementation details are provided in Appendix~\ref{sc:implementation}.

\section{Problem statement}
\label{sc:problem}

In this section, the governing equations that describe the mechanical behaviour of a deformable solid within the infinitesimal strain theory are introduced and the corresponding formulation using the Voigt notation for symmetric tensors is recalled.
For a complete introduction to this subject, the interested reader is referred to~\cite{gould1993introduction, MR1262126, MR936420}.

\subsection{Strong form of the linear elastic problem} \label{sc:elasticity}

Let $\Omega\subset\mathbb{R}^{\nsd}$ be an open bounded domain in $\nsd$ spatial dimensions with boundary $\partial\Omega = \overline{\Gamma}_D \cup \overline{\Gamma}_N$, $\overline{\Gamma}_D \cap \overline{\Gamma}_N = \emptyset$ and $\Gamma_D$ featuring positive $(\nsd-1)$--dimensional Hausdorff measure. 
The mechanical behaviour of a deformable solid $\Omega$ within the infinitesimal strain theory is the described by
\begin{equation} \label{eq:elasticity}
\left\{\begin{aligned}
-\grad \cdot \stress  &= \bf       &&\text{in $\Omega$,}\\
\stress &= \stress^T       &&\text{in $\Omega$,}\\
\bu &= \bu_D  &&\text{on $\Gamma_D$,}\\
\bn\cdot\stress &= \bg        &&\text{on $\Gamma_N$,}\\
\end{aligned}\right.
\end{equation}
where $\bu$ is the displacement field and $\stress$ is the Cauchy stress tensor.
The elastic structure $\Omega$ under analysis is thus subject to a volume force $\bf$, a tension $\bg$ on the surface $\Gamma_N$ and an imposed displacement $\bu_D$ on $\Gamma_D$.

Equation~\eqref{eq:elasticity} is the strong form of the linear elastic problem and states two conservation laws, namely the balance of momentum and the balance of angular momentum.
Remark that the latter implies the symmetry of the stress tensor, that is $\stress$ belongs to the space $\SSym^{\nsd}$ of $\nsd \times \nsd$ symmetric matrices.
The full set of equations is closed by a material law that describes the relationship among the variables at play and depends on the type of solid under analysis. In particular, a linear elastic material is considered. 
Within this context, the so--called Hooke's law establishes a linear dependency between the stress tensor $\stress$ and the linearised strain tensor $\defo{\bu} := \left( \grad \bu + \grad \bu^T \right)/2$ via the fourth--order tensor $A : \Omega \rightarrow \SSym^{\nsd}$ known as the elasticity tensor.
In this work, only homogeneous isotropic materials are considered, whence the elasticity tensor $A$ depends neither on the spatial coordinate $\bx$ nor on the direction of the main strains.
The mechanical properties of a linear elastic homogeneous isotropic material are determined by the pair $(E,\nu)$, respectively known as Young's modulus and Poisson's ratio (cf. e.g.~\cite{MR1262126}).
Within the range of physically admissible values of these constants (i.e. $\nu \in (-1,0.5)$), the relationship between the stress tensor and the linearised strain tensor reads 
\begin{equation} \label{eq:elas_hooke_lame} 
\stress = A \defo{\bu} = \frac{E}{1+\nu} \defo{\bu} + \frac{E \nu}{(1+\nu)(1-2\nu)} \tr(\defo{\bu})\Insd, 
\end{equation}
where $\Insd$ is the $\nsd \times \nsd$ identity matrix and $\tr(\cdot) := \cdot \bddot \Insd$ is the trace operator, being $\bddot$ the Frobenius product, also known as double contraction.
For the purpose of the current work, only non--auxetic materials are considered, that is the Poisson's ratio $\nu$ is assumed to be non--negative.

By plugging~\eqref{eq:elas_hooke_lame} into~\eqref{eq:elasticity}, the stress tensor may be expressed in terms of the displacement field and the pure displacement formulation of the linear elastic problem is retrieved:
\begin{equation} \label{eq:elasticityDisplacement}
\left\{\begin{aligned}
-\grad \cdot \left( A \defo{\bu} \right)  &= \bf       &&\text{in $\Omega$,}\\
\bu &= \bu_D  &&\text{on $\Gamma_D$,}\\
\bn\cdot \left( A \defo{\bu} \right) &= \bg        &&\text{on $\Gamma_N$.}\\
\end{aligned}\right.
\end{equation}

\begin{remark}
	The elasticity tensor exists and is invertible as long as $\nu < 0.5$. 
	It is straightforward to observe that when $\nu \rightarrow 0.5$, the divergence of the displacement field in~\eqref{eq:elas_hooke_lame} has to vanish, that is, the material under analysis is incompressible.
	\label{remark:incompressibility}
\end{remark}
This case cannot be properly handled by the pure displacement formulation since the elasticity tensor deteriorates and $A$ fails to exist in the incompressible limit, thus preventing the stress tensor to be expressed in terms of the displacement field.
A possible remedy is represented by mixed formulations in which both the displacement field and the stress tensor act as unknown of the problem.
%
%
The associated first--order problem is thus obtained by considering the following system of equations:
\begin{equation} \label{eq:elasticitySystem}
\left\{\begin{aligned}
-\grad \cdot \stress  &= \bf       &&\text{in $\Omega$,}\\
\stress &= \stress^T       &&\text{in $\Omega$,}\\
\stress &= A \defo{\bu}      &&\text{in $\Omega$,}\\
\bu &= \bu_D  &&\text{on $\Gamma_D$,}\\
\bn\cdot\stress &= \bg       &&\text{on $\Gamma_N$.}\\
\end{aligned}\right.
\end{equation}

\subsection{Strong enforcement of the symmetry of the stress tensor} \label{sc:Voigt}

Consider the classical theory of linear elasticity~\cite{FishBelytschko2007}.
Let $\bu := \bigl[ u_i \bigr]^T \in \RR^{\nsd}, \ i=1,\dotsc,\nsd$ be the vector field describing the displacement.
The strain tensor may be divided into its diagonal components (namely, the extensional strains $\varepsilon_{ii}$) and its off--diagonal terms $\gamma_{ij}$ known as shear strains
\begin{equation} \label{eq:NormalShear}
\varepsilon_{ii} := \frac{\partial u_i}{\partial x_i} , \quad \gamma_{ij} := \frac{\partial u_i}{\partial x_j} + \frac{\partial u_j}{\partial x_i}, \quad \text{for } i,j = 1,\dotsc,\nsd .
\end{equation}
Owing to its symmetry, only three components in 2D (two extensional and one shear strains) and six components in 3D (three extensional and three shear strains) need to be stored.
More precisely, according to the so--called Voigt notation, the components of the strain may be arranged as a column vector in $\RR^{\msd}$ as follows:
\begin{equation} \label{eq:strainVoigt}
\strainV :=\begin{cases}
\bigl[\varepsilon_{11} ,\; \varepsilon_{22} ,\; \gamma_{12} \bigr]^T
&\text{in 2D,} \\
\bigl[\varepsilon_{11} ,\; \varepsilon_{22} ,\; \varepsilon_{33} ,\; \gamma_{12} ,\; \gamma_{13} ,\; \gamma_{23} \bigr]^T
&\text{in 3D,} 
\end{cases}
\end{equation}
where $\msd = \nsd(\nsd+1)/2$.

\begin{remark}
	The linearised strain tensor $\defo{\bu} \in \SSym^{\nsd}$ differs from its Voigt counterpart $\strainV \in \RR^{\msd}$ by a factor $1/2$ in the shear components, that is:
	\begin{equation} \label{eq:strainClassicVoigt}
	\defo{\bu} :=\begin{cases}
	\begin{bmatrix}
	\varepsilon_{11} & \gamma_{12}/2 \\
	\gamma_{12}/2 & \varepsilon_{22}
	\end{bmatrix} 
	&\text{in 2D,} \\
	\begin{bmatrix}
	\varepsilon_{11} & \gamma_{12}/2 & \gamma_{13}/2 \\
	\gamma_{12}/2 & \varepsilon_{22} & \gamma_{23}/2 \\
	\gamma_{13}/2 & \gamma_{23}/2 & \varepsilon_{33}
	\end{bmatrix}
	&\text{in 3D.} 
	\end{cases}
	\end{equation}
\end{remark}

Following the framework described by Fish and Belytschko in~\cite{FishBelytschko2007}, the matrix $\gradS \in \RR^{\msd \times \nsd}$ accounting for the symmetric gradient operator is introduced:
\begin{equation} \label{eq:symmGrad}
\gradS :=\begin{cases}
\begin{bmatrix}
\partial/\partial x_1 & 0 & \partial/\partial x_2 \\
0 & \partial/\partial x_2 & \partial/\partial x_1
\end{bmatrix}^T
&\text{in 2D,} \\
\begin{bmatrix}
\partial/\partial x_1 & 0 & 0 & \partial/\partial x_2 & \partial/\partial x_3 & 0 \\
0 & \partial/\partial x_2 & 0 & \partial/\partial x_1 & 0 & \partial/\partial x_3 \\
0 & 0 & \partial/\partial x_3 & 0 & \partial/\partial x_1 & \partial/\partial x_2
\end{bmatrix}^T
&\text{in 3D.} 
\end{cases}
\end{equation}
Thus, the components of $\strainV$ may be expressed in terms of the displacements $\bu$ by means of a single matrix equation $\strainV = \gradS \bu$.

In a similar fashion, within the stress tensor $\stress$ two normal components $\sigma_{ii}$ and one shear component $\tau_{ij}$ in 2D (respectively, three and three in 3D) may be identified owing to the symmetry arising form the balance of angular momentum (cf. equation~\eqref{eq:elasticitySystem}).
Thus, according to Voigt notation, the stress tensor may be written as the following column vector in $\RR^{\msd}$:
\begin{equation} \label{eq:stressVoigt}
\stressV :=\begin{cases}
\bigl[\sigma_{11} ,\; \sigma_{22} ,\; \tau_{12} \bigr]^T
&\text{in 2D,} \\
\bigl[\sigma_{11} ,\; \sigma_{22} ,\; \sigma_{33} ,\; \tau_{12} ,\; \tau_{13} ,\; \tau_{23} \bigr]^T
&\text{in 3D.} 
\end{cases}
\end{equation}
%

\subsection{The linear elastic problem using Voigt notation} \label{sc:ElastVoigt}

In this section, the previously introduced Voigt notation is exploited to rewrite the linear elastic problem~\eqref{eq:elasticitySystem} by strongly enforcing the symmetry of the stress tensor.
The second equation in~\eqref{eq:elasticitySystem} is thus verified in a straightforward manner.
The balance of momentum may be rewritten as a matrix equation by exploiting the notation introduced in~\eqref{eq:symmGrad} for the symmetric gradient operator:
\begin{equation} \label{eq:momentumVoigt}
-\gradS^T \stressV = \bf .
\end{equation}
Moreover, the constitutive equation~\eqref{eq:elas_hooke_lame} may be expressed as $\stressV = \bD \strainV$, where $\bD$ is an $\msd \times \msd$ symmetric positive definite matrix describing the generalised Hooke's law.
\begin{remark}
	In two dimensions, the structure of the matrix $\bD$ depends on the assumption made to simplify the three dimensional model. 
	On the one hand, according to the \emph{plain strain} model, the body is thick with respect to the plane $x_1 x_2$ and consequently the extensional strain along $x_3$ and the shear strains $\gamma_{i3}, i =1,2$ vanish.
	On the other hand, the \emph{plane stress} model is based on the assumption that the body is thin relative to the dimensions in the $x_1 x_2$ plane. Thus, no loads are applied along the $x_3$ direction and the component $\sigma_{33}$ of the stress tensor is assumed to vanish.
\end{remark}
\begin{equation} \label{eq:LawVoigt}
\bD{:=}\begin{cases}
\displaystyle\frac{E}{(1+\nu)(1-2\nu)}
\begin{bmatrix}
1-\nu & \nu & 0 \\
\nu & 1-\nu & 0 \\
0 & 0 & (1-2\nu)/2
\end{bmatrix}
&\text{in 2D (plain strain),} \\
\displaystyle\frac{E}{1- \nu^2}
\begin{bmatrix}
1 & \nu & 0 \\
\nu & 1 & 0 \\
0 & 0 & (1-\nu)/2
\end{bmatrix}
&\text{in 2D (plain stress),} \\
\displaystyle\frac{E}{(1+\nu)(1-2\nu)}
\begin{bmatrix}
1-\nu & \nu & \nu & \\
\nu & 1-\nu & \nu & \bm{0}_{\nsd} \\
\nu & \nu & 1-\nu & \\
& \bm{0}_{\nsd} & & \tfrac{1-2\nu}{2}\Insd
\end{bmatrix}
&\text{in 3D.} 
\end{cases}
\end{equation}
Following the same rationale discussed above, an $\nsd \times \msd$ matrix accounting for the normal direction to the boundary is introduced:
\begin{equation} \label{eq:normalVoigt}
\bN :=\begin{cases}
\begin{bmatrix}
n_1 & 0 & n_2 \\
0 & n_2 & n_1
\end{bmatrix}^T
&\text{in 2D,} \\
\begin{bmatrix}
n_1 & 0 & 0 & n_2 & n_3 & 0\\
0 & n_2 & 0 & n_1 & 0 & n_3 \\
0 & 0 & n_3 & 0 & n_1 & n_2
\end{bmatrix}^T
&\text{in 3D,} 
\end{cases}
\end{equation}
and the matrix counterpart of the traction boundary conditions is imposed on $\Gamma_N$, that is $\bN^T \stressV = \bg$.

Hence, the linear elastic problem~\eqref{eq:elasticitySystem} using Voigt notation reads as follows:
\begin{equation} \label{eq:elasticitySystemVoigt}
\left\{\begin{aligned}
-\gradS^T \stressV  &= \bf       &&\text{in $\Omega$,}\\
\stressV &= \bD \strainV      &&\text{in $\Omega$,}\\
\bu &= \bu_D  &&\text{on $\Gamma_D$,}\\
\bN^T \stressV &= \bg        &&\text{on $\Gamma_N$.}\\
\end{aligned}\right.
\end{equation}

\subsection{Generalised Gauss's and Stokes' theorems} \label{sc:DivTheoVoigt}

In order to state the variational formulation of Equation~\eqref{eq:elasticitySystemVoigt}, a counterpart of the classical Gauss's theorem using the Voigt matrices introduced in the previous sections is required. The following result holds:
\begin{lemma}[Generalised Gauss's theorem] \label{theo:Gauss}
	Consider a vector $\bv \in \RR^{\nsd}$ and a symmetric tensor $\bTens \in \SSym^{\nsd}$ whose counterpart in Voigt notation is $\bTensV$. It holds:
	\begin{equation} \label{eq:GenGauss}
	\int_{\partial\Omega} \left( \bN^T \bTensV \right) \cdot \bv \ d\Gamma = \int_{\Omega} \bTensV \cdot \left( \gradS \bv \right) d\Omega + \int_{\Omega} \left( \gradS^T \bTensV \right) \cdot \bv \ d\Omega .
	\end{equation}
	\begin{proof}
		Rewrite each term in~\eqref{eq:GenGauss} in terms of the operators associated with the matrices introduced by the Voigt notation:
		\begin{subequations} \label{eq:GaussProof}
			\begin{align}
			\int_{\partial\Omega} \left( \bN^T \bTensV \right) \cdot \bv \ d\Gamma &= \int_{\partial\Omega} \left( \bn \cdot \bTens \right) \cdot \bv \ d\Gamma ,
			\label{eq:gaussFirst}\\
			\int_{\Omega} \bTensV \cdot \left( \gradS \bv \right) \ d\Omega & = \int_{\Omega} \bTens \bddot \defo{\bv} \ d\Omega ,
			\label{eq:gaussSecond}\\
			\int_{\Omega} \left( \gradS^T \bTensV \right) \cdot \bv  \ d\Omega &= \int_{\Omega} \left( \grad \cdot \bTens \right) \cdot \bv \ d\Omega .
			\label{eq:gaussThird}
			\end{align}
		\end{subequations}
		%
		%
		By summing the right hand sides of~\eqref{eq:GaussProof}, the classical statement of Gauss's theorem is retrieved and consequently~\eqref{eq:GenGauss} holds.
	\end{proof}
\end{lemma}

The aforementioned result allows to derive the formulation of the HDG method which will be discussed in Section~\ref{sc:HDG}.
Moreover, in Section~\ref{sc:superconvergence}, a novel post--process procedure of the HDG solution which relies on a condition on the $\curl$ operator will be introduced.
In order to properly state the aforementioned results, first consider the infinitesimal rotation of a vector field using Voigt notation.
Consider $\bR \in \RR^{\nrr \times \nsd}$, with $\nrr$ the number of rigid body rotations in the space (one in 2D and three in 3D).
Within this rationale, $\curl(\bu) := \grad \times \bu$ may be written as the matrix equation $\curlV(\bu) = \bR \bu$, where 
\begin{equation} \label{eq:curlVoigt}
\bR :=\begin{cases}
\bigl[-\partial/\partial x_2 ,\; \partial/\partial x_1 \bigr]
&\text{in 2D,} \\
\begin{bmatrix}
0 & -\partial/\partial x_3 & \partial/\partial x_2 \\
\partial/\partial x_3 & 0 & -\partial/\partial x_1 \\
-\partial/\partial x_2 & \partial/\partial x_1 & 0
\end{bmatrix}
&\text{in 3D.} 
\end{cases}
\end{equation}
\begin{remark}
	Recall that the $\curl$ of a vector field $\bv \in \RR^2$ exists solely as a scalar quantity, namely
	\begin{equation} \label{eq:curl2D}
	\grad \times \bv = \frac{\partial v_2}{\partial x_1} - \frac{\partial v_1}{\partial x_2} .
	\end{equation}
	Nevertheless, by embedding $\bv$ in $\RR^3$ and setting its third component equal to zero, the $\curl$ may be interpreted as a vector pointing entirely in the direction $x_3$ with magnitude given by $\bR \bv$, that is, the value on the right hand side of Equation~\eqref{eq:curl2D}.
\end{remark}
Moreover, consider the following matrix $\bT \in \RR^{\nrr \times \nsd}$ describing the tangent direction to the boundary of $\Omega \subset \RR^{\nsd}$, that is a tangent line in 2D and a tangent surface in 3D:
\begin{equation} \label{eq:tangentVoigt}
\bT :=\begin{cases}
\bigl[n_2 ,\; -n_1 \bigr]^T
&\text{in 2D,} \\
\begin{bmatrix}
0 & -n_3 & n_2 \\
n_3 & 0 & -n_1 \\
-n_2 & n_1 & 0
\end{bmatrix}
&\text{in 3D.} 
\end{cases}
\end{equation}
As previously done for the Gauss's theorem, a generalised Stokes' theorem using the Voigt matrices is stated:
\begin{lemma}[Generalised Stokes' theorem] \label{theo:Stokes}
	Consider a vector $\bv \in \RR^{\nsd}$. It holds:
	\begin{equation} \label{eq:GenStokes}
	\int_{\Omega} \bR \bv \ d\Omega = \int_{\partial\Omega} \bv^T \bT \ d\Gamma .
	\end{equation}
	\begin{proof}
		Following the same rationale used in Lemma~\ref{theo:Gauss}, each term in~\eqref{eq:GenStokes} may be rewritten as follows:
		\begin{subequations} \label{eq:StokesProof}
			\begin{align}
			\int_{\Omega} \bR \bv \ d\Omega &= \int_{\Omega} \grad \times \bv \ d\Omega ,
			\label{eq:stokesFirst} \\
			\int_{\partial\Omega} \bv^T \bT \ d\Gamma &= \int_{\partial\Omega} \bv \cdot \bt \ d\Gamma ,
			\label{eq:stokesSecond}
			\end{align}
		\end{subequations}
		where $\bt$ is the tangential direction to the boundary $\partial\Omega$.
		By plugging~\eqref{eq:StokesProof} into~\eqref{eq:GenStokes}, the classical statement of Stokes' theorem is retrieved and consequently~\eqref{eq:GenStokes} holds.
	\end{proof}
\end{lemma}

\section{Hybridisable discontinuous Galerkin formulation}
\label{sc:HDG}

Consider a partition of the domain $\Omega$ in $\numel$ disjoint subdomains $\Omega_e$ with boundaries $\partial\Omega_e$.
The internal interface $\Gamma$ is defined as
\begin{equation} \label{eq:skeleton}
\Gamma := \left[ \bigcup_{e=1}^{\numel} \partial\Omega_e \right] \setminus \partial\Omega .
\end{equation}

The second--order elliptic problem of Equation~\eqref{eq:elasticitySystemVoigt} can be written in mixed form, in the so--called broken computational domain, as a system of first--order equations, namely
\begin{equation} \label{eq:elasticityBrokenFirstOrder}
\left\{\begin{aligned}
\bL + \bDHalf \gradS \bu &= \bm{0}    &&\text{in $\Omega_e$, and for $e=1,\dotsc ,\numel$,}\\	
\gradS^T \bDHalf \bL &= \bm{f}          &&\text{in $\Omega_e$, and for $e=1,\dotsc ,\numel$,}\\
\bu &= \bu_D     &&\text{on $\Gamma_D$,}\\
\bN^T \bDHalf \bL  &= -\bg         &&\text{on $\Gamma_N$,}\\
\jump{\bu \otimes \bn} &= \bm{0}  &&\text{on $\Gamma$,}\\
\jump{\bN^T \bDHalf \bL} &= \bm{0}  &&\text{on $\Gamma$,}\\
\end{aligned} \right.
\end{equation}
where $\jump{\cdot}$ denotes the \emph{jump} operator, defined along each portion of the interface according to~\cite{AdM-MFH:08} as the sum of the values from the element on the right and left, say $\Omega_e$ and $\Omega_l$:
\begin{equation}
\jump{\odot} = \odot_e + \odot_l .
\end{equation}
Therefore, the last two equations in~\eqref{eq:elasticityBrokenFirstOrder} enforce the continuity of respectively the primal variable - i.e. the displacement field - and the normal trace of the stress across the interface $\Gamma$.

\subsection{Strong form of the local and global problems}
\label{sc:HDGstrong}

The HDG formulation solves the problem of Equation~\eqref{eq:elasticityBrokenFirstOrder} in two stages~\cite{Jay-CGL:09,Cockburn-CDG:08,Nguyen-NPC:09,Nguyen-NPC:09b,Nguyen-NPC:10,Nguyen-NPC:11}. First a local pure Dirichlet problem is defined to compute $(\bL_e,\bu_e)$ element--by--element in terms of the unknown hybrid variable $\bhu$, namely
\begin{equation} \label{eq:elasticityStrongLocal}
\left\{\begin{aligned}
\bL_e + \bDHalf \gradS \bu_e &= \bm{0}    &&\text{in $\Omega_e$}\\	
\gradS^T \bDHalf \bL_e &= \bm{f}          &&\text{in $\Omega_e$}\\
\bu_e &= \bu_D     &&\text{on $\partial\Omega_e \cap \Gamma_D$,}\\
\bu_e &= \bhu  &&\text{on $\partial\Omega_e \setminus \Gamma_D$,}\\
\end{aligned} \right.
\end{equation}
for $e=1,\dotsc ,\numel$.

Second, the global problem is defined to determine the hybrid variable (i.e. the trace of the displacement field on the mesh skeleton $\Gamma \cup \Gamma_N$), namely
\begin{equation} \label{eq:elasticityStrongGlobal}
\left\{\begin{aligned}
\jump{\bu \otimes \bn} &= \bm{0}  &&\text{on $\Gamma$,}\\
\jump{\bN^T \bDHalf \bL} &= \bm{0}  &&\text{on $\Gamma$,}\\
\bN^T \bDHalf \bL  &= -\bg         &&\text{on $\Gamma_N$.}\\
\end{aligned} \right.
\end{equation}
As usual in an HDG context, the first equation in~\eqref{eq:elasticityStrongGlobal} is automatically satisfied due to the unique definition of the hybrid variable $\bhu$ on each face and the Dirichlet boundary condition $\bu_e = \bhu$ imposed in the local problems.

\subsection{Weak form of the local and global problems}
\label{sc:HDGweak}

Following the notation in~\cite{RS-SH:16}, the discrete functional spaces 
\begin{subequations}\label{eq:HDG-Spaces}
	\begin{align} 
	\Vh(\Omega) & := \left\{ v \in \eltwo(\Omega) : v \vert_{\Omega_e}\in \Pk(\Omega_e) \;\forall\Omega_e \, , \, e=1,\dotsc ,\numel \right\} , \label{eq:spaceScalarElem} \\
	\VhHat(S) & := \left\{ \hv \in \eltwo(S) : \hv\vert_{\Gamma_i}\in \Pk(\Gamma_i)
	\;\forall\Gamma_i\subset S\subseteq\Gamma\cup\partial\Omega \right\}, \label{eq:spaceScalarFace}
	\end{align}
\end{subequations}
are introduced, where $\mathcal{P}^{k}(\Omega_e)$ and $\mathcal{P}^{k}(\Gamma_i)$ are the spaces of polynomial functions of complete degree at most $k$ in $\Omega_e$ and on $\Gamma_i$ respectively. 
In addition, the classical internal products of vector functions in $\eltwo(\Omega_e)$ and $\eltwo(\Gamma_i)$
\begin{equation} \label{eq:innerScalar}
(\bm{p},\bm{q})_{\Omega_e} := \int_{\Omega_e} \bm{p} \cdot\bm{q} \ d\Omega  , \qquad \langle \hat{\bm{p}}, \hat{\bm{q}} \rangle_{\partial\Omega_e} := \sum_{\Gamma_i \subset \partial\Omega_e} \int_{\Gamma_i} \hat{\bm{p}} \cdot \hat{\bm{q}} \ d\Gamma 
\end{equation}
are considered. 

For each element $\Omega_e, \ e=1,\dotsc ,\numel$, the discrete weak formulation of~\eqref{eq:elasticityStrongLocal} reads as follows: given $\bu_D$ on $\Gamma_D$ and $\bhu $ on $\Gamma\cup\Gamma_N$,
find $(\bL^h_e ,\bu^h_e) \in [\Vh(\Omega_e)]^{\msd} \times [\Vh(\Omega_e)]^{\nsd}$ that satisfies 
\begin{subequations}\label{eq:HDGElasticityWeakLocalPre}
	\begin{align}
	&
	- (\bv,\bL^h_e)_{\Omega_e} + (\gradS^T \bDHalf \bv, \bu^h_e)_{\Omega_e} =   \langle \bN_e^T \bDHalf \bv , \bu_D\rangle_{\partial\Omega_e\cap\Gamma_D} + \langle \bN_e^T \bDHalf \bv , \bhu^h \rangle_{\partial\Omega_e\setminus\Gamma_D} , \label{eq:HDGElasticityWeakLocalLPre}
	\\
	&
	-(\gradS \bw, \bDHalf \bL^h_e)_{\Omega_e}  + \langle \bw ,\bN_e^T \widehat{\bDHalf \bL^h_e} \rangle_{\partial\Omega_e} =  (\bw,\bm{f})_{\Omega_e} , \label{eq:HDGElasticityWeakLocalUPre}
	\end{align}
\end{subequations}
for all $(\bv ,\bw) \in [\Vh(\Omega_e)]^{\msd} \times [\Vh(\Omega_e)]^{\nsd}$.

Integrating by parts Equation~\eqref{eq:HDGElasticityWeakLocalUPre} and introducing the following definition of the trace of the numerical stress featuring a stabilisation parameter $\btau_e$
\begin{equation} \label{eq:traceElasticity}
\bN_e^T \widehat{\bDHalf \bL^h_e} := 
\begin{cases}
\bN_e^T \bDHalf \bL^h_e + \btau_e (\bu^h_e - \bu_D) & \text{on $\partial\Omega_e\cap\Gamma_D$,} \\
\bN_e^T \bDHalf \bL^h_e + \btau_e (\bu^h_e - \bhu^h) & \text{elsewhere,}  
\end{cases}
\end{equation}
leads to the symmetric form of the discrete weak local problem: for $e=1,\dotsc ,\numel$, given $\bu_D$ on $\Gamma_D$ and $\bhu $ on $\Gamma\cup\Gamma_N$,
find $(\bL^h_e ,\bu^h_e) \in [\Vh(\Omega_e)]^{\msd} \times [\Vh(\Omega_e)]^{\nsd}$ that satisfies 
\begin{subequations}\label{eq:HDGElasticityWeakLocal}
	\begin{gather}
	- (\bv,\bL^h_e)_{\Omega_e} + (\gradS^T \bDHalf \bv, \bu^h_e)_{\Omega_e} =   \langle \bN_e^T \bDHalf \bv , \bu_D\rangle_{\partial\Omega_e\cap\Gamma_D} + \langle \bN_e^T \bDHalf \bv , \bhu^h \rangle_{\partial\Omega_e\setminus\Gamma_D} , \label{eq:HDGElasticityWeakLocalL}
	\\
  (\bw, \gradS^T \bDHalf \bL^h_e)_{\Omega_e}  {+} 
  \langle \bw , \btau_e \bu^h_e \rangle_{\partial\Omega_e} {=}  
  (\bw,\bm{f})_{\Omega_e}  + \langle \bw , \btau_e \bu_D\rangle_{\partial\Omega_e\cap\Gamma_D} {+} 
  \langle \bw , \btau_e \bhu^h \rangle_{\partial\Omega_e\setminus\Gamma_D} ,
  \label{eq:HDGElasticityWeakLocalU}
  \end{gather}
\end{subequations}
for all $(\bv ,\bw) \in [\Vh(\Omega_e)]^{\msd} \times [\Vh(\Omega_e)]^{\nsd}$.

Similarly, the discrete weak form of the global problem that accounts for the transmission conditions and the Neumann boundary condition is: find $\bhu^h\in[\VhHat(\Gamma\cup\Gamma_N)]^{\nsd}$ such that

\begin{multline} \label{eq:HDGElasticityWeakGlobal}
  \sum_{e=1}^{\numel}\Bigl\{
  \langle \bhw, \bN_e^T \bDHalf \bL^h_e \rangle_{\partial\Omega_e\setminus\Gamma_D}
  + \langle \bhw,\btau_e\, \bu^h_e \rangle_{\partial\Omega_e\setminus\Gamma_D} 
  - \langle \bhw,\btau_e\,\bhu^h \rangle_{\partial\Omega_e\setminus\Gamma_D}\Bigr\}  = \\
  -\sum_{e=1}^{\numel} \langle \bhw, \bm{g} \rangle_{\partial\Omega_e\cap\Gamma_N},
\end{multline}
for all $\bhw\in[\VhHat(\Gamma\cup\Gamma_N)]^{\nsd}$.

\subsection{Spatial discretisation}
\label{sc:HDGdiscretisation}

The discretisation of the weak form of the local problem given by Equation~\eqref{eq:HDGElasticityWeakLocal} using an isoparametric formulation for the primal and mixed variables leads to a linear system with the following structure
\begin{equation} \label{eq:localProblemSystem}
\begin{bmatrix}
\mat{A}_{LL} & \mat{A}_{Lu} \\
\mat{A}_{Lu}^T & \mat{A}_{uu}\\
\end{bmatrix}_e
\begin{Bmatrix}
\nodalLV_e \\
\nodaluV_e 
\end{Bmatrix} 
= 
\begin{Bmatrix}
\vect{f}_L\\
\vect{f}_u
\end{Bmatrix}_e 
+
\begin{bmatrix}
\mat{A}_{L\hat{u}}\\
\mat{A}_{u\hat{u}}
\end{bmatrix}_e
\nodaluhV_e,
\end{equation}
for $e=1,\dotsc ,\numel$. 

Similarly, using an isoparametric formulation for the hybrid variable produce the following system of equations
\begin{equation}\label{eq:globalProblemSystem}
\sum_{e=1}^{\numel}\Big\{
\begin{bmatrix} \mat{A}_{L \hat{u}}^T & \mat{A}_{u \hat{u}}^T \end{bmatrix}_e
\begin{Bmatrix} \nodalLV_e \\ \nodaluV_e \end{Bmatrix}
+
[\mat{A}_{\hat{u}\hat{u}}]_e \: \nodaluhV_e \Big\}
= 
\sum_{i=e}^{\numel} [\vect{f}_{\hat{u}}]_e.
\end{equation}

The expressions of the matrices and vectors appearing in Equation~\eqref{eq:localProblemSystem}-\eqref{eq:globalProblemSystem} are detailed in Appendix~\ref{sc:implementation}.

After replacing the solution of the local problem of Equation~\eqref{eq:localProblemSystem} in Equation~\eqref{eq:globalProblemSystem}, the global problem becomes
\begin{equation}\label{eq:globalProblemSystemFinal}
\mat{\widehat{K}}\vect{\hat{u}}=\vect{\hat{f}},
\end{equation}
with
\begin{subequations}\label{eq:Dglobal-matvect}
	\begin{gather}
	\mat{\widehat{K}} =
	\Assem_{e=1}^{\numel}
	\begin{bmatrix} \mat{A}_{L \hat{u}}^T & \mat{A}_{u \hat{u}}^T \end{bmatrix}_e
	\begin{bmatrix}
	\mat{A}_{LL} & \mat{A}_{Lu} \\
	\mat{A}_{Lu}^T & \mat{A}_{uu}\\
	\end{bmatrix}_e^{-1}
	\begin{bmatrix}
	\mat{A}_{L\hat{u}}\\
	\mat{A}_{u\hat{u}}
	\end{bmatrix}_e
	+
	[\mat{A}_{\hat{u}\hat{u}}]_e \\
	\intertext{and}
	\vect{\hat{f}}=
	\Assem_{e=1}^{\numel} [\vect{f}_{\hat{u}}]_e
	-
	\begin{bmatrix} \mat{A}_{L \hat{u}}^T & \mat{A}_{u \hat{u}}^T \end{bmatrix}_e
	\begin{bmatrix}
	\mat{A}_{LL} & \mat{A}_{Lu} \\
	\mat{A}_{Lu}^T & \mat{A}_{uu}\\
	\end{bmatrix}_e^{-1}
	\begin{Bmatrix} \vect{f}_L \\ \vect{f}_u \end{Bmatrix}_e .
	\end{gather}
\end{subequations}

\subsection{A remark on the $\eltwo$ convergence rates for the primal and mixed variables}
\label{sc:HDGconvergence}

Differently from the classical results for HDG~\cite{Jay-CGL:09,Cockburn-CDG:08,Nguyen-NPC:09,Nguyen-NPC:09b,Nguyen-NPC:10,Nguyen-NPC:11}, the best convergence rates proved by Cockburn and co-workers for the linear elasticity equation~\cite{MR3340089} only achieve a convergence of order $k$ for the gradient of the displacement field.
The convergence rate of both the strain and stress tensors achieves order $k+1/2$ but remains sub--optimal with respect to the one of the displacement field (order $k+1$).
This issue vanishes when moving to high--order approximations in which the optimal convergence of the gradient of the displacement field is retrieved.
Nevertheless, the aforementioned limitation represents a major drawback for the application of the classical HDG formulation using polynomials of degree less than 3.

The formulation based on Voigt notation discussed in this article outperforms the convergence rates proved in~\cite{MR3340089} by always achieving order $k+1$ for all the variables (cf. Section~\ref{sc:examples}).
The possibility of deriving a sharper \emph{a priori} bound for the mixed variable exploiting the rationale introduced by the Voigt notation will be investigated in a future work.
In next Section, the optimal numerical convergence of the mixed variable is exploited to construct a post--processed displacement field which super--converges with order $k+2$.


\section{Super--convergent post--process of the displacement field}
\label{sc:superconvergence}

As previously mentioned, a known feature of the HDG method is the possibility to exploit the accuracy granted by the convergence of order $k+1$ of the mixed variable (i.e. the stress tensor) to perform a local post--process of the primal variable and construct element--by--element a displacement field $\bu^\star$ super--converging with order $k+2$.
Nevertheless, for the linear elastic problem under analysis the classical approach in~\cite{soon2009hybridizable} shows some issues resulting in a loss of super--convergence of the post--processed solution for low--order approximations.
Following~\cite{RS-SH:16}, in this section a novel post--process procedure is discussed and the super--convergence of $\bu^\star$ is retrieved.

Introduce the space $\Vh_\star(\Omega)$ of the polynomials of complete degree at most $k+1$ on each element $\Omega_e$:
\begin{equation} \label{eq:Vstar}
\Vh_\star(\Omega) := \left\{ v \in \eltwo(\Omega) : v \vert_{\Omega_e}\in \mathcal{P}^{k+1}(\Omega_e) \;\forall\Omega_e \, , \, e=1,\dotsc ,\numel \right\} .
\end{equation}
For each element $\Omega_e, \ e=1,\dotsc,\numel$, consider the definition of the mixed variable in~\eqref{eq:elasticityStrongLocal}:
\begin{equation} \label{eq:mixedVar}
\bL_e + \bDHalf \gradS \bu_e = \bm{0} .
\end{equation}
The post--processed solution $\bu^\star$ is sought in the richer space $\left[ \Vh_\star(\Omega) \right]^{\nsd}$ and fulfils the following element--by--element problem:
\begin{equation} \label{eq:postVoigt}
\left\{\begin{aligned}
\gradS^T \bDHalf \gradS \bu_e^\star  &= - \gradS^T \bL_e       &&\text{in $\Omega_e, \ e=1,\dotsc,\numel$,}\\
\bN^T \bDHalf \gradS \bu_e^\star &= - \bN^T \bL_e        &&\text{on $\partial\Omega_e$.}\\
\end{aligned}\right.
\end{equation}
\begin{remark}
	The solution of Equation~\eqref{eq:postVoigt} is not uniquely identified in $\left[ \Vh_\star(\Omega) \right]^{\nsd}$.
	More precisely, it is unique excluding \emph{rigid motions}, that is up to a family of functions $\bv^\star$ such that $\gradS \bv^\star = \bm{0}$.
	From a practical point of view, $\bu^\star$ is identified up to three (respectively, six) constants in two (respectively, three) dimensions. 
	Each constant is associated with one rigid motion, namely two translations and one rotation (respectively, three and three) in 2D (respectively, 3D).
\end{remark}

Consider the classical solvability constraint added in the HDG literature to close Equation~\eqref{eq:postVoigt}:
\begin{equation} \label{eq:postprocessMean}
\int_{\Omega_e} \bu^\star_e \ d\Omega = \int_{\Omega_e} \bu^h_e \ d\Omega .
\end{equation} 
It is straightforward to observe that condition~\eqref{eq:postprocessMean} removes the under--determination related to the translational modes.
Nonetheless, one additional constraint is required in 2D and three in 3D in order to remove the rotational modes of the element.
In~\cite{soon2009hybridizable}, Soon \etal \ decompose the post--processed solution in two components, the first one arising from the projection of the HDG solution onto the space of rigid motion displacements and the second one from the solution of~\eqref{eq:postVoigt} in the space of polynomials with no rigid motion.
This approach may be interpreted as constraining~\eqref{eq:postVoigt}-\eqref{eq:postprocessMean} with the following additional condition accounting for the rigid rotation of the displacement field with respect to the barycentre of the element:
\begin{equation} \label{eq:postprocess1}
\int_{\Omega_e} (\bx - \bx_b) \times \bu^\star_e \ d\Omega  =  \int_{\Omega_e} (\bx - \bx_b) \times \bu^h_e \  d\Omega ,
\end{equation}
where $\bx$ is the position vector and $\bx_b$ is the barycentre of the element $\Omega_e$ under analysis.
This post--process technique is inspired by the work of Stenberg~\cite{MR954768} on mixed finite elements and allows to retrieve the uniqueness of the post--processed solution but the super--convergence is lost for low--order approximations.

To remedy this issue without resorting to the extremely elegant, but rather complicated, framework of the $\bm{M}$-decomposition discussed in~\cite{arXiv2017BCGF}, a novel constraint which has not been previously considered in the literature is proposed to substitute~\eqref{eq:postprocess1}.
More precisely, a constraint on the mean value of the $\curl$ inside the element $\Omega_e$ is introduced
\begin{equation} \label{eq:postprocess2}
\int_{\Omega_e} \grad \times \bu^\star_e \ d\Omega  =  \int_{\Omega_e} \grad \times \bu^h_e \ d\Omega .
\end{equation}
By applying the Stokes' theorem to the right--hand side of~\eqref{eq:postprocess2}, an alternative formulation which exploits the hybrid variable $\bhu$ is obtained
\begin{equation} \label{eq:postprocess3}
\int_{\Omega_e} \grad \times \bu^\star_e \ d\Omega  =  \int_{\partial \Omega_e} \bhu^h \cdot \bt_e \ d\Gamma ,
\end{equation}
where $\bt_e$ is the tangential direction to the boundary $\partial\Omega_e$.
By exploiting the Voigt notation introduced in Section~\ref{sc:DivTheoVoigt},~\eqref{eq:postprocess2}-\eqref{eq:postprocess3} may be written as
\begin{equation} \label{eq:postprocessVoigt}
\int_{\Omega_e} \bR \bu^\star_e \ d\Omega = 
\begin{cases}
\displaystyle\int_{\Omega_e} \bR \bu^h_e \ d\Omega
&\text{according to~\eqref{eq:postprocess2},} \\
\displaystyle\int_{\partial\Omega_e} \bigl[ \bhu^h \bigr]^T \bT \ d\Gamma 
&\text{according to~\eqref{eq:postprocess3}.} 
\end{cases}
\end{equation}
By intuition, the second formulation in~\eqref{eq:postprocessVoigt} guarantees a higher accuracy owing to the convergence rate of order $k+1$ of the hybrid variable $\bhu^h$.
An alternative physical interpretation of conditions~\eqref{eq:postprocess2}-\eqref{eq:postprocess3} is given in~\cite{preprintVoigtStokes}, exploiting the definition of the vorticity of a fluid as the curl of its velocity field.
For additional details on the post--process procedures inspired by the velocity--pressure--vorticity formulation of the Stokes equation, the interested reader is referred to~\cite{Nguyen-CNP:10, MR2772094}.

\section{Numerical examples}
\label{sc:examples}

\subsection{Optimal order of convergence}
\label{sc:convergence}

This section considers two examples, in two and three dimensions, with known analytical solution to test the optimal convergence properties of the error of the primal and mixed variables, $\bu$ and $\bL$ respectively, measured in the $\eltwo(\Omega)$ norm and for different types of elements.

\subsubsection{Two dimensional example}
\label{sc:convergence2D}

The first example considers the model problem of Equation~\eqref{eq:elasticity} in the domain $\Omega = [0,1]^2$. The external load is selected so that the analytical solution is 
\begin{equation}
\bu(\bx) = \frac{1}{100} \Big( x_2 \sin(\pi x_1), \quad x_1^3 + \cos(\pi x_2) \Big).
\end{equation}
Neumann boundary conditions, corresponding to the analytical normal stress, are imposed on $\Gamma_N = \{(x_1,x_2) \in \mathbb{R}^2 \; | \; x_2=0\}$ and Dirichlet boundary conditions, corresponding to the analytical solution, are imposed on $\Gamma_D = \partial \Omega \setminus \Gamma_N$. The Young's modulus is taken as $E=1$ and the Poisson's ratio is $\nu = 0.25$.

\begin{figure}[!tb]
	\centering
	\subfigure[Quadrilateral mesh 1]{\includegraphics[width=0.24\textwidth]{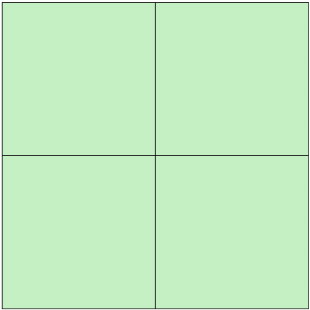}}
	\subfigure[Quadrilateral mesh 2]{\includegraphics[width=0.24\textwidth]{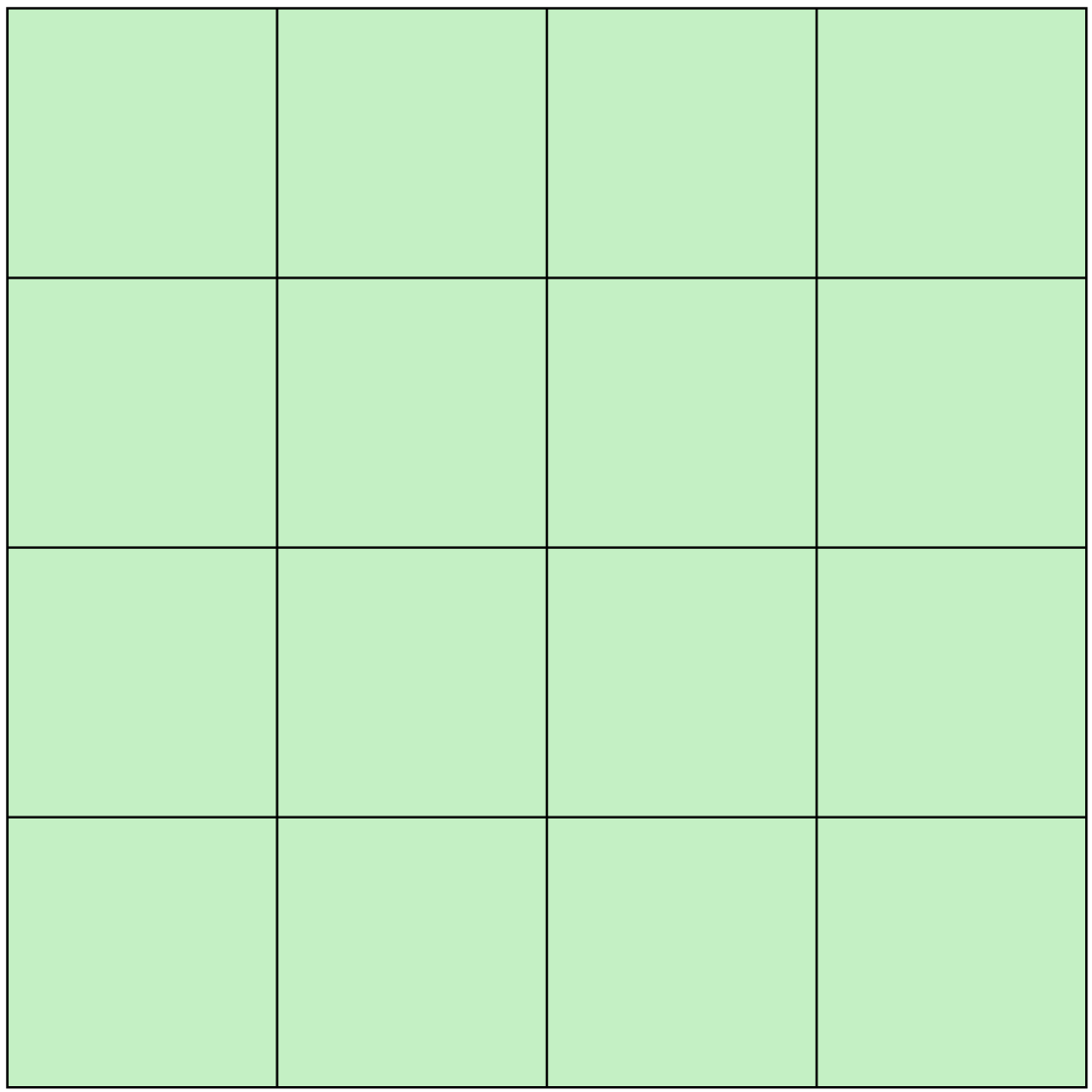}}
	\subfigure[Triangular mesh 1]{\includegraphics[width=0.24\textwidth]{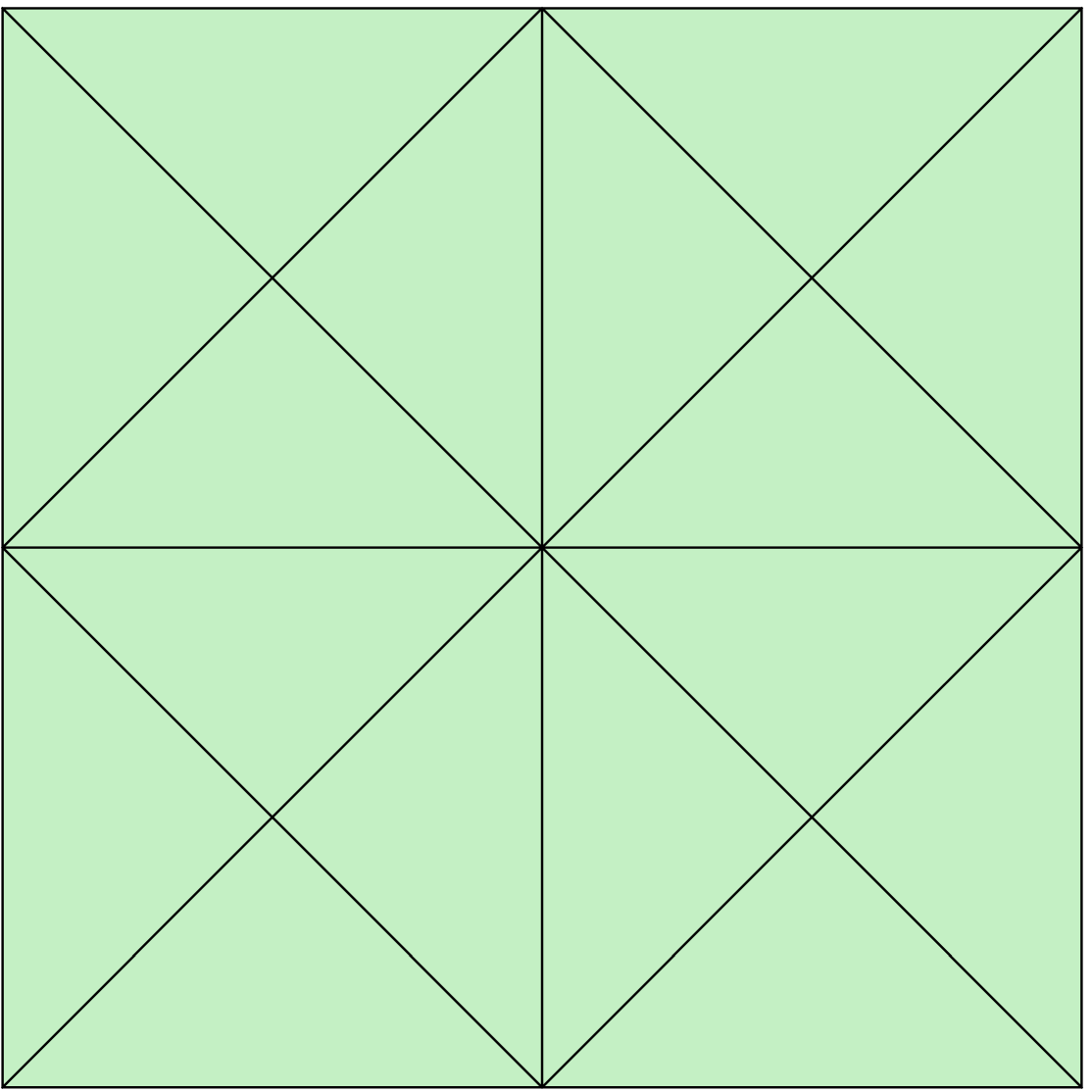}}
	\subfigure[Triangular mesh 2]{\includegraphics[width=0.24\textwidth]{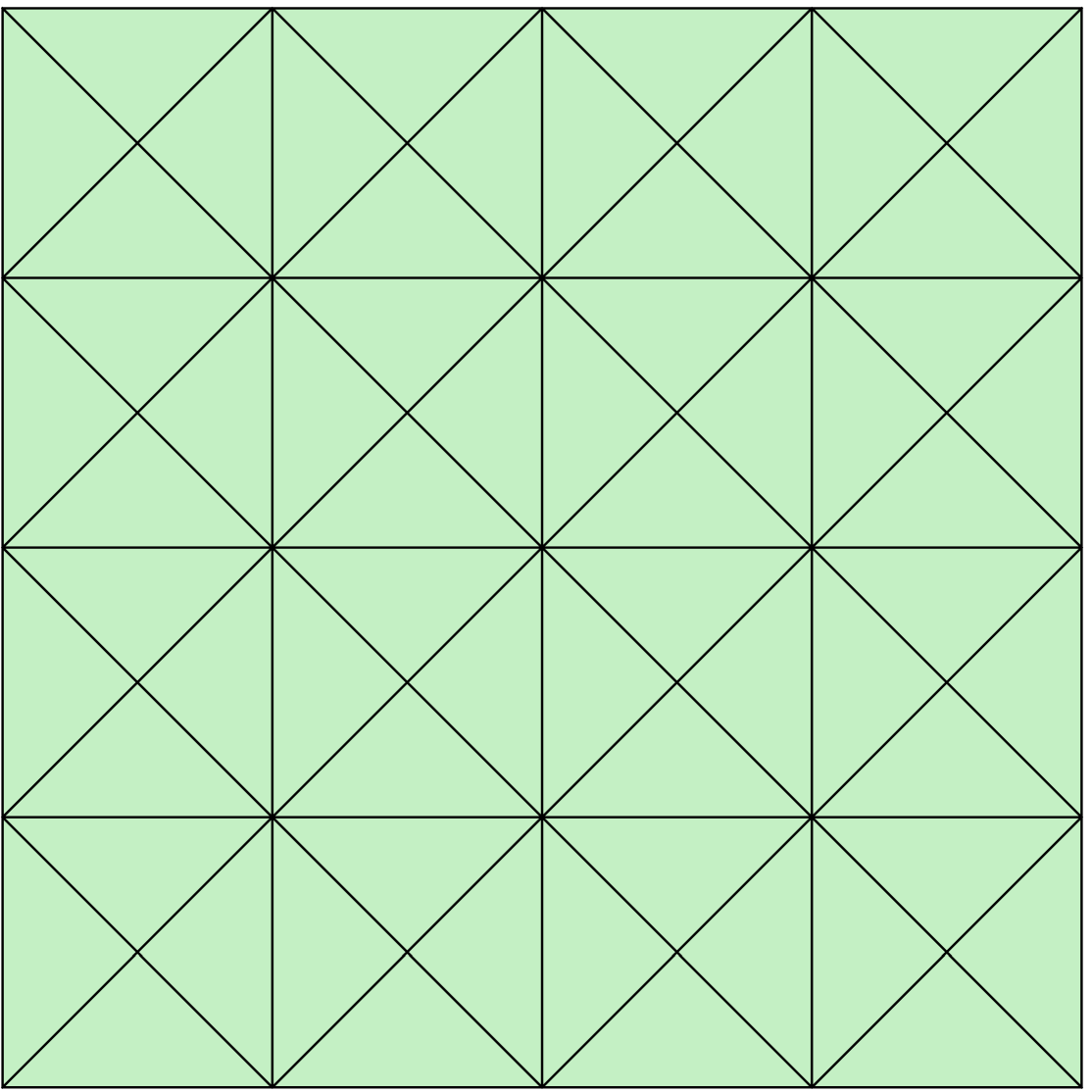}}
	\caption{Two dimensional meshes of $\Omega=[0,1]^2$ for the mesh convergence study.}
	\label{fig:2Dmeshes}
\end{figure}

Uniform meshes of quadrilateral and triangular elements are considered to perform an $h$--convergence study. The first two quadrilateral and triangular meshes are shown in Figure~\ref{fig:2Dmeshes}.

\begin{figure}[!tb]
	\centering
	\subfigure[$u_1$]{\includegraphics[width=0.32\textwidth]{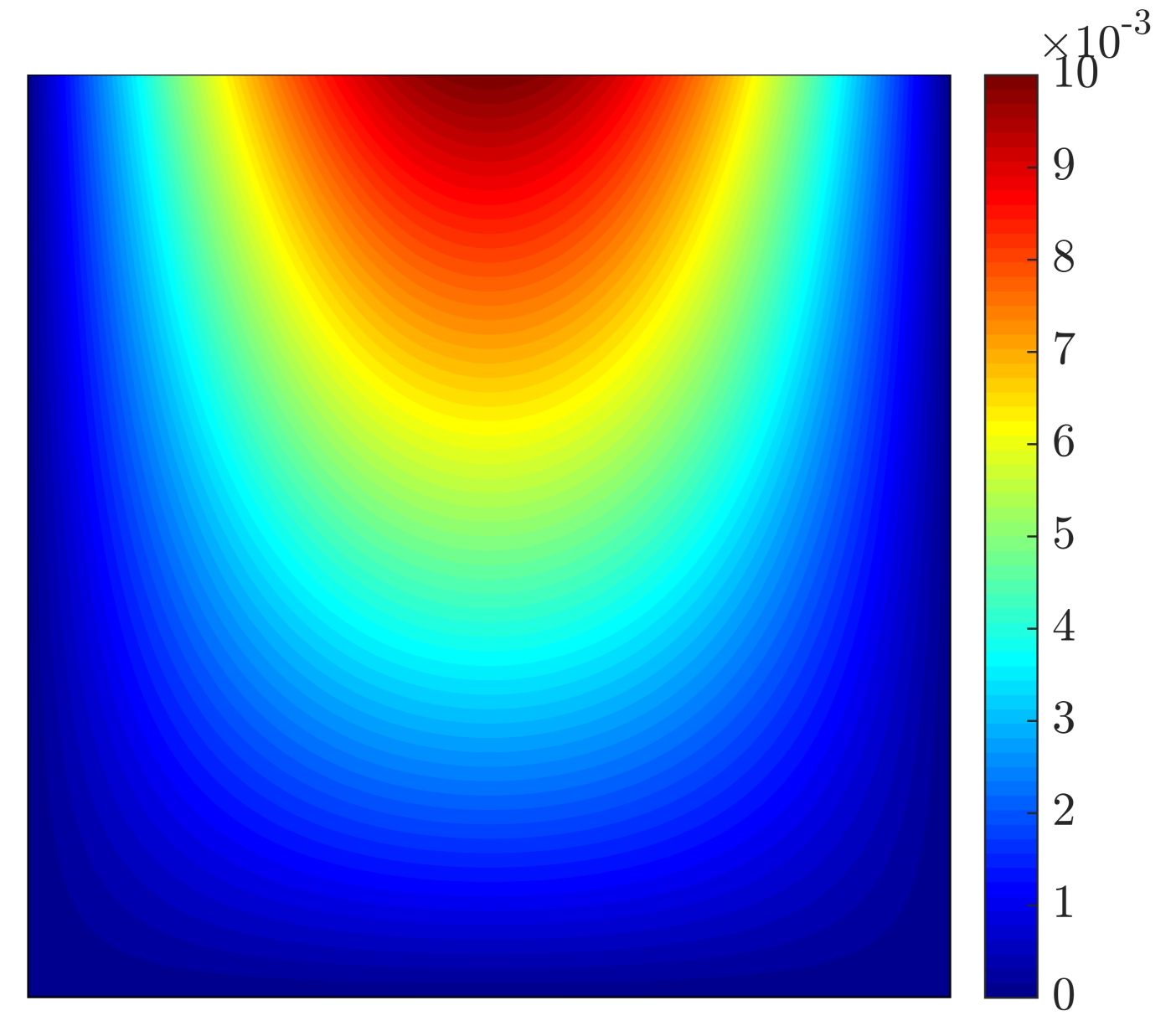}}
	\subfigure[$u_2$]{\includegraphics[width=0.32\textwidth]{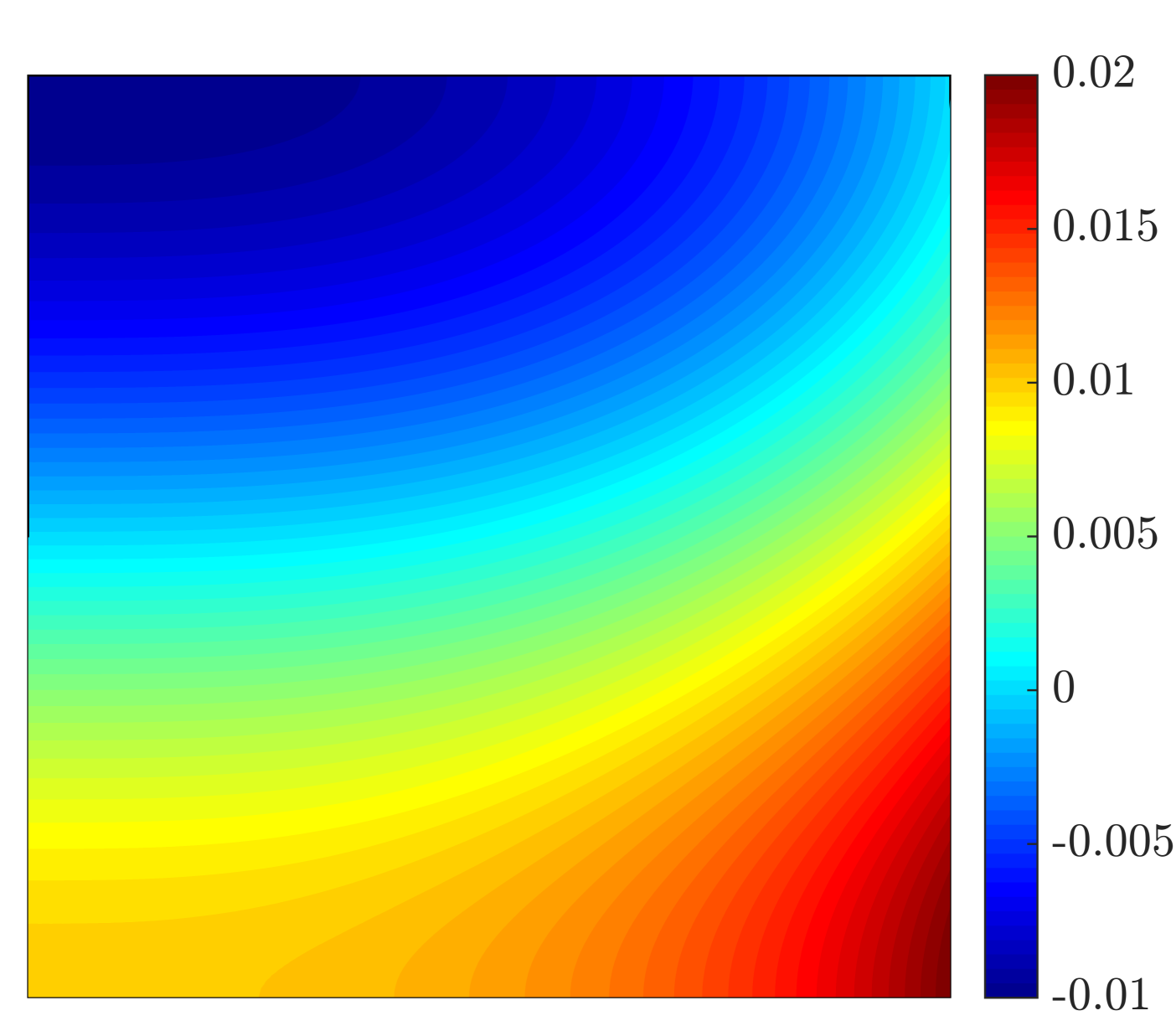}}
	\subfigure[$\sigma_{\texttt{VM}}$]{\includegraphics[width=0.32\textwidth]{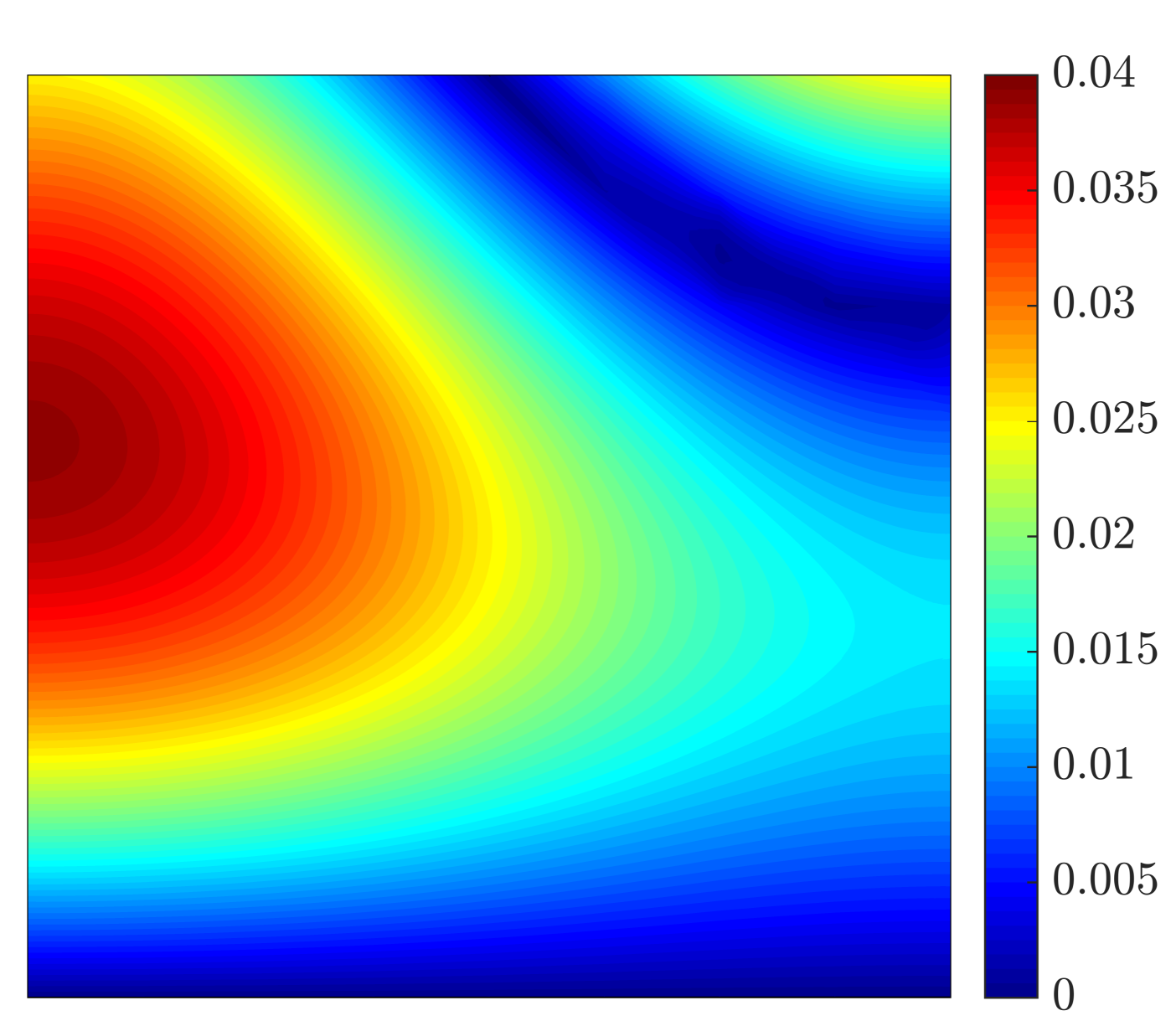}}
	\caption{Two dimensional problem: HDG approximation of the displacement field and the Von Mises stress using the third triangular mesh and $k=2$.}
	\label{fig:2Dsol}
\end{figure}

The displacement field and the Von Mises stress computed on the third triangular mesh and using a quadratic degree of approximation are depicted in Figure~\ref{fig:2Dsol}.

\begin{figure}[!tb]
	\centering
	\subfigure[Quadrilaterals]{\includegraphics[width=0.4\textwidth]{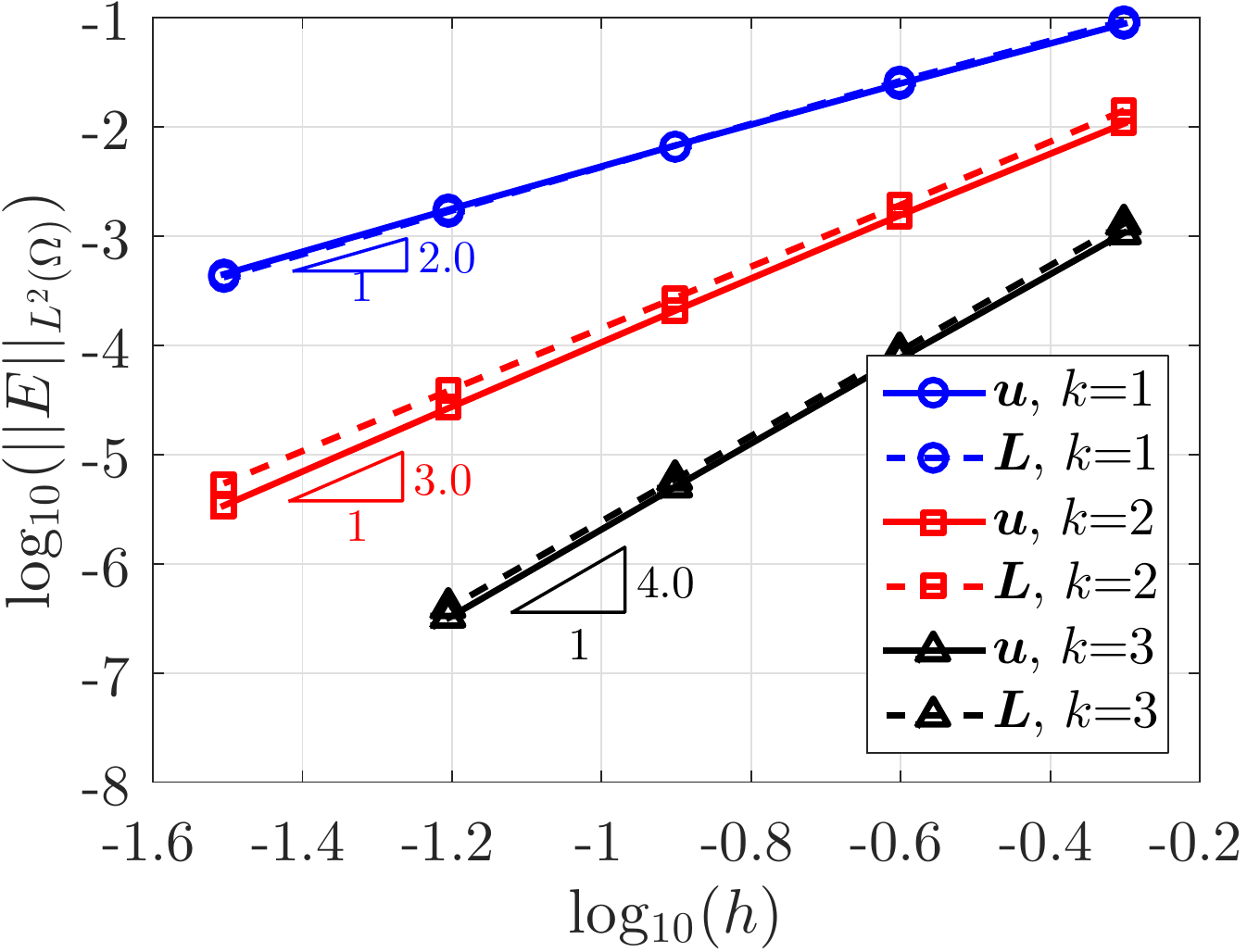}}
	\subfigure[Triangles]{\includegraphics[width=0.4\textwidth]{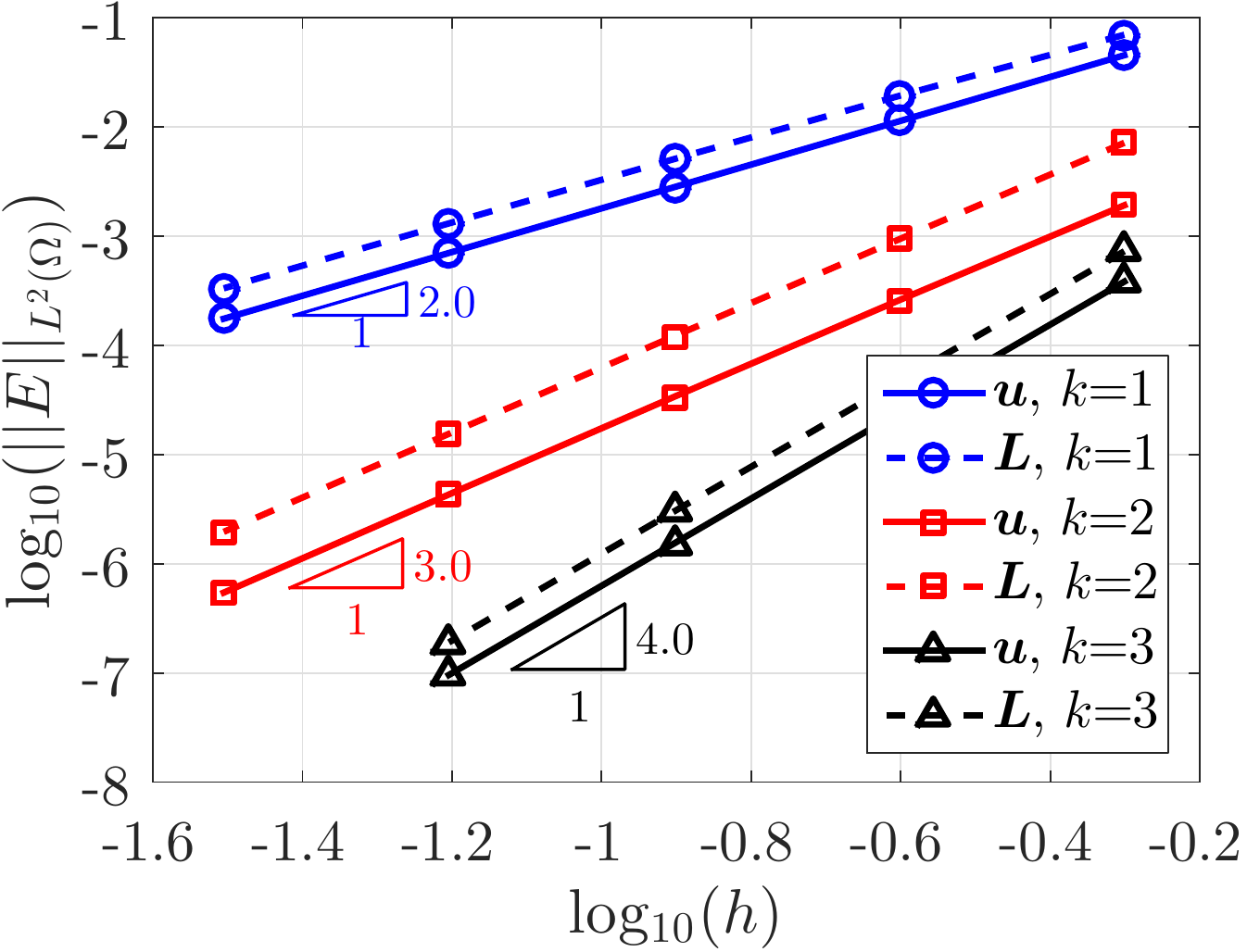}}	
	\caption{Two dimensional problem: $h$--convergence of the error of the primal and mixed variables, $\bu$ and $\bL$ in the $\eltwo(\Omega)$ norm for quadrilateral and triangular meshes with different orders of approximation.}
	\label{fig:hConv2D}
\end{figure}

The convergence of the error of the primal and mixed variables $\bu$ and $\bL$, measured in the $\eltwo(\Omega)$ norm, as a function of the characteristic element size $h$ is represented in Figure~\ref{fig:hConv2D} for both quadrilateral and triangular elements and for a degree of approximation ranging from $k=1$ up to $k=3$. 
It can be observed that the optimal rate of convergence $h^{k+1}$ is obtained for all the element types and degrees of approximation considered. 
It is worth noting that for the same characteristic element size, the triangular meshes have four times more internal faces than the quadrilateral mesh with the same element size. Therefore, despite the results in Figure~\ref{fig:hConv2D} indicate that for the same element size triangular elements provide more accurate results, when a comparison in terms of the number of degrees of freedom is performed both elements provide similar accuracy.

\subsubsection{Three dimensional example}
\label{sc:convergence3D}

The next example considers the model problem of Equation~\eqref{eq:elasticity} in the domain $\Omega = [0,1]^3$. The external load is selected so that the analytical solution is 
\begin{equation}
  \bu(\bx) = \frac{1}{100} 
  \begin{Bmatrix} x_1\sin(2\pi x_2) + x_2\cos(2\pi x_3) \\[1ex] 
                           x_2\sin(2\pi x_3) + x_3\cos(2\pi x_1) \\[1ex]
                           x_3\sin(2\pi x_1) + x_1\cos(2\pi x_2) \end{Bmatrix}.
\end{equation}
Neumann boundary conditions, corresponding to the analytical normal stress, are imposed on $\Gamma_N = \{(x_1,x_2,x_3) \in \mathbb{R}^3 \; | \; x_3=0\}$ and Dirichlet boundary conditions, corresponding to the analytical solution, are imposed on $\Gamma_D = \partial \Omega \setminus \Gamma_N$. The Young's modulus is taken as $E=1$ and the Poisson's ratio is $\nu = 0.25$.

\begin{figure}[!tb]
	\centering
	\subfigure[Hexahedral mesh 3]{\includegraphics[width=0.24\textwidth]{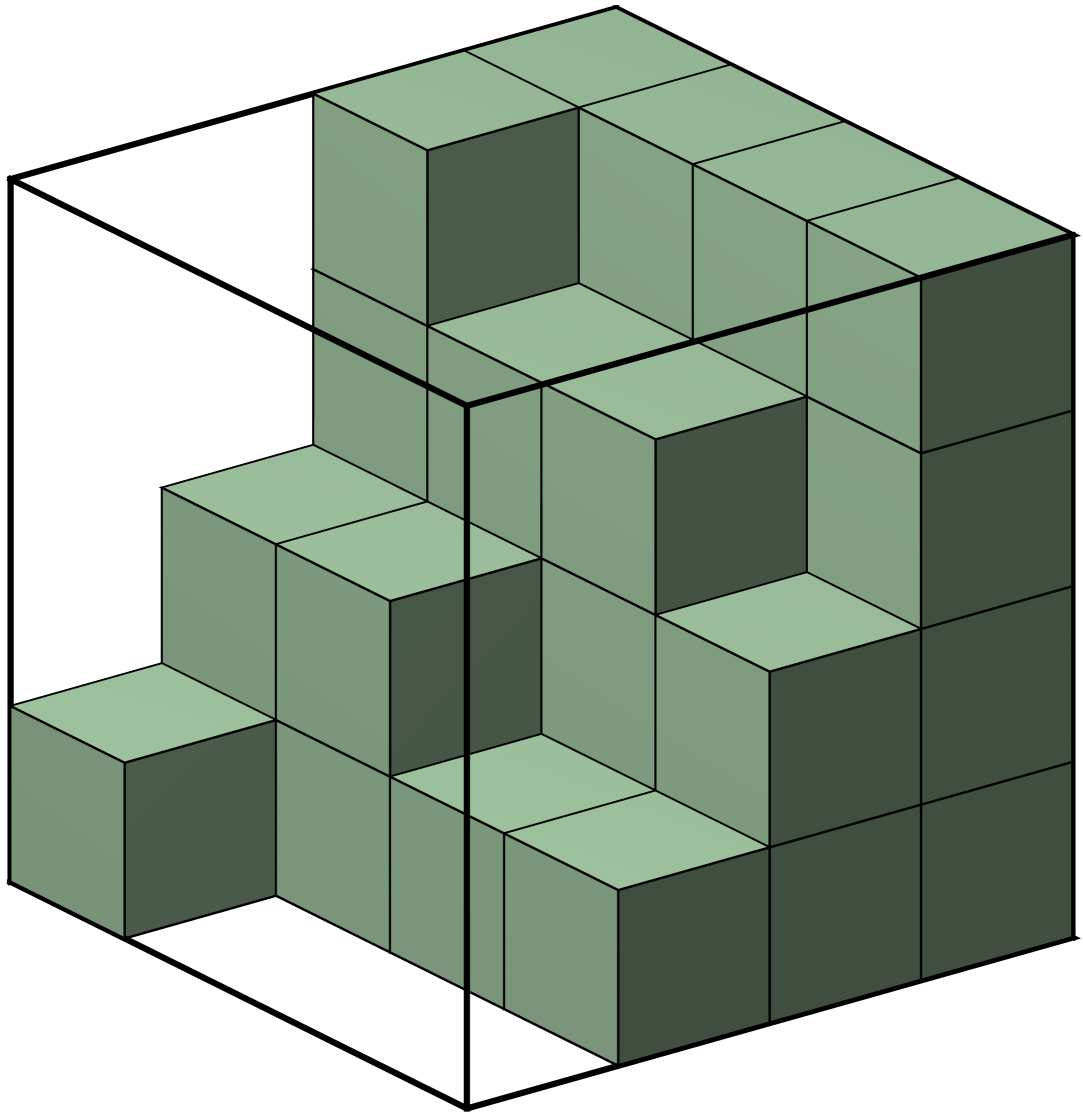}}
	\subfigure[Tetrahedral mesh 3]{\includegraphics[width=0.24\textwidth]{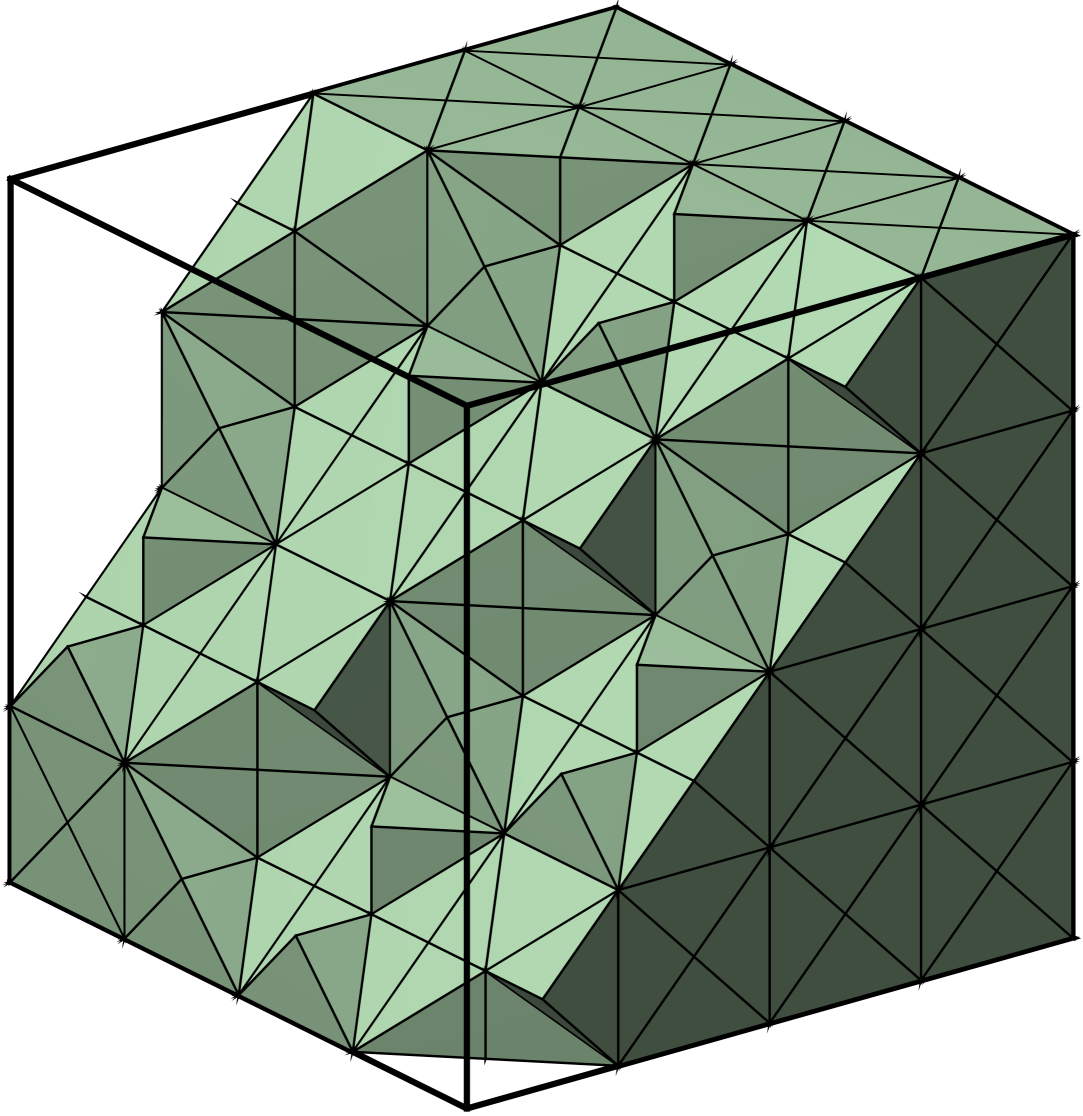}}		
	\subfigure[Prismatic mesh 3]{\includegraphics[width=0.24\textwidth]{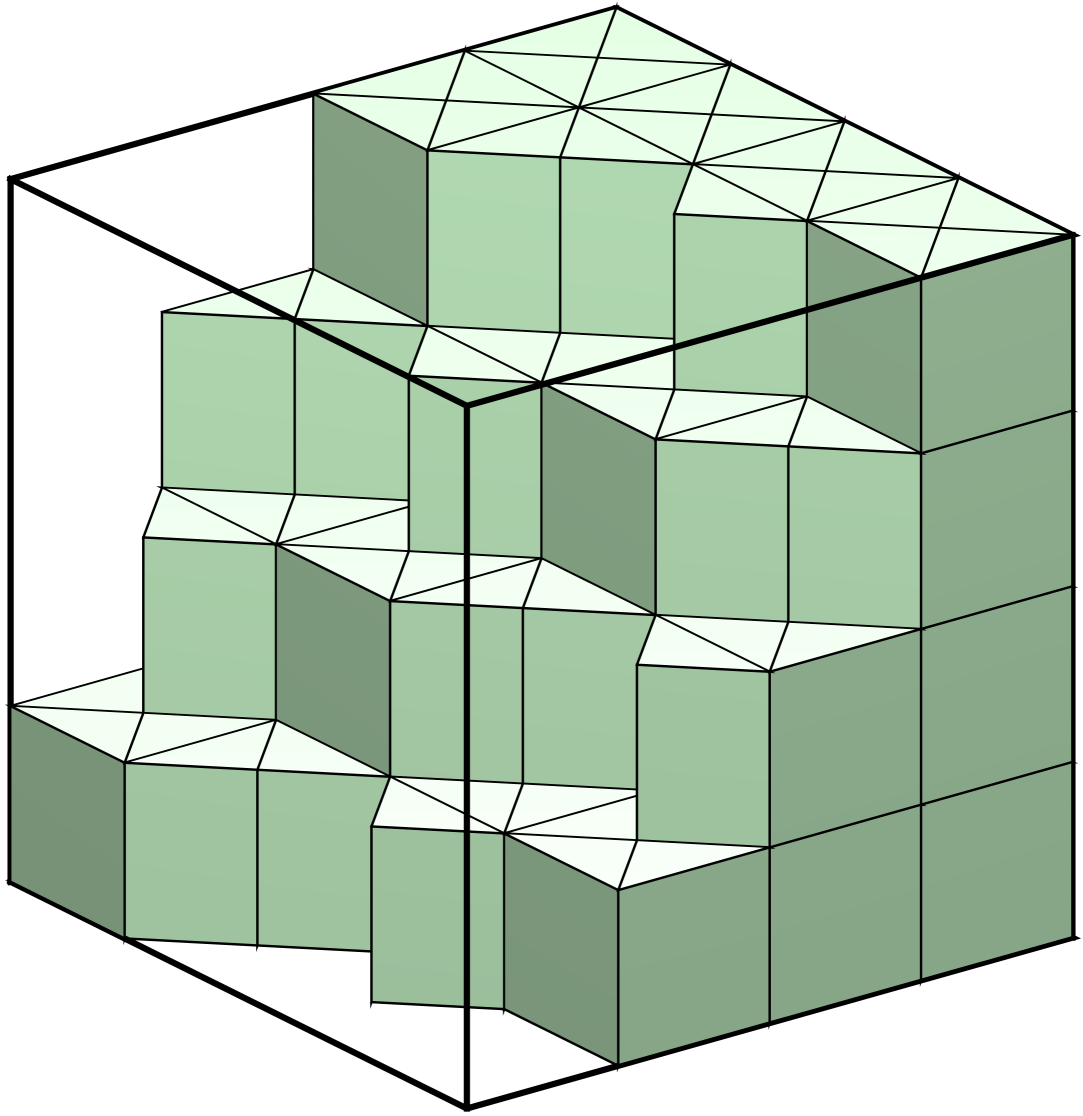}}
	\subfigure[Pyramidal mesh 3]{\includegraphics[width=0.24\textwidth]{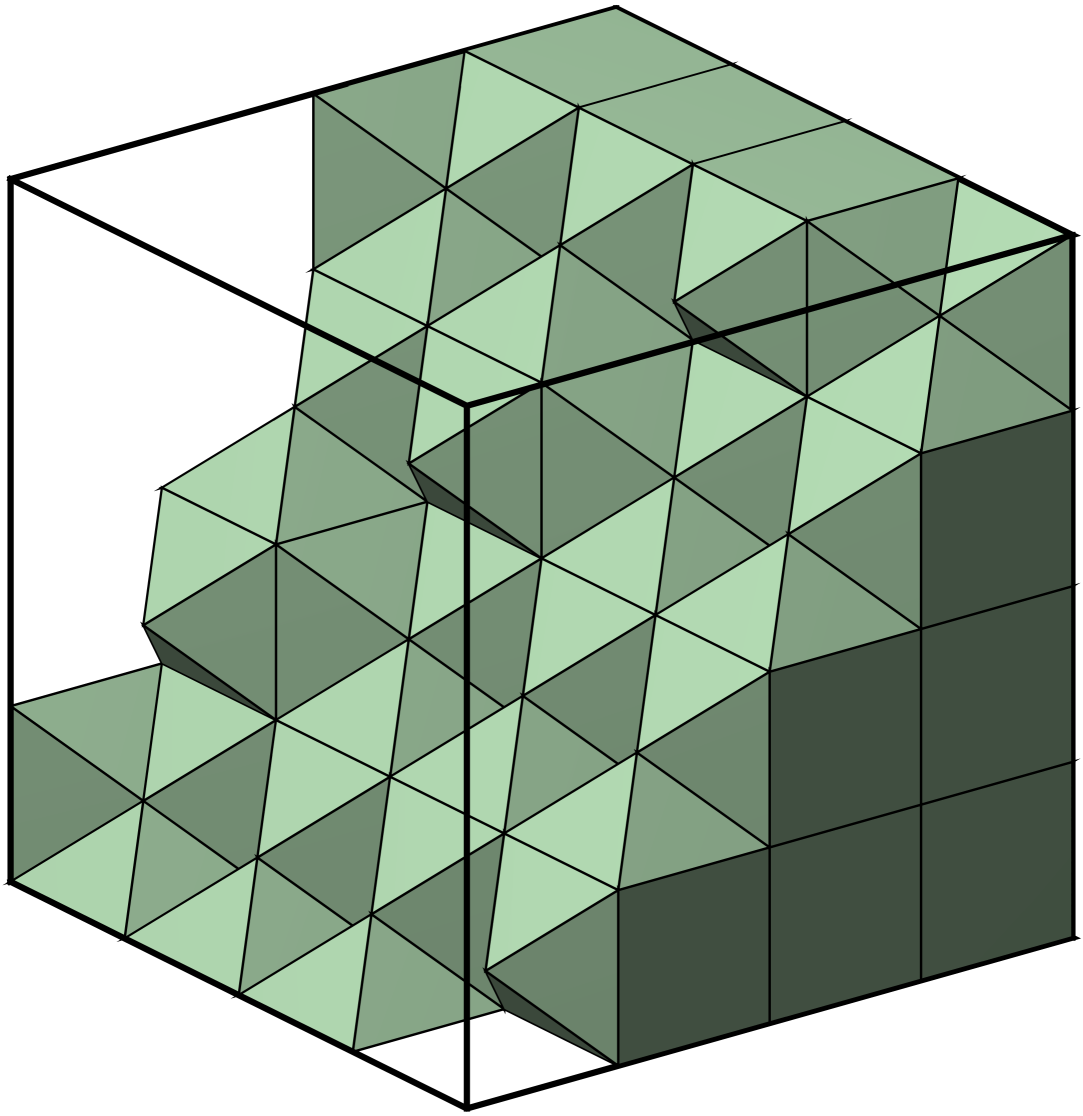}}
	\caption{Three dimensional meshes of $\Omega=[0,1]^3$ for the mesh convergence study.}
	\label{fig:3Dmeshes}
\end{figure}

Uniform meshes of hexahedral, tetrahedral, prismatic and pyramidal elements are considered to perform an $h$--convergence study. A cut through the meshes of the domain for the third level of refinement considered is represented in Figure~\ref{fig:3Dmeshes} for all the element types.

\begin{figure}[!tb]
	\centering
	\subfigure[$u_1$]{\includegraphics[width=0.24\textwidth]{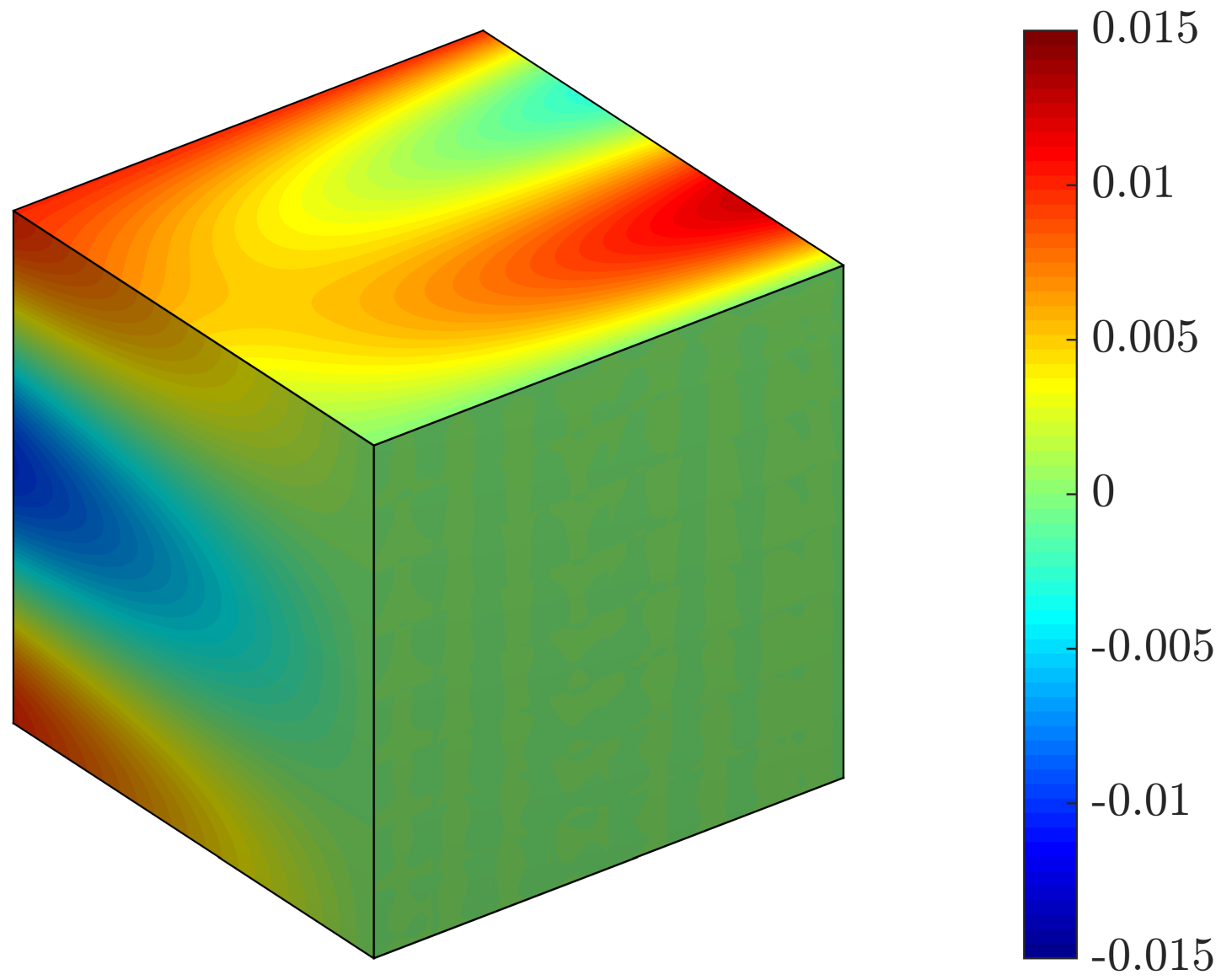}}
	\subfigure[$u_2$]{\includegraphics[width=0.24\textwidth]{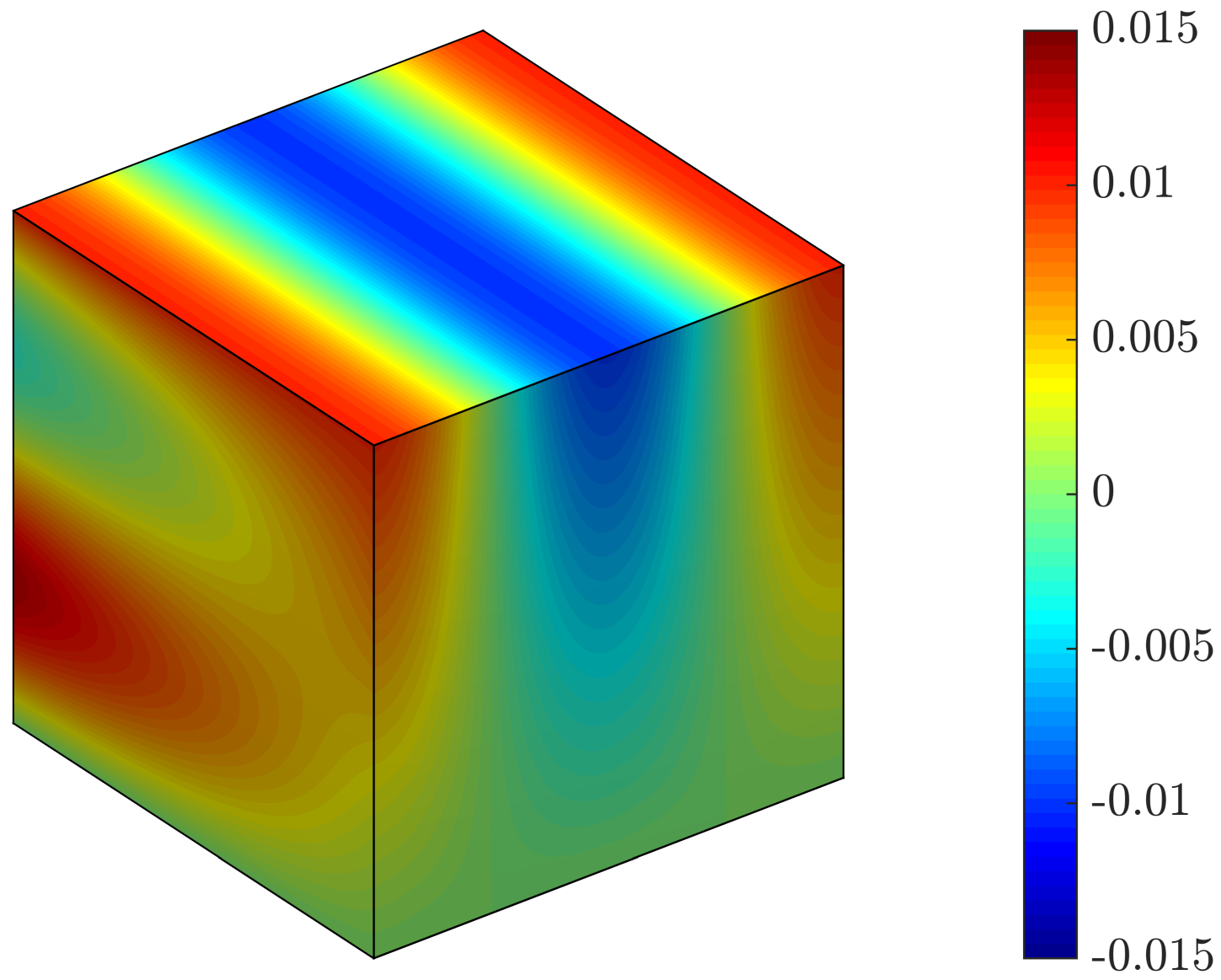}}
	\subfigure[$u_3$]{\includegraphics[width=0.24\textwidth]{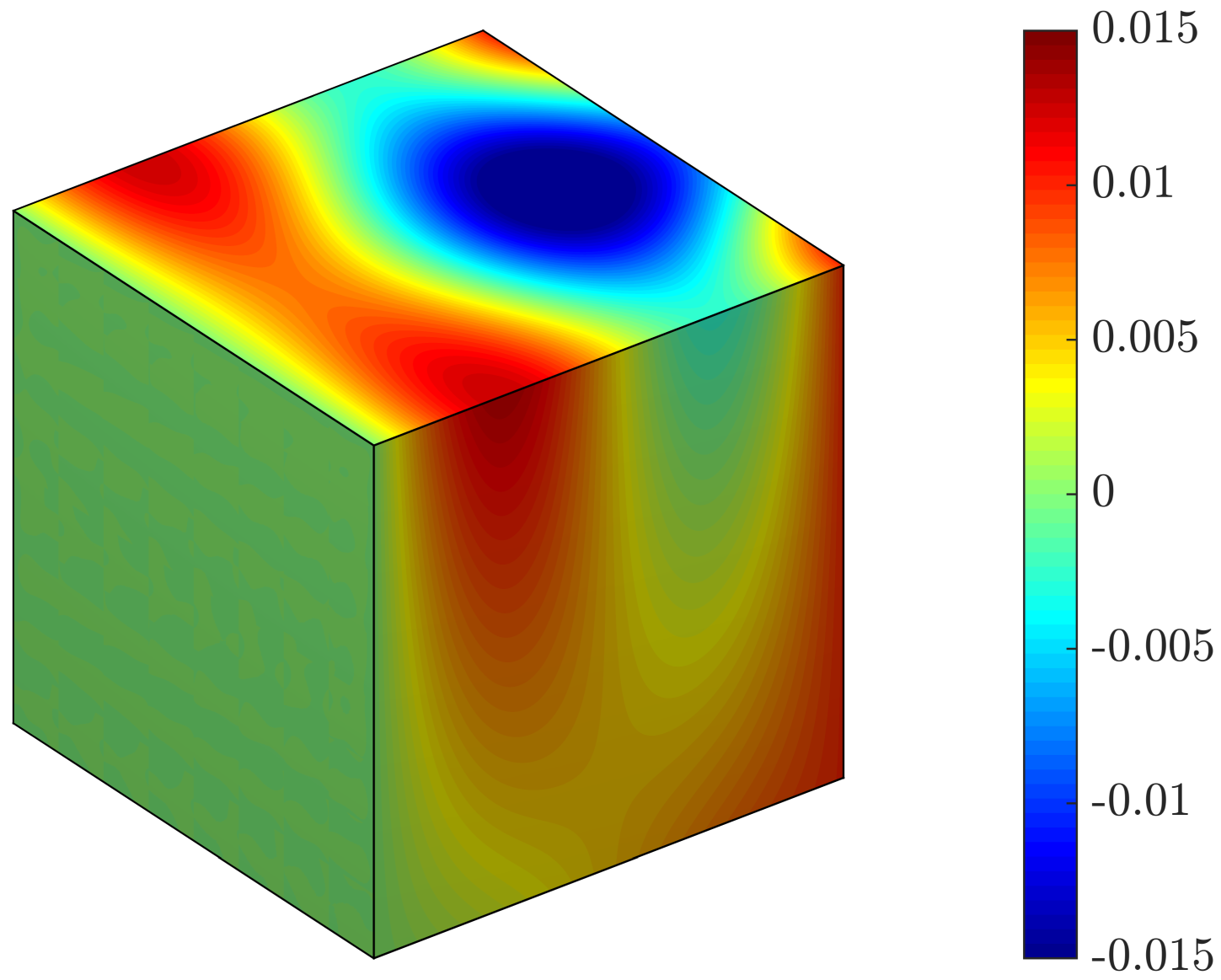}}
	\subfigure[$\sigma_{\texttt{VM}}$]{\includegraphics[width=0.24\textwidth]{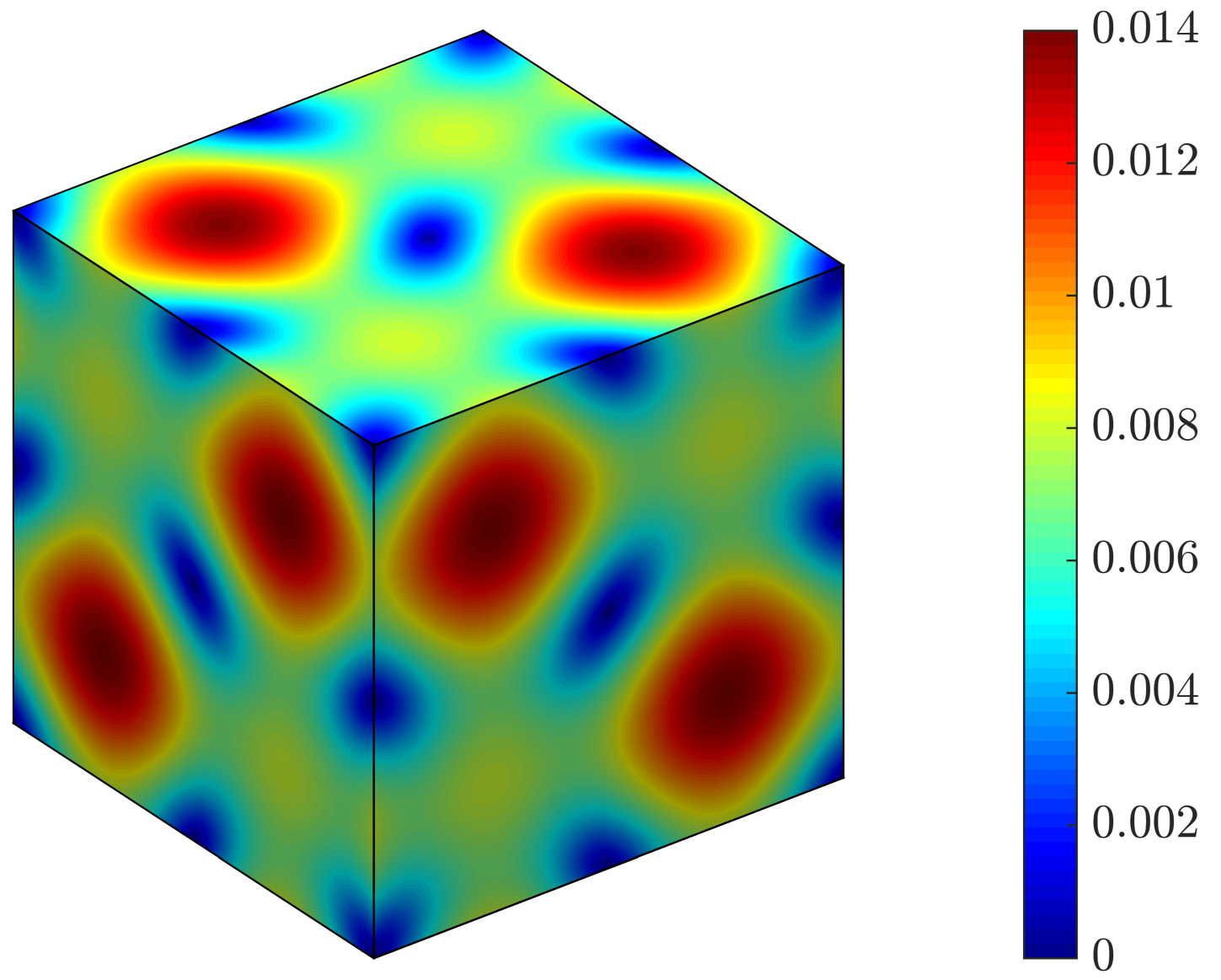}}
	\caption{Three dimensional problem: HDG approximation of the displacement field and the Von Mises stress using the fourth hexahedral mesh and $k=3$.}
	\label{fig:3Dsol}
\end{figure}

The displacement field and the Von Mises stress computed on the fourth hexahedral mesh and using a cubic degree of approximation are depicted in Figure~\ref{fig:3Dsol}.

\begin{figure}[!tb]
	\centering
	\subfigure[Hexahedrons]{\includegraphics[width=0.4\textwidth]{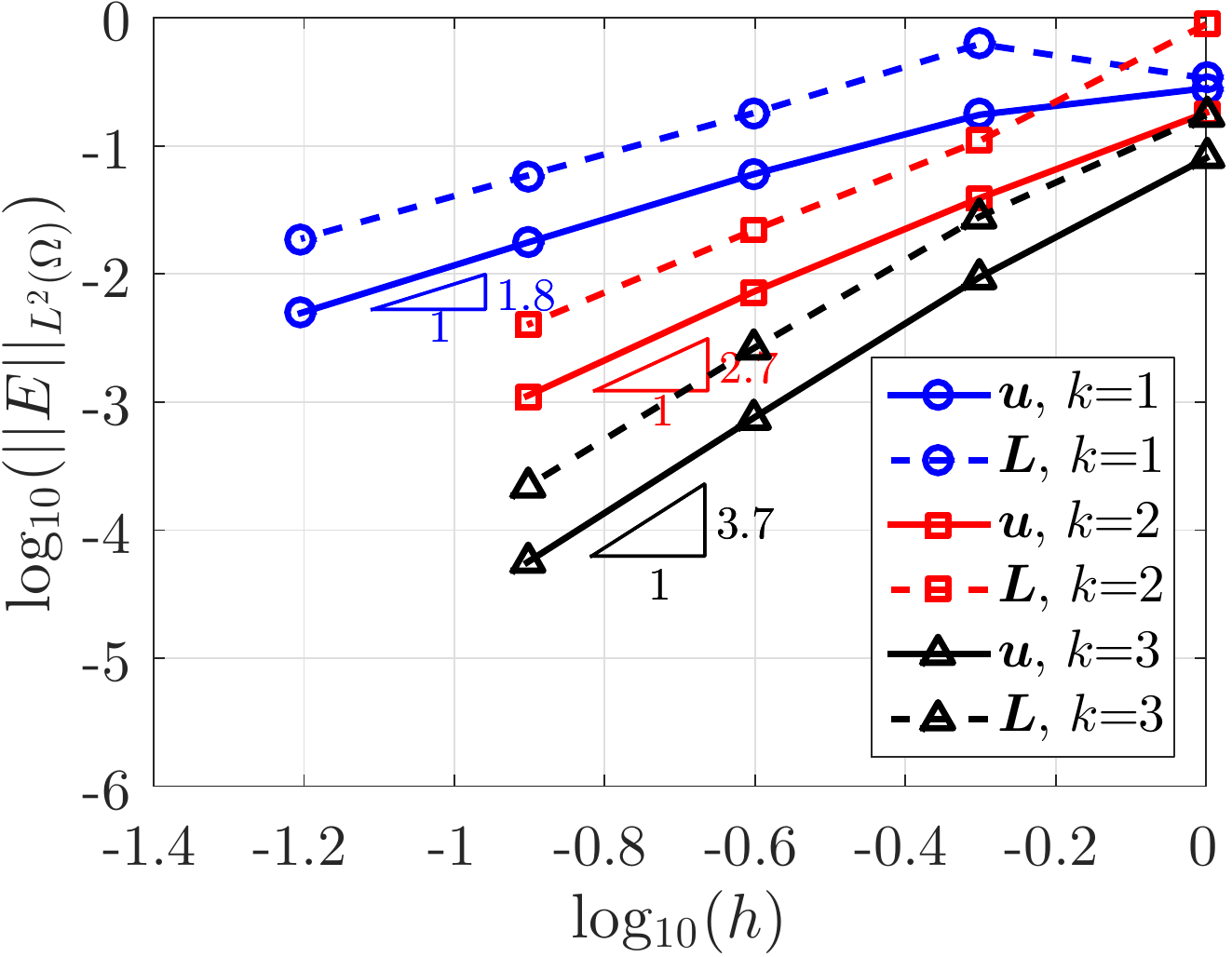}}
	\subfigure[Tetrahedrons]{\includegraphics[width=0.4\textwidth]{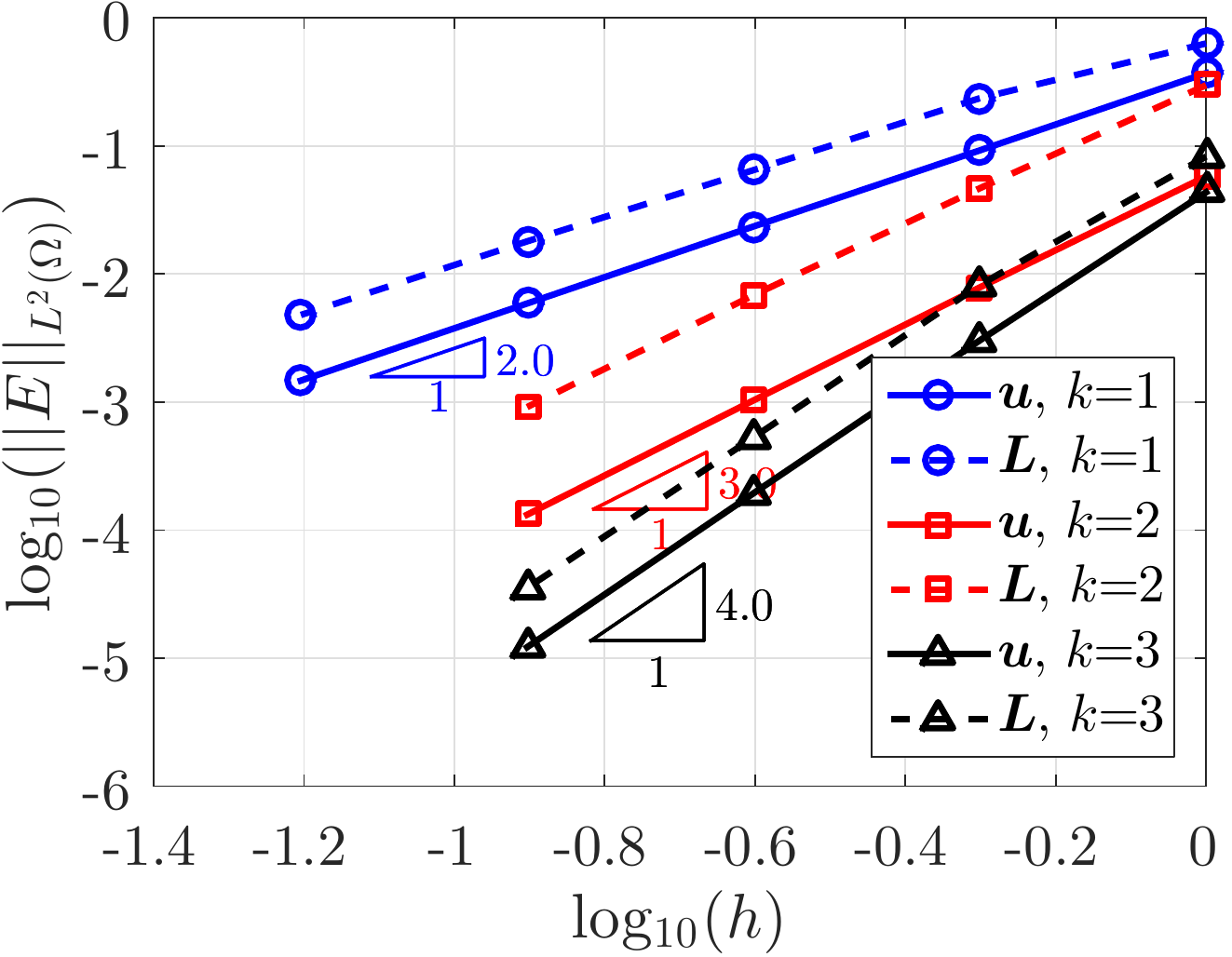}}
	\subfigure[Prisms]{\includegraphics[width=0.4\textwidth]{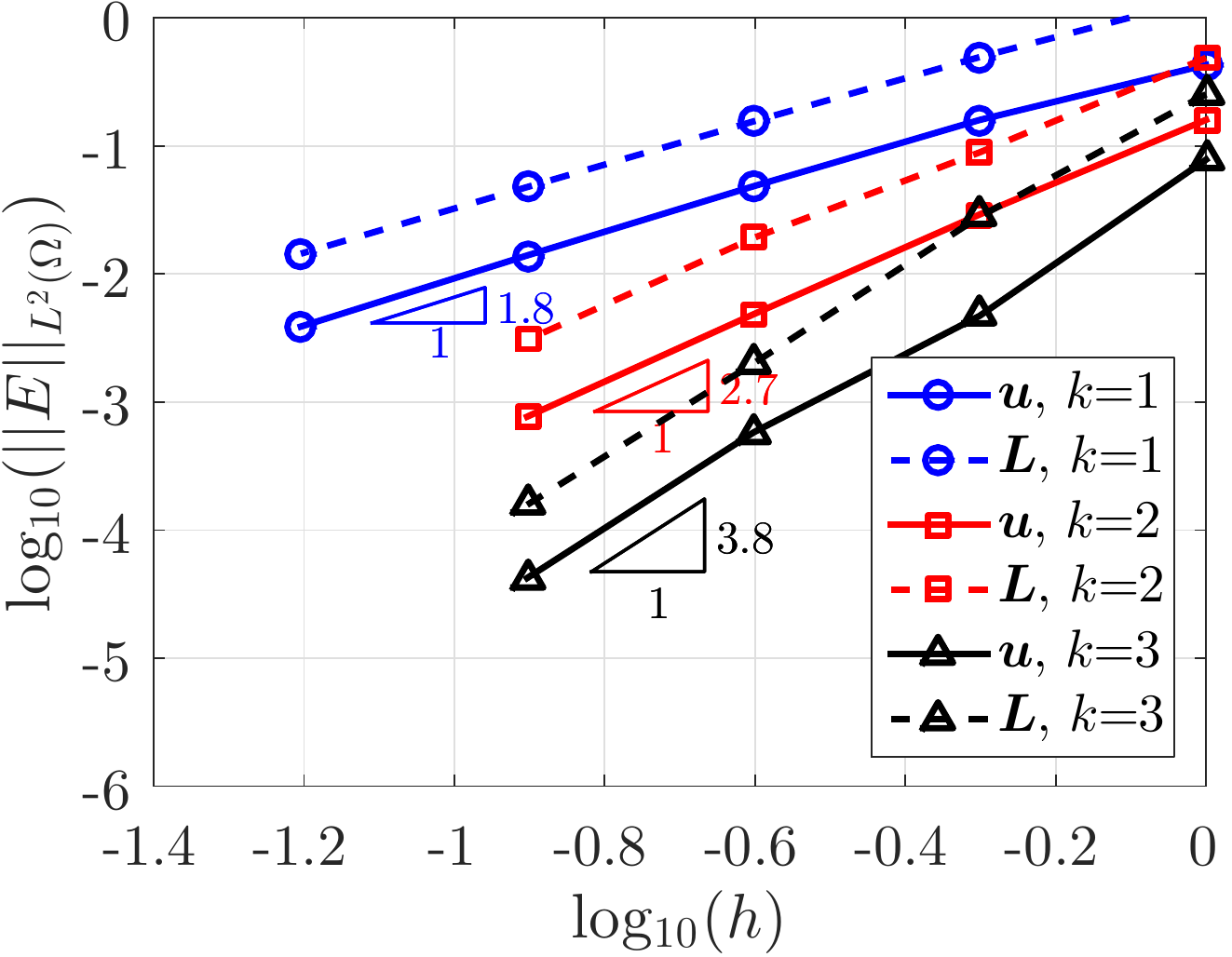}}
	\subfigure[Pyramids]{\includegraphics[width=0.4\textwidth]{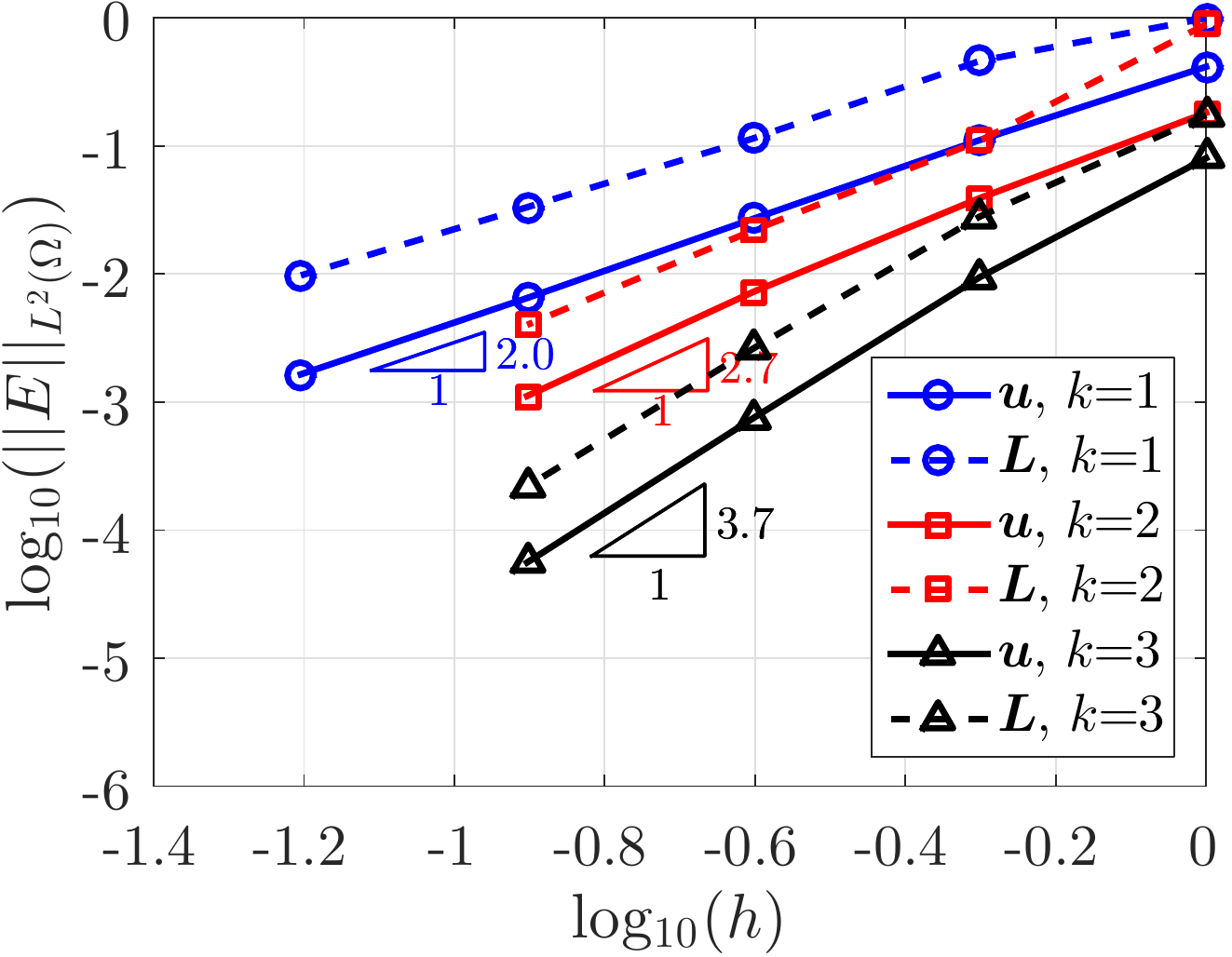}}
	\caption{Three dimensional problem: $h$--convergence of the error of the primal and mixed variables, $\bu$ and $\bL$ in the $\eltwo(\Omega)$ norm for hexahedral, tetrahedral, prismatic and pyramidal  meshes with different orders of approximation.}
	\label{fig:hConv3D}
\end{figure}

Analogously to the previous example, the convergence of the error of the primal and mixed variables $\bu$ and $\bL$, measured in the $\eltwo(\Omega)$ norm, as a function of the characteristic element size $h$ is represented in Figure~\ref{fig:hConv3D} for all the element types and for a degree of approximation ranging from $k=1$ up to $k=3$. 
It can be observed that a near optimal rate of convergence $h^{k+1}$ is obtained for all the element types and degrees of approximation considered. 

\subsection{Super--convergence of the displacement field}
\label{sc:superconvergenceStudy}

In this section, the three post--process procedures described in Section~\ref{sc:superconvergence} are tested using numerical examples. It is worth recalling that the three post--process options differ in the condition used to remove the indeterminacy related to the rigid rotational modes.

\subsubsection{Two dimensional example}
\label{sc:superconvergence2D}

The different post--process techniques are applied to the two dimensional example of Section~\ref{sc:convergence2D}.

\begin{figure}[!tb]
	\centering
	\subfigure[Quadrilaterals]{\includegraphics[width=0.4\textwidth]{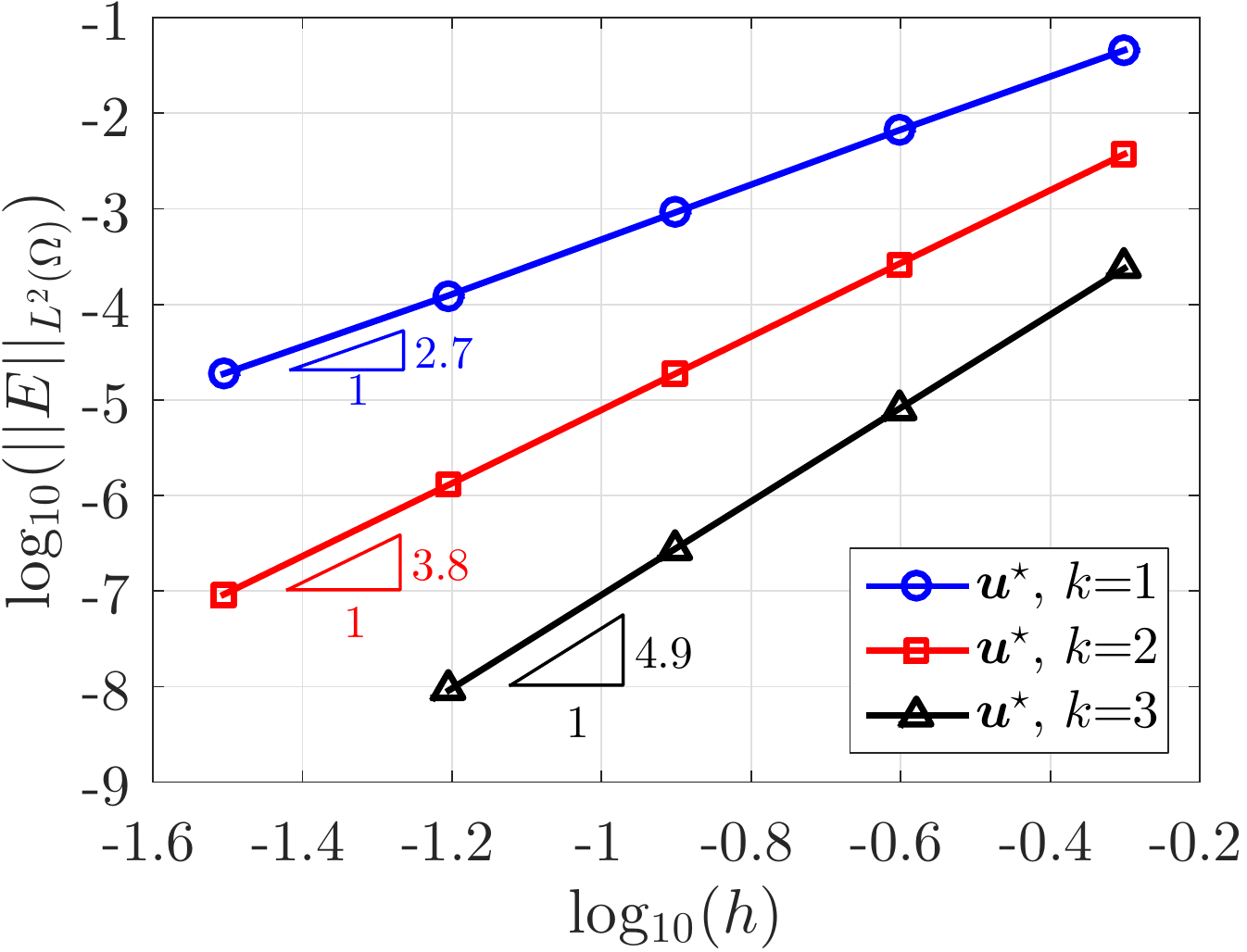}}
	\subfigure[Triangles]{\includegraphics[width=0.4\textwidth]{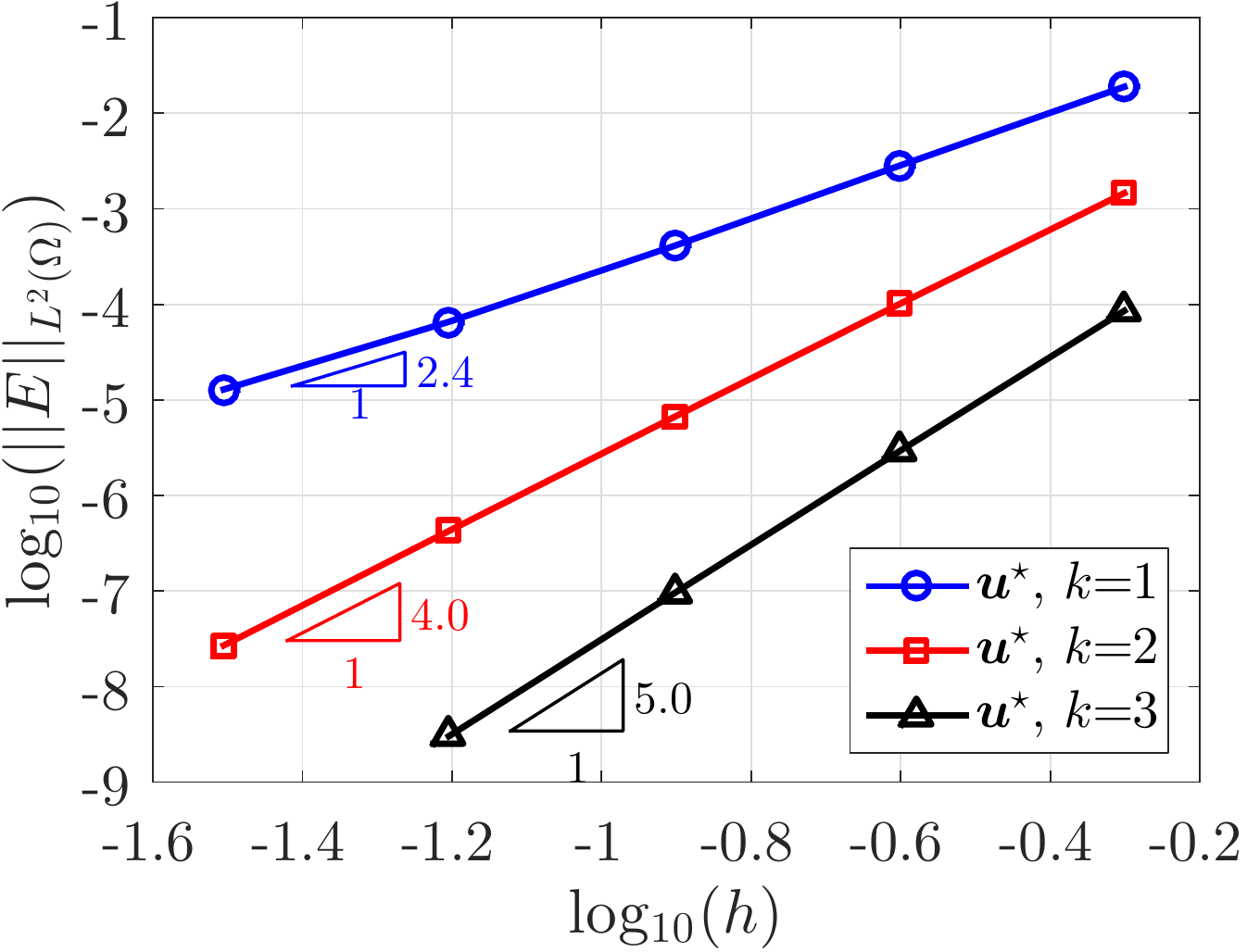}}
	\caption{Two dimensional problem: $h$--convergence of the error of the post--processed solution in the $\eltwo(\Omega)$ norm for quadrilateral and triangular meshes with different orders of approximation using the post--process technique of Equation~\eqref{eq:postprocess1}.}
	\label{fig:hConv2DstarOpt1}
\end{figure}

The first post--process considers the condition of Equation~\eqref{eq:postprocess1}. The convergence of the error of the post--processed variable $\bu^\star$, measured in the $\eltwo(\Omega)$ norm, as a function of the characteristic element size $h$ is represented in Figure~\ref{fig:hConv2DstarOpt1} for both quadrilateral and triangular elements and for a degree of approximation ranging from $k=1$ up to $k=3$. 
The results indicate that, as other HDG methods for linear elasticity~\cite{soon2009hybridizable}, super--convergence of the post--processed solution is obtained for $k \geq 2$. When a linear approximation is used, quadrilateral elements show almost optimal convergence but for triangular elements a sub--optimal rate of 2.4 is observed. 

Comparing the errors of the post--processed solution to the errors of the HDG solution in Figure~\ref{fig:hConv2D}, it is apparent that, despite no super--convergent results are provided by the first post--process technique, the post--processed solution is substantially more accurate than the HDG solution for both quadrilateral and triangular elements.
It is worth noting that this post--process was utilised in a different HDG formulation of the linear elastic problem for linear triangles~\cite{soon2009hybridizable} and sub--optimal convergence was also observed.

\begin{figure}[!tb]
	\centering
	\subfigure[Quadrilaterals]{\includegraphics[width=0.4\textwidth]{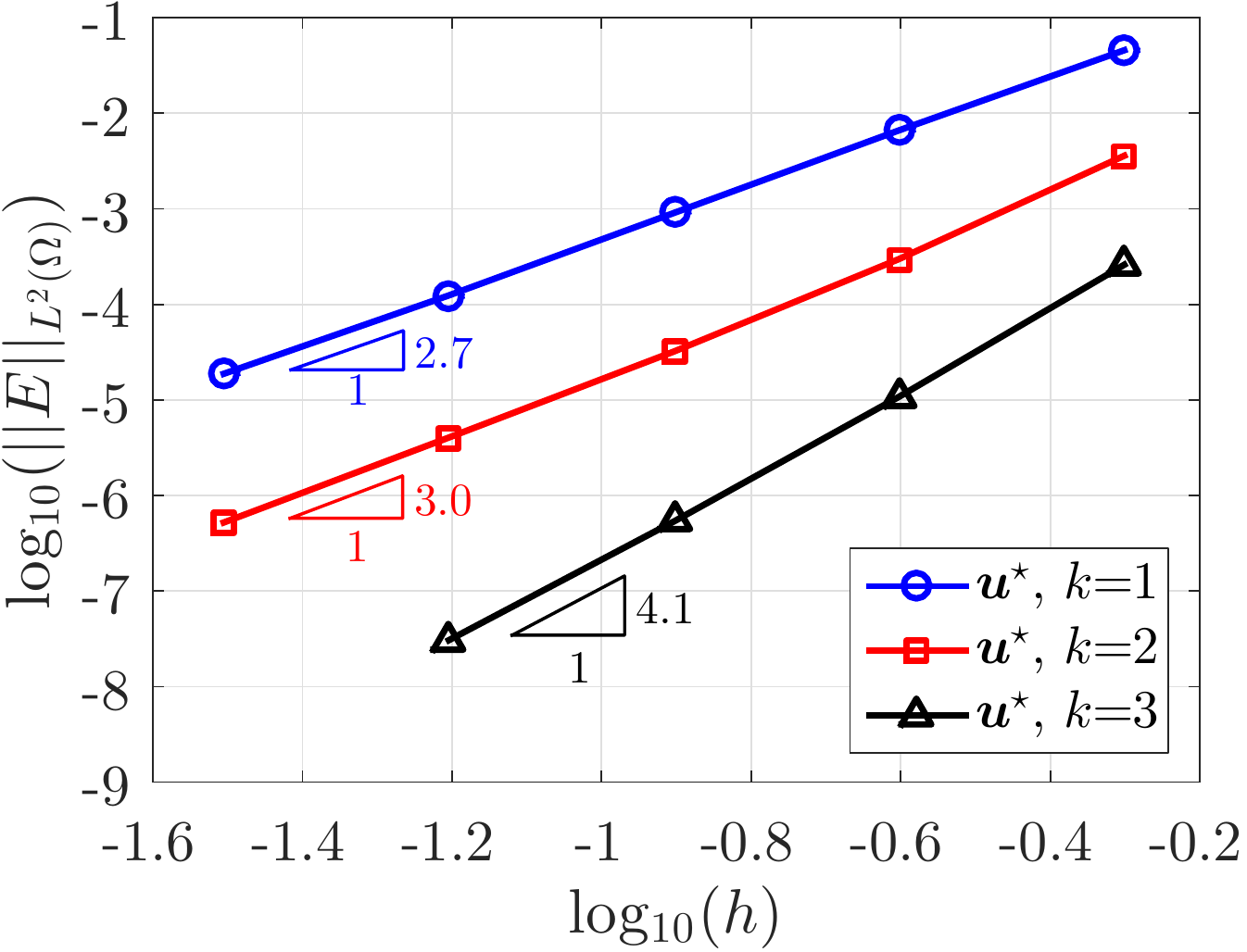}}
	\subfigure[Triangles]{\includegraphics[width=0.4\textwidth]{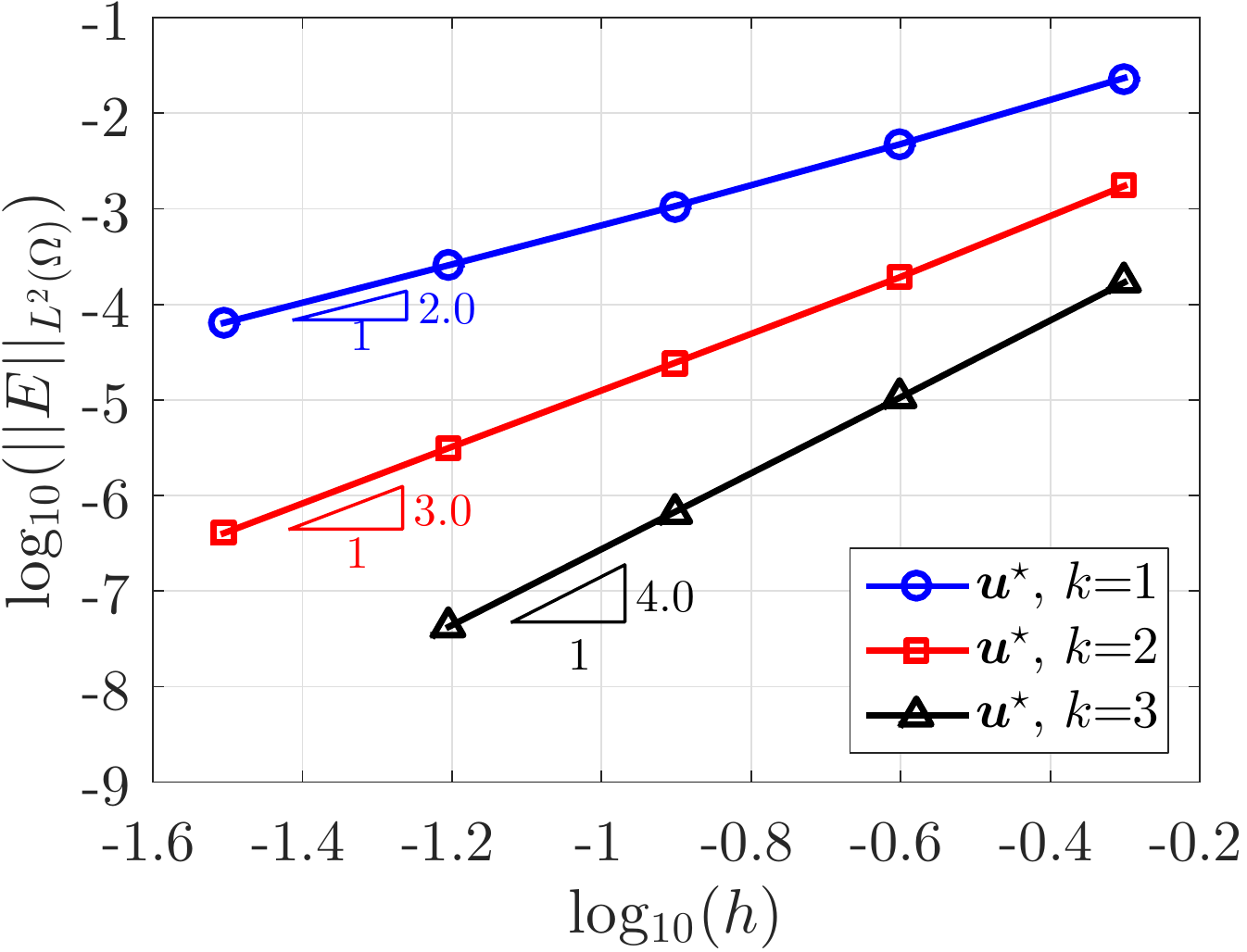}}	
	\caption{Two dimensional problem: $h$--convergence of the error of the post--processed solution in the $\eltwo(\Omega)$ norm for quadrilateral and triangular meshes with different orders of approximation using the post--process technique of Equation~\eqref{eq:postprocess2}.}
	\label{fig:hConv2DstarOpt2}
\end{figure}

Next, the post--process that considers the condition of Equation~\eqref{eq:postprocess2} is tested. Figure~\ref{fig:hConv2DstarOpt2} shows the convergence study for the error of the post--processed variable $\bu^\star$ measured in the $\eltwo(\Omega)$ norm.
The results for quadrilateral elements are almost identical to the results obtained with the first technique, whereas, for triangles, a sub--optimal order $k+1$ is observed for all degrees of approximation. 

Comparing the errors of the post--processed solution with triangles to the errors of the HDG solution in Figure~\ref{fig:hConv2D}, it is apparent that little gain in accuracy is obtained with the post--processed solution. This is crucial when the super--convergent solution is sought to devise automatic degree adaptive processes~\cite{giorgiani2014hybridizable, RS-AH:18} and suggests that the post--process provided by the second option cannot be used to produce an accurate error estimator with triangles.

\begin{figure}[!tb]
	\centering
	\subfigure[Quadrilaterals]{\includegraphics[width=0.4\textwidth]{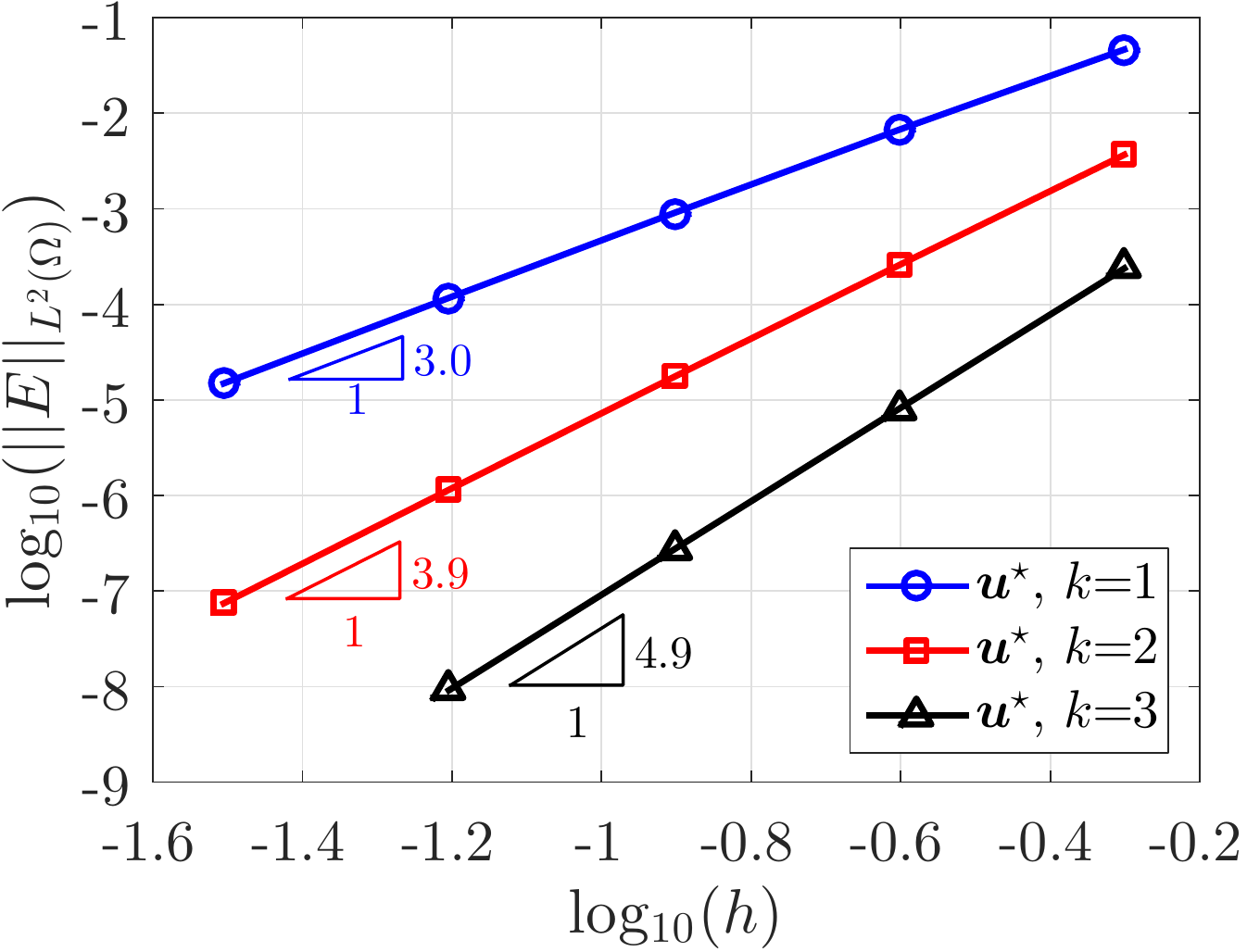}}
	\subfigure[Triangles]{\includegraphics[width=0.4\textwidth]{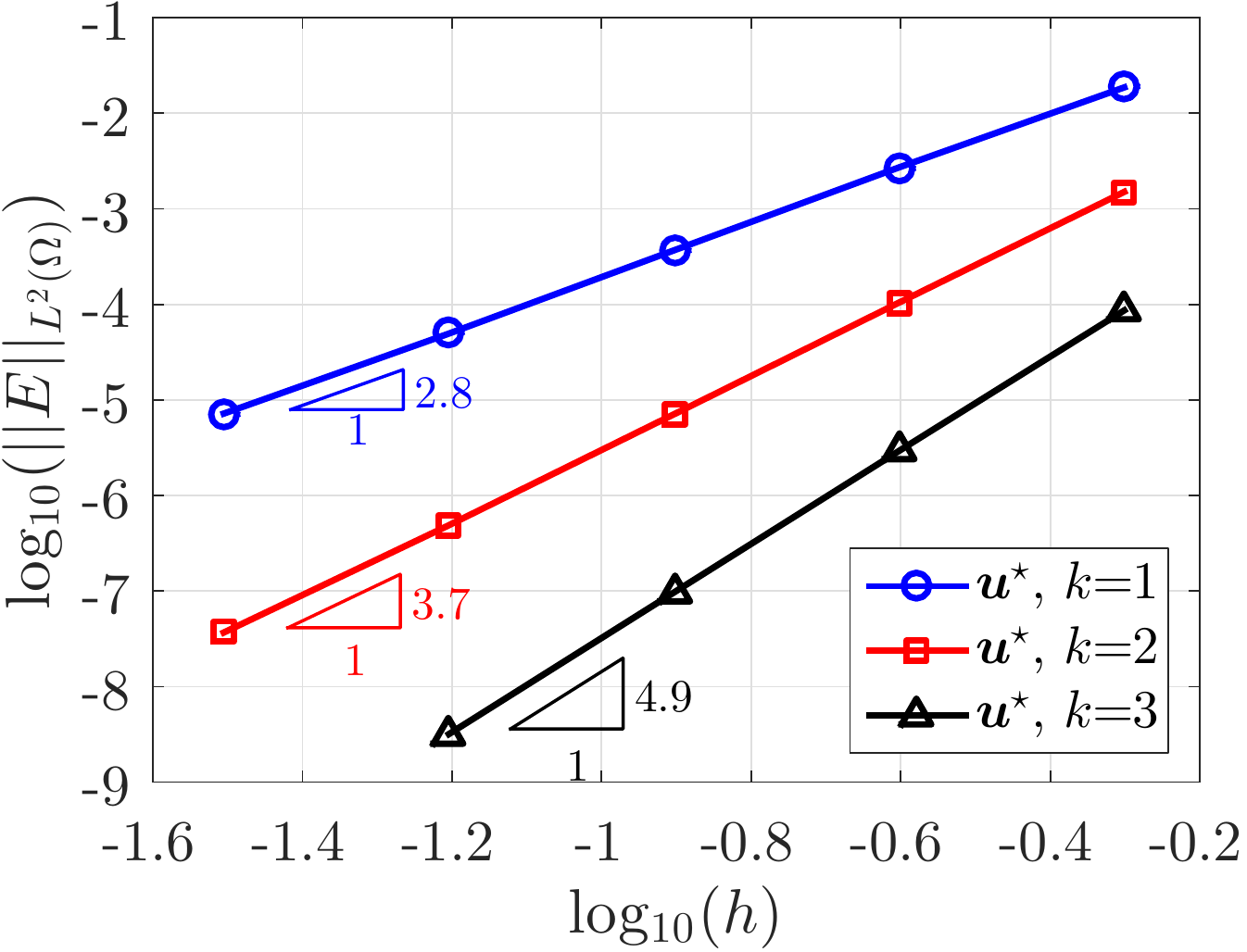}}	
	\caption{Two dimensional problem: $h$--convergence of the error of the post--processed solution in the $\eltwo(\Omega)$ norm for quadrilateral and triangular meshes with different orders of approximation using the post--process technique of Equation~\eqref{eq:postprocess3}.}
	\label{fig:hConv2DstarOpt3}
\end{figure}

The last post--process technique proposed in this paper is considered, consisting of imposing the condition of Equation~\eqref{eq:postprocess3}. The convergence of the error of the post--processed variable $\bu^\star$, measured in the $\eltwo(\Omega)$ norm, as a function of the characteristic element size $h$ is represented in Figure~\ref{fig:hConv2DstarOpt3}.
The results reveal that almost the optimal rate of convergence is attained for both quadrilateral and triangular elements and for all degrees of approximation. This indicates that the average of the hybrid variable on the boundary leads to super--convergent results. 
It is worth noting that the error of the post--processed solution obtained with the third post--process technique, proposed here, is not only showing the optimal rate but it also provides an extra gain in accuracy when compared to the first post--process technique, previously used in an HDG context. 

\subsubsection{Three dimensional example}
\label{sc:superconvergence3D}

The different post--process techniques are considered in the three dimensional example of Section~\ref{sc:convergence3D}.

\begin{figure}[!tb]
	\centering
	\subfigure[Hexahedrons]{\includegraphics[width=0.4\textwidth]{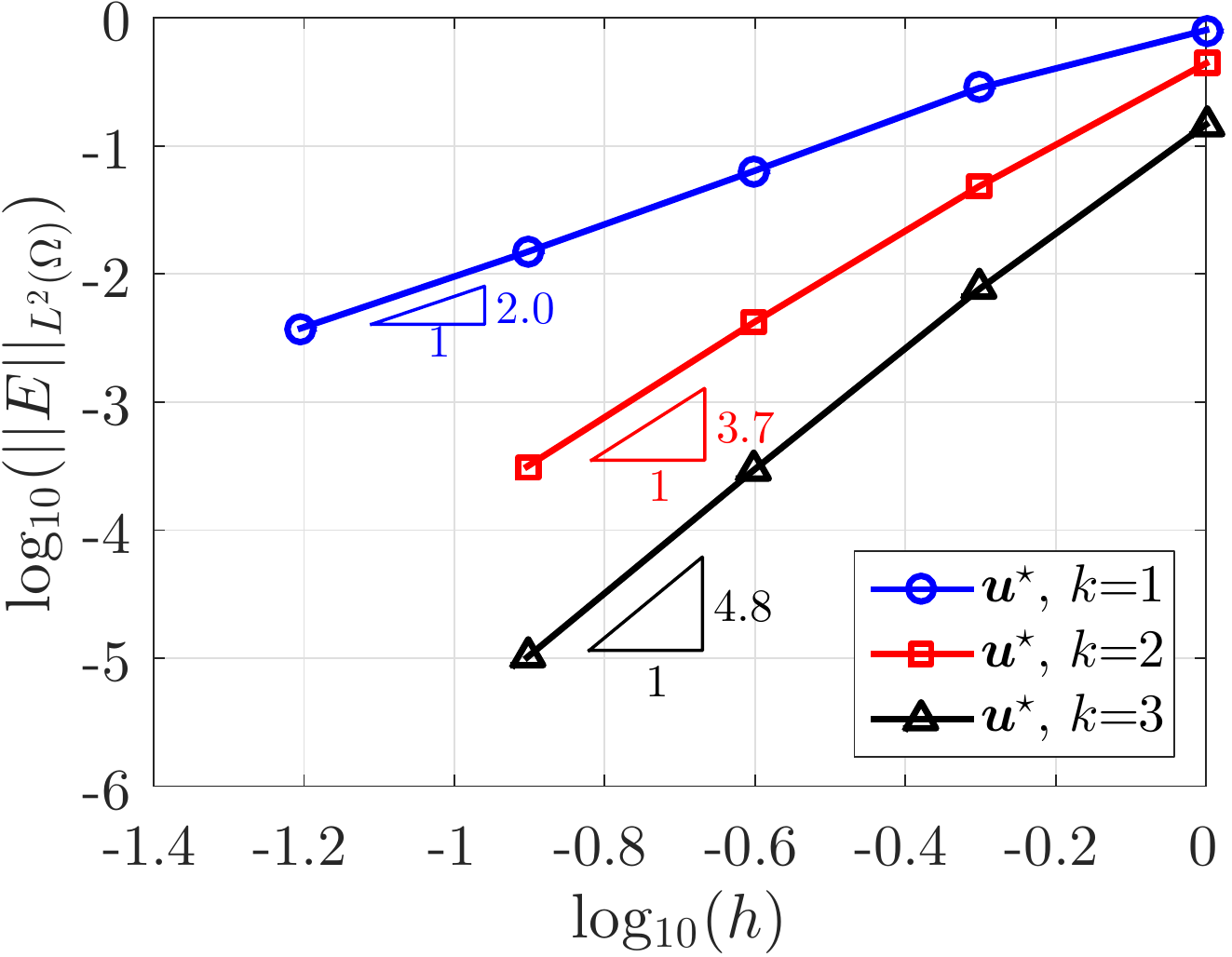}}
	\subfigure[Tetrahedrons]{\includegraphics[width=0.4\textwidth]{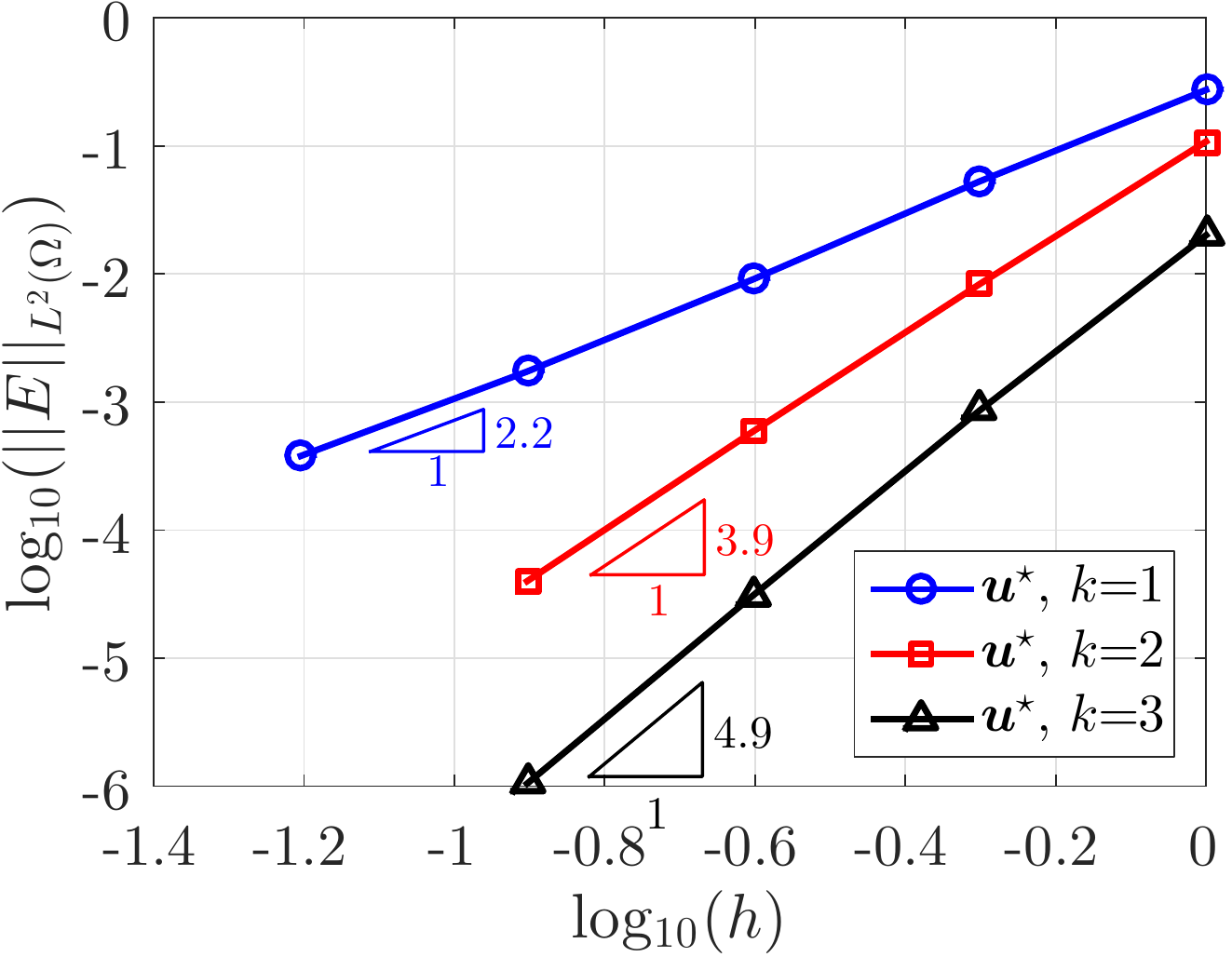}}
	\subfigure[Prisms]{\includegraphics[width=0.4\textwidth]{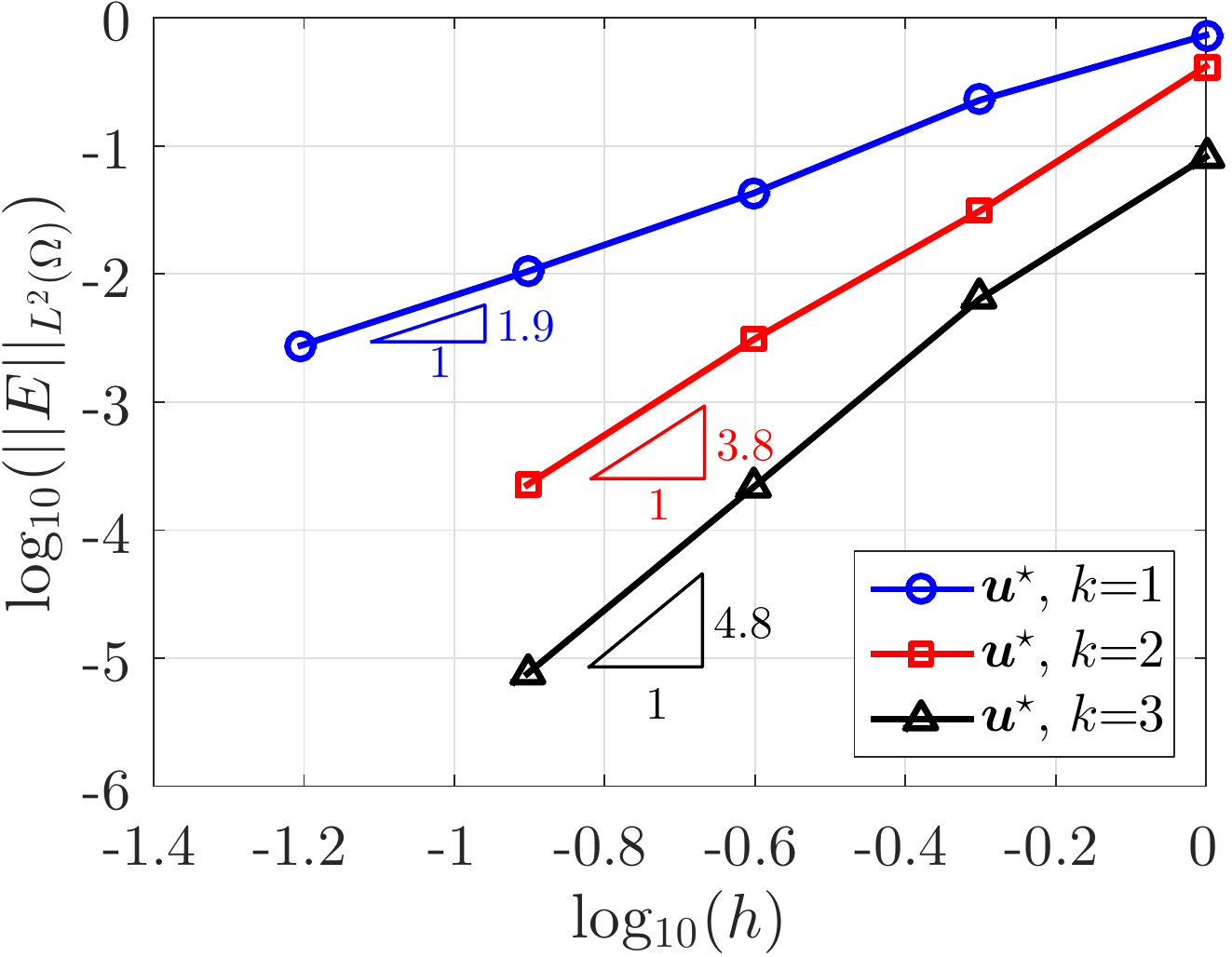}}
	\subfigure[Pyramids]{\includegraphics[width=0.4\textwidth]{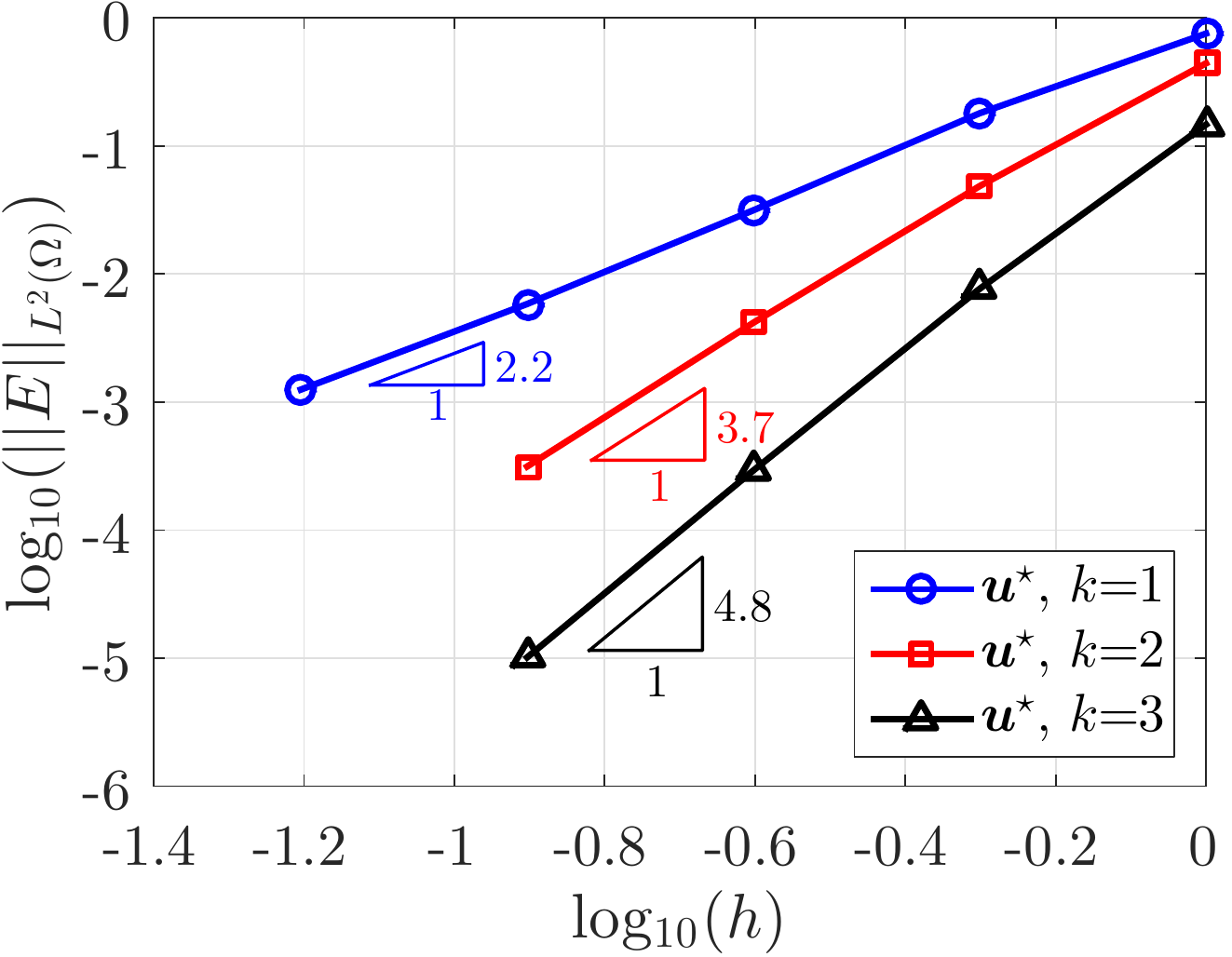}}
	\caption{Three dimensional problem: $h$--convergence of the error of the post--processed solution in the $\eltwo(\Omega)$ norm for hexahedral, tetrahedral, prismatic and pyramidal meshes with different orders of approximation using the post--process technique of Equation~\eqref{eq:postprocess1}.}
	\label{fig:hConv3DstarOpt1}
\end{figure}

The convergence of the error of the post--processed variable $\bu^\star$, measured in the $\eltwo(\Omega)$ norm, as a function of the characteristic element size $h$ is represented in Figure~\ref{fig:hConv3DstarOpt1} when using the post--process technique of Equation~\eqref{eq:postprocess1}. 
The results reveal that super--convergent results are obtained with $k\geq 2$ whereas with linear elements sub--optimal convergence is attained. It is worth noting that in the two dimensional example almost super--convergent results where obtained with quadrilateral elements whereas in three dimensions sub--optimal convergence of order $k+1$ is observed for all the different element types considered.

\begin{figure}[!tb]
	\centering
	\subfigure[Hexahedrons]{\includegraphics[width=0.4\textwidth]{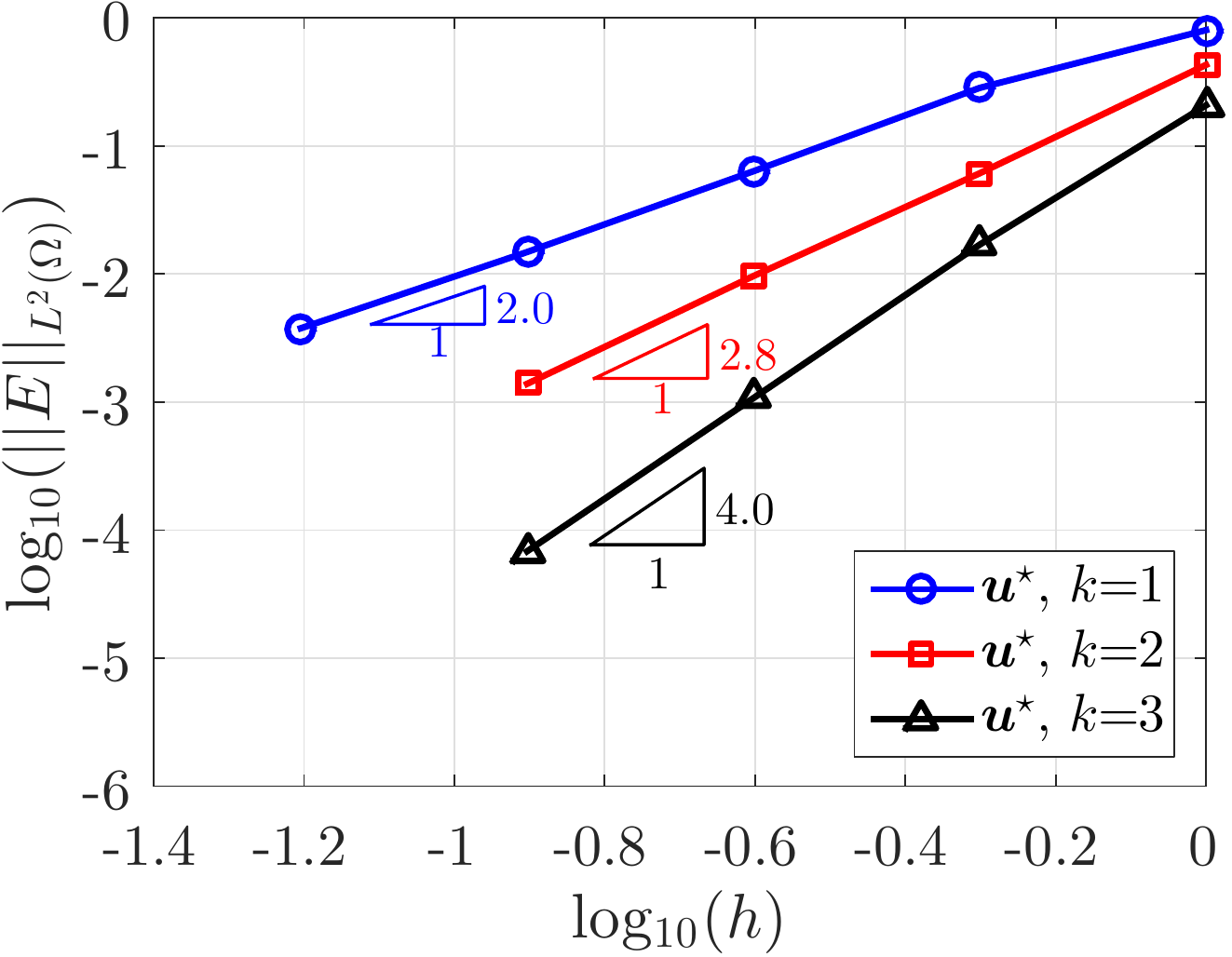}}
	\subfigure[Tetrahedrons]{\includegraphics[width=0.4\textwidth]{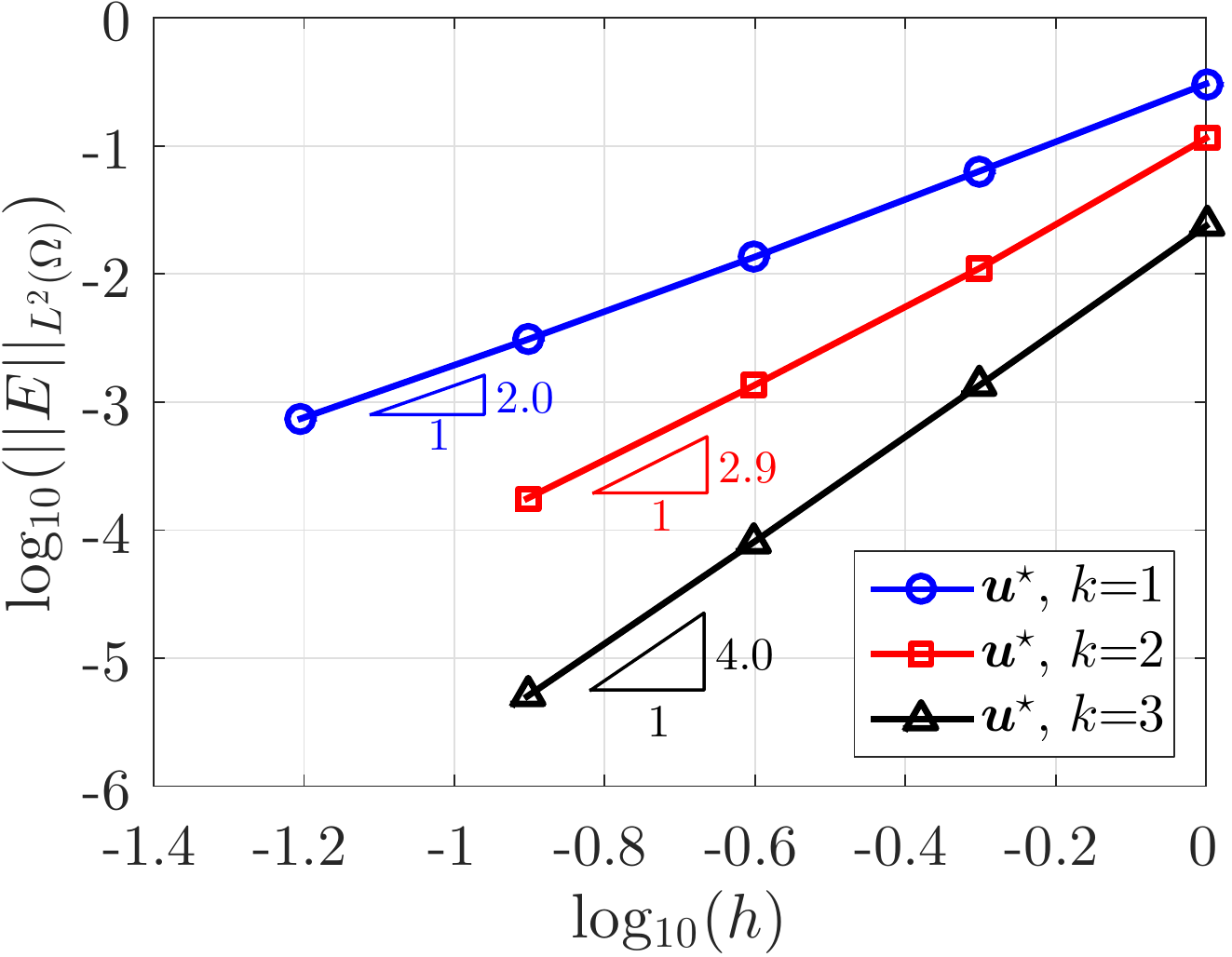}}
	\caption{Three dimensional problem: $h$--convergence of the error of the post--processed solution in the $\eltwo(\Omega)$ norm for hexahedral and tetrahedral meshes with different orders of approximation using the post--process technique of Equation~\eqref{eq:postprocess2}.}
	\label{fig:hConv3DstarOpt2}
\end{figure}

Next, the post--process that considers the condition of Equation~\eqref{eq:postprocess2} is tested. Figure~\ref{fig:hConv3DstarOpt2} shows the convergence study for the error of the post--processed variable $\bu^\star$ measured in the $\eltwo(\Omega)$ norm.
Only the results for hexahedral and tetrahedral elements are reported in Figure~\ref{fig:hConv3DstarOpt2} because, analogously to the two dimensional example, this post--process leads to a sub--optimal rate $k+1$ for all the different elements types and degrees of approximation.

\begin{figure}[!tb]
	\centering
	\subfigure[Hexahedrons]{\includegraphics[width=0.4\textwidth]{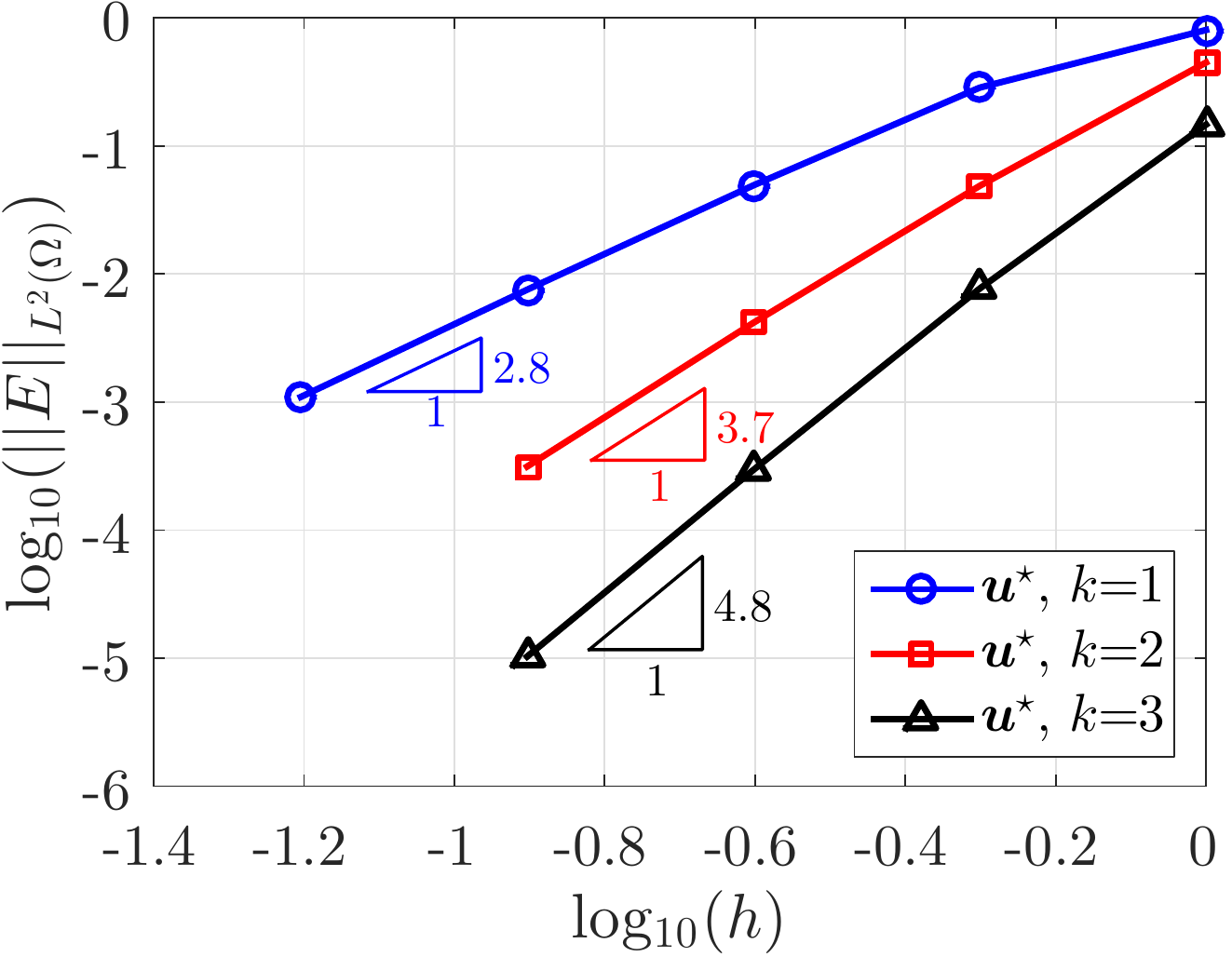}}
	\subfigure[Tetrahedrons]{\includegraphics[width=0.4\textwidth]{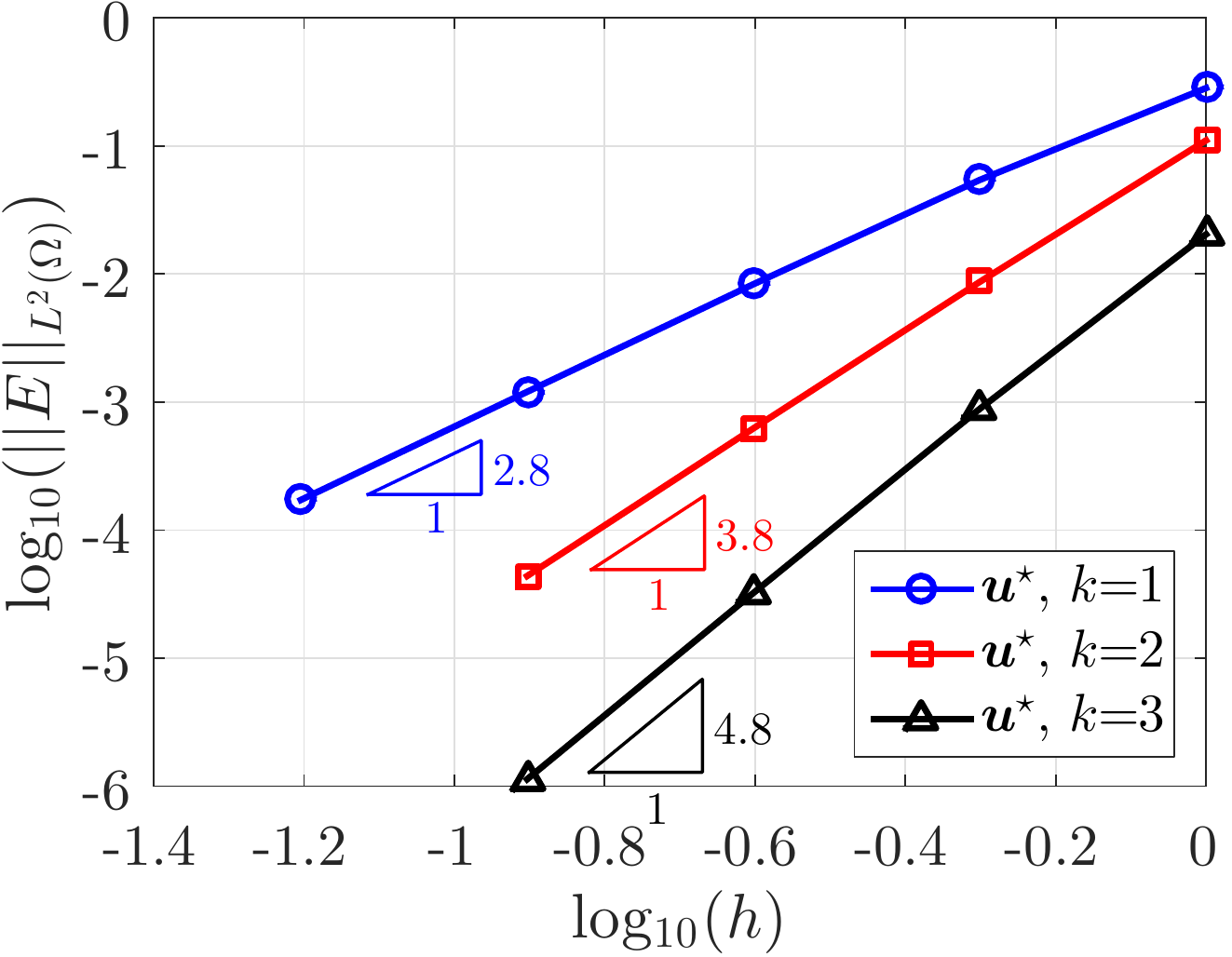}}
	\subfigure[Prisms]{\includegraphics[width=0.4\textwidth]{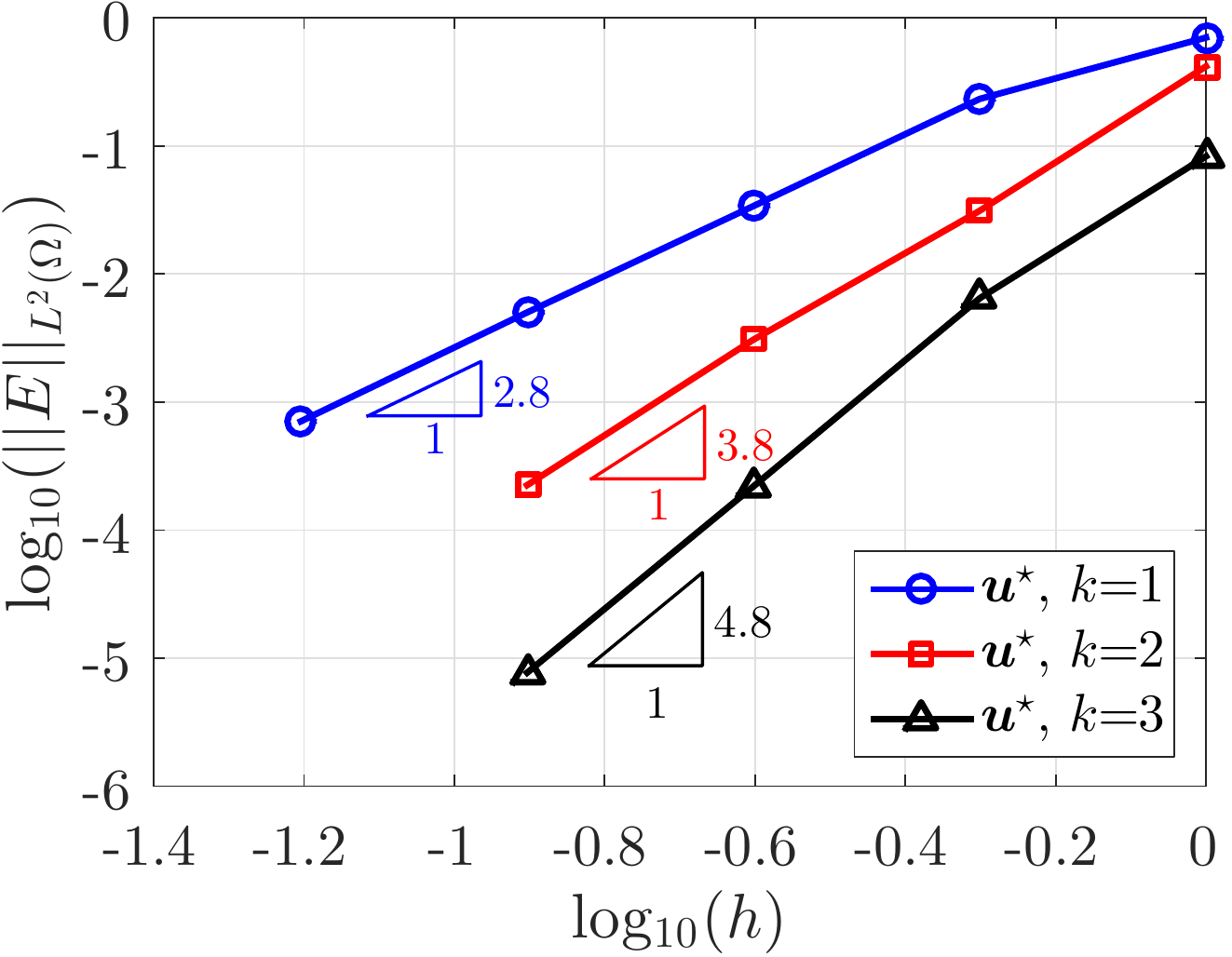}}
	\subfigure[Pyramids]{\includegraphics[width=0.4\textwidth]{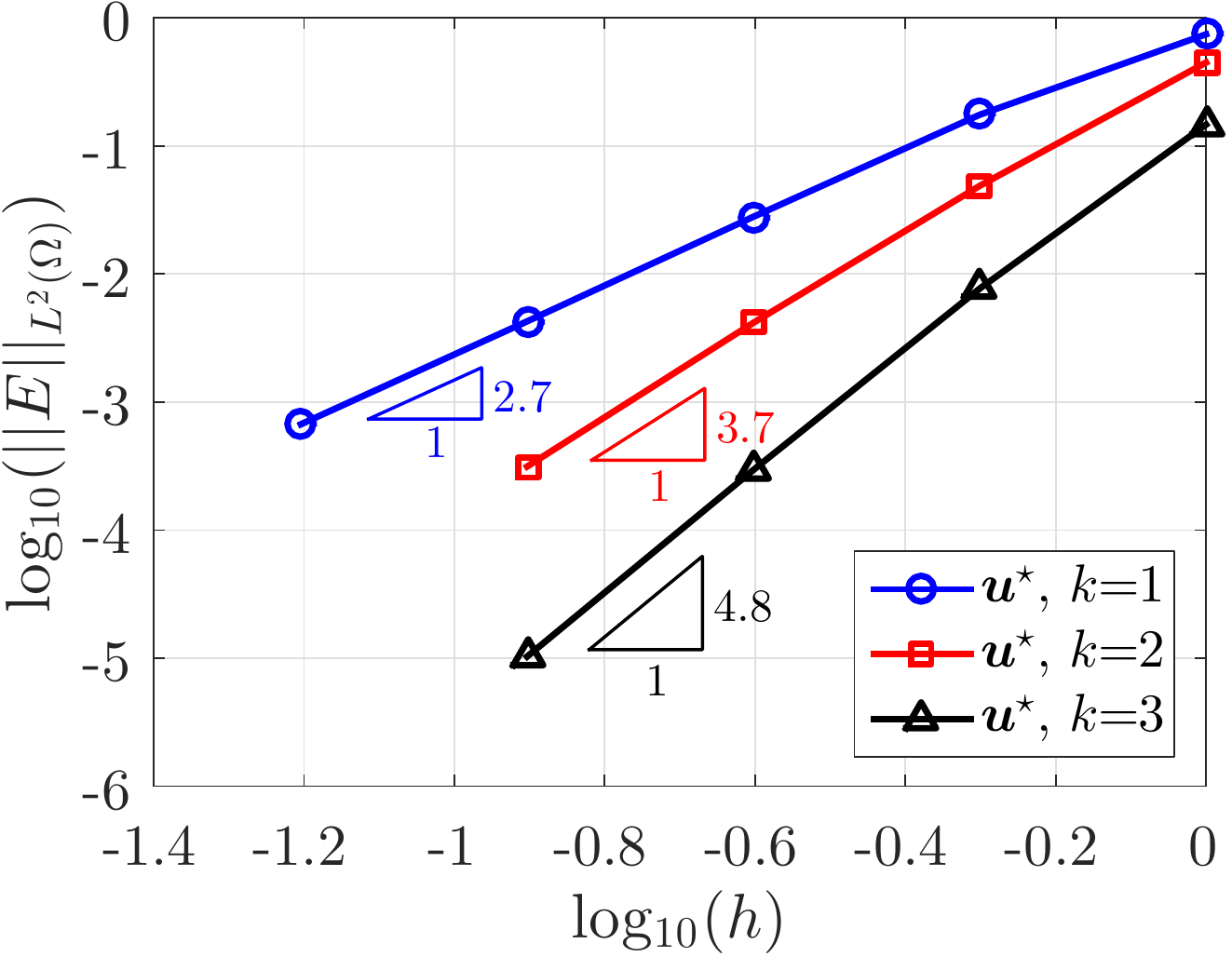}}
	\caption{Three dimensional problem: $h$--convergence of the error of the post--processed solution in the $\eltwo(\Omega)$ norm for hexahedral, tetrahedral, prismatic and pyramidal meshes meshes with different orders of approximation using the post--process technique of Equation~\eqref{eq:postprocess3}.}
	\label{fig:hConv3DstarOpt3}
\end{figure}

Finally, the last post--process technique proposed in this paper is considered. The convergence of the error of the post--processed variable $\bu^\star$, measured in the $\eltwo(\Omega)$ norm, as a function of the characteristic element size $h$ is represented in Figure~\ref{fig:hConv3DstarOpt3}.
The results show that almost the optimal rate of convergence is attained for all the element types and for all degrees of approximation considered. 
The numerical experiments performed in two and three dimensions confirm that the post--process technique proposed in this paper for the first time lead to optimal super--convergent results of the primal variable. 

As in the two dimensional example, the error of the post--processed solution obtained with the third technique, proposed here, is not only showing the optimal rate but it also provides an extra gain in accuracy when compared to the first post--process technique, previously used in an HDG context. When compared to the error of the HDG solution, represented in Figure~\ref{fig:hConv3D}, the post--process proposed here provides a solution that is almost one order of magnitude more accurate than the HDG solution, even for linear approximation of the solution.

\subsection{Influence of the stabilisation parameter}
\label{sc:tau}

The stabilisation tensor $\btau$ is known to have an important effect on the stability, accuracy and convergence properties of the resulting HDG method~\cite{Jay-CGL:09,Cockburn-CDG:08,soon2009hybridizable}. This section presents a numerical study to assess the influence of the stabilisation parameter on the accuracy of the results. For simplicity, it is assumed that $\btau = \tau \Insd$ and the influence of the scalar stabilisation parameter $\tau$ is investigated.

\subsubsection{Two dimensional example}
\label{sc:tau2D}

\begin{figure}[!tb]
	\centering
	\subfigure[Quadrilaterals]{\includegraphics[width=0.4\textwidth]{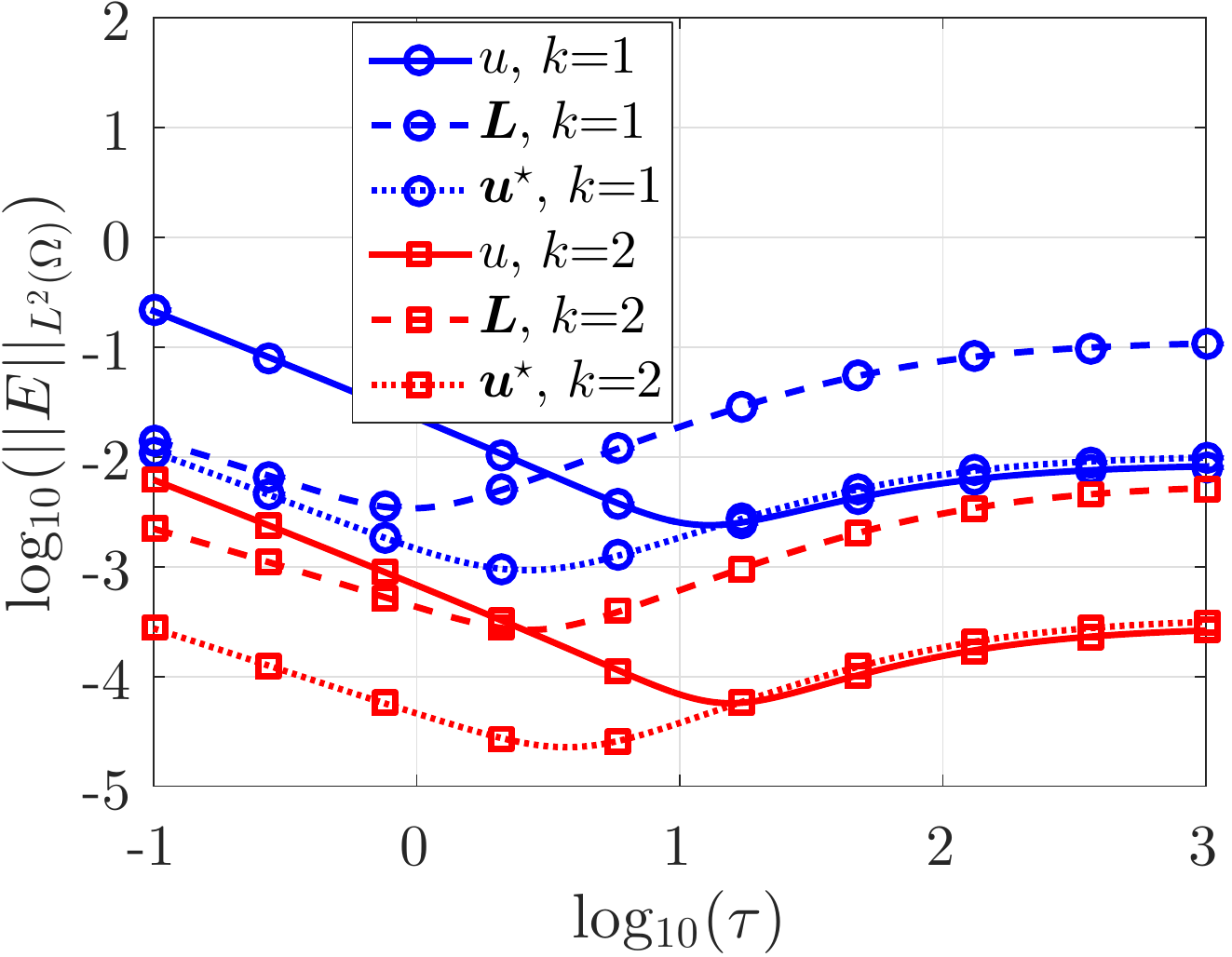}}
	\subfigure[Triangles]{\includegraphics[width=0.4\textwidth]{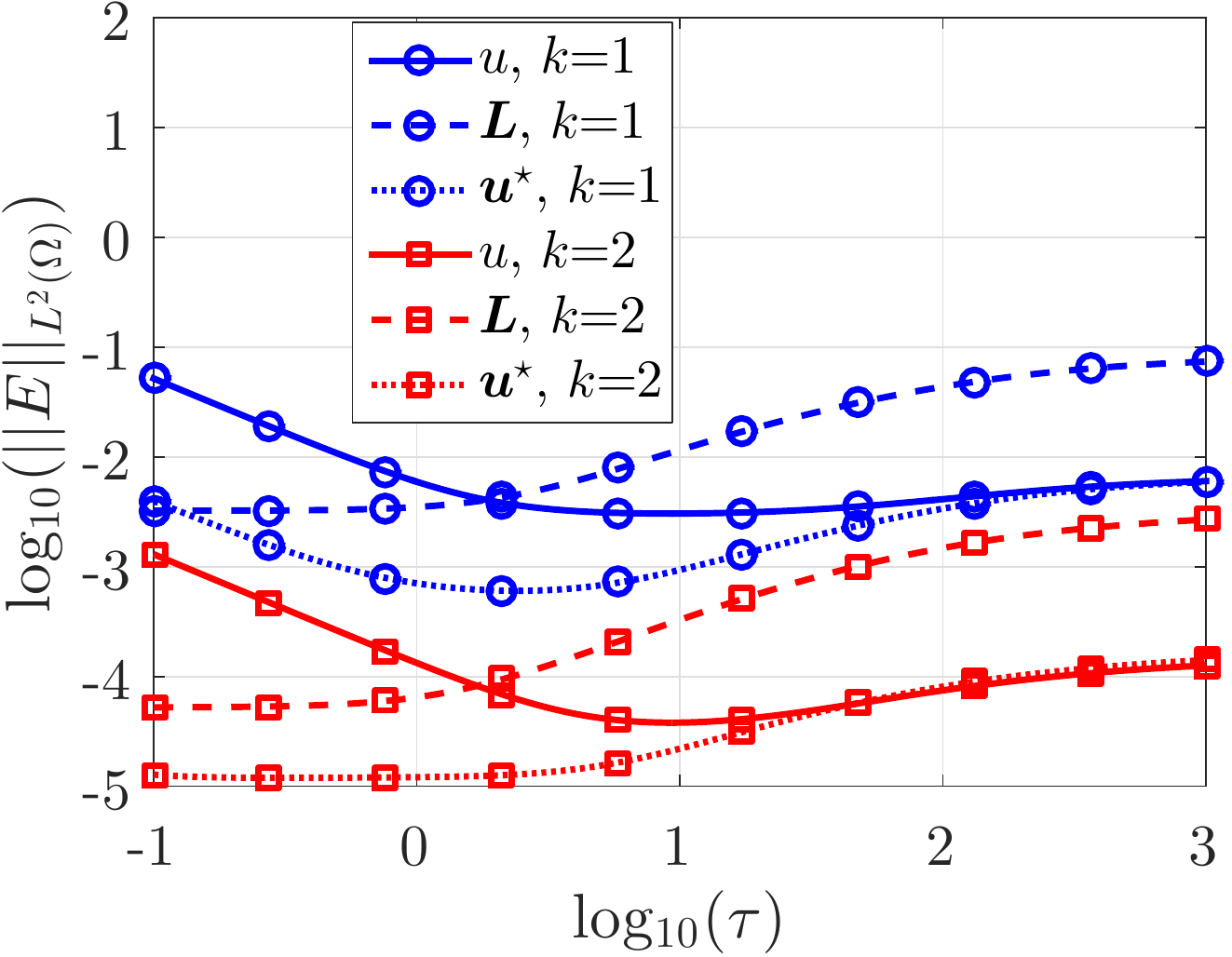}}	
	\caption{Two dimensional problem: error of the primal, mixed and post--processed variables, $\bu$, $\bL$ and $\bu^\star$ respectively, in the $\eltwo(\Omega)$ norm as a function of the stabilisation parameter $\tau$.}
	\label{fig:tauInfluence2D}
\end{figure}

Figure~\ref{fig:tauInfluence2D} shows the evolution of the error of the primal, mixed and post--processed variables, $\bu$, $\bL$ and $\bu^\star$ respectively, in the $\eltwo(\Omega)$ norm as a function of the stabilisation parameter $\tau$ for the two dimensional example studied in Section~\ref{sc:convergence2D}. 
The numerical experiment is performed using linear and quadratic approximations and for quadrilateral and triangular meshes and the value of $\tau$ varies from 0.1 to 1,000.

The results reveal that there is a value of $\tau$ for which the error of the primal solution is minimum. For both quadrilateral and triangular meshes with linear and quadratic approximation, this value is near $\tau = 10$. However, it is worth noting that for $\tau = 10$ the post--process of the displacement field offers little or no extra gain in accuracy. When the error of the mixed variable is of interest, the minimum error is achieved for a different value of the stabilisation parameter, near $\tau = 3$. It is worth noting that the value of $\tau = 3$ also provides the best accurate results for the post--processed variable. As a result, the value of $\tau = 3$ is considered the optimum value in this experiment as it provides the most accurate solution for both the displacement (i.e. the post--processed variable) and the stress (i.e. the mixed variable). It is also interesting to observe that for $\tau = 3$ the accuracy on the primal and mixed variables is almost identical.

\subsubsection{Three dimensional example}
\label{sc:tau3D}

\begin{figure}[!tb]
	\centering
	\subfigure[Hexahedrons]{\includegraphics[width=0.4\textwidth]{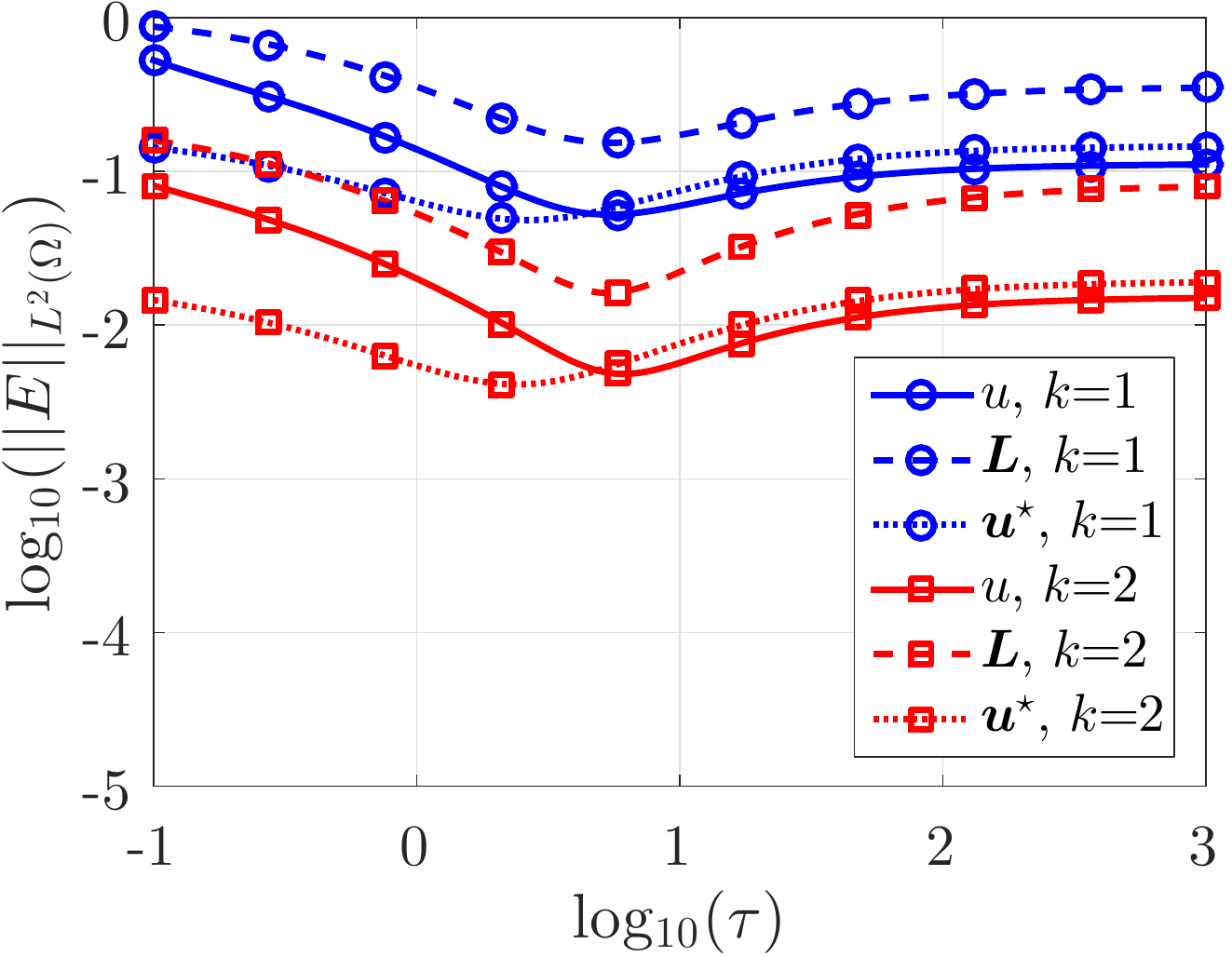}}
	\subfigure[Tetrahedrons]{\includegraphics[width=0.4\textwidth]{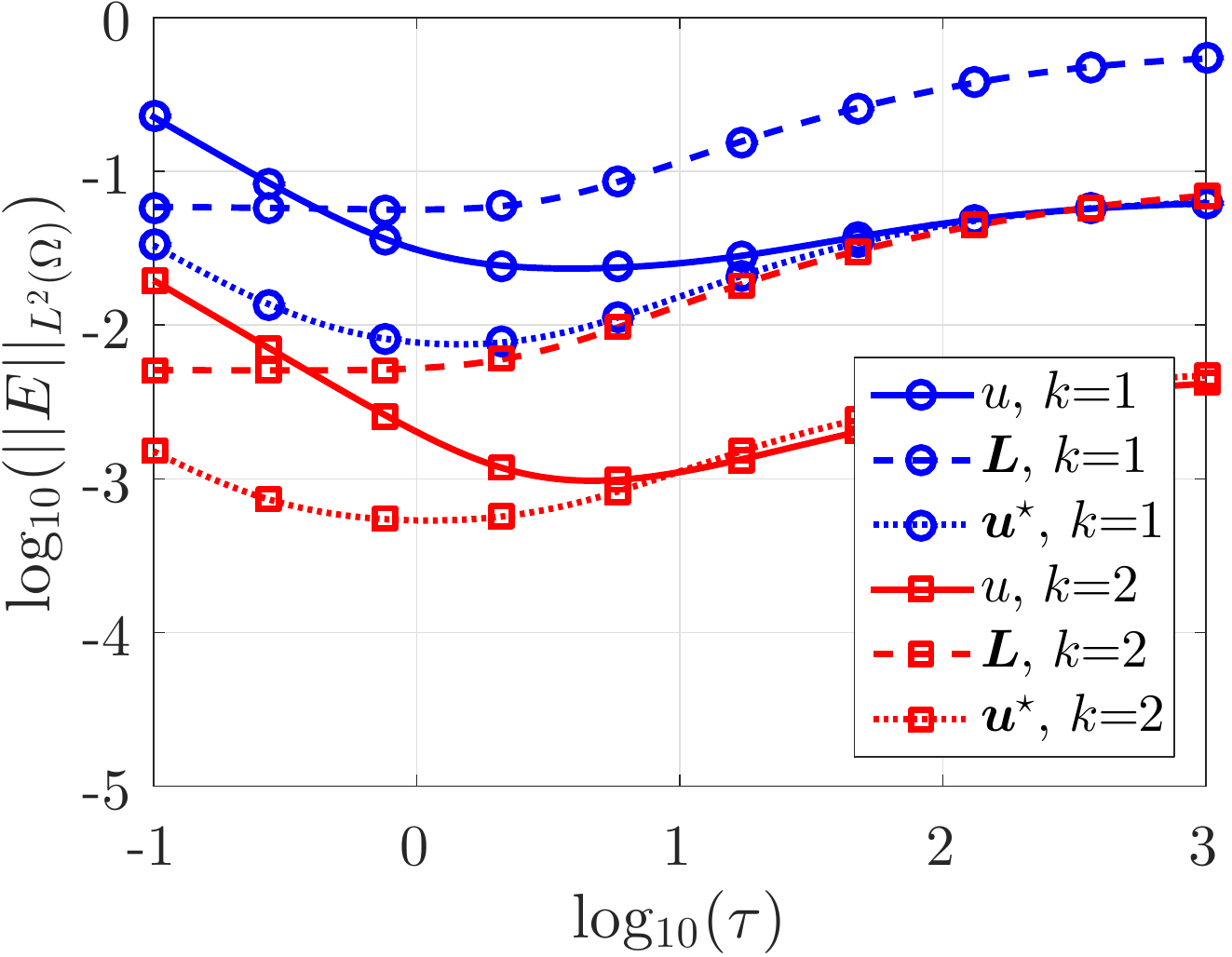}}
	\subfigure[Prisms]{\includegraphics[width=0.4\textwidth]{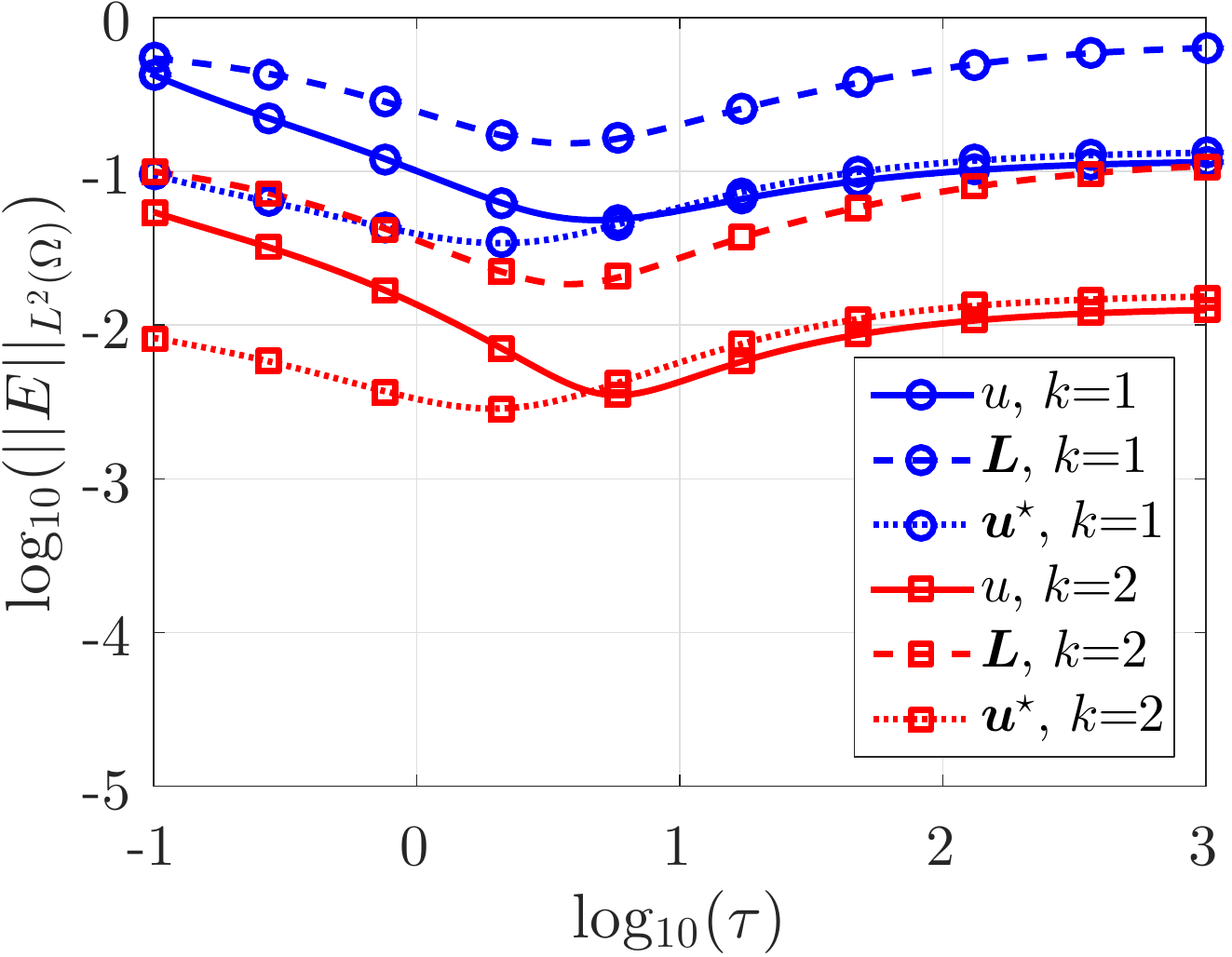}}
	\subfigure[Pyramids]{\includegraphics[width=0.4\textwidth]{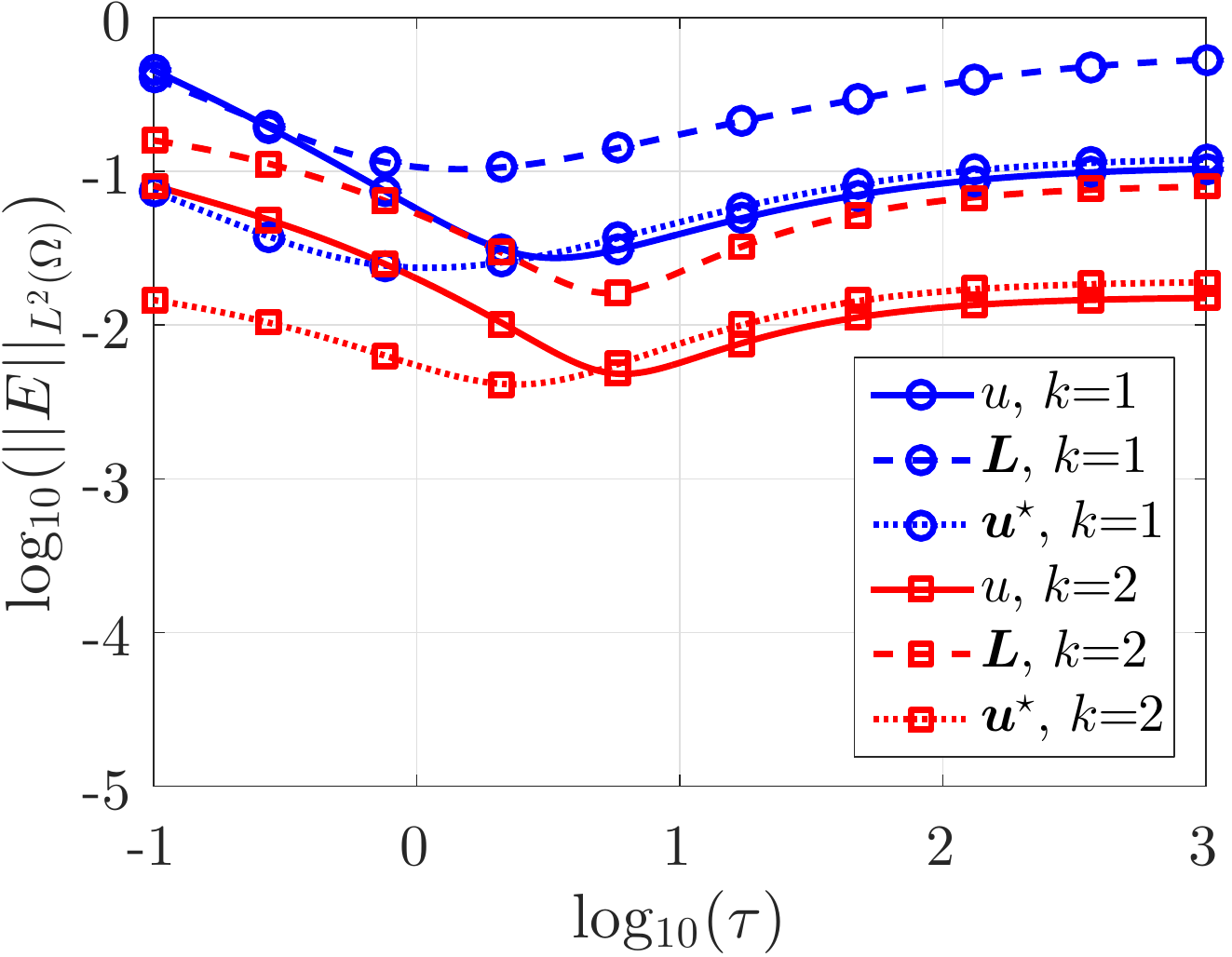}}
	\caption{Three dimensional problem: error of the primal, mixed and post--processed variables in the $\eltwo(\Omega)$ norm as a function of the stabilisation parameter.}
	\label{fig:tauInfluence3D}
\end{figure}

A similar study is performed next for the three dimensional example of Section~\ref{sc:convergence3D}. Figure~\ref{fig:tauInfluence3D} shows the evolution of the error of the primal, mixed and post--processed variables, $\bu$, $\bL$ and $\bu^\star$ respectively, in the $\eltwo(\Omega)$ norm as a function of the stabilisation parameter $\tau$.
As in the two dimensional example a value of the stabilisation parameter near $\tau=10$ provides the minimum error for the primal variable but with no extra gain in accuracy when the post--process is performed. Also, as in the two dimensional case, the value of $\tau$ that provides the most accurate results for both the mixed and the post--processed variables is near $\tau = 3$. It is worth noting that for hexahedral, prismatic and pyramidal elements, and contrary to the results obtained in the two dimensional problem, the value that provides the most accurate results for the primal and mixed variable is almost identical.

The conclusions that are extracted from this study are similar to the ones obtained in the two dimensional example and show that the optimal value of the stabilisation parameter is not dependent upon the degree of approximation, the type of element or the dimensionality of the problem.

\subsection{Locking--free behaviour in the incompressible limit}
\label{sc:incompressible}

The last example considers a problem with a nearly incompressible material (i.e. $\nu \sim 0.5$) that is commonly used in the literature~\cite{soon2009hybridizable}. The problem, defined in $\Omega = [0,1]^2$, has analytical solution given by 
\begin{equation}
\bu(\bx) = \Big( -x_1^2 x_2 (x_1-1)^2(x_2-1)(2x_2-1), \quad x_1^2 x_2 (x_1-1)^2(x_2-1)(2x_2-1) \Big).
\end{equation}
The external load and boundary conditions are derived from the exact solution. The Young's modulus is taken as $E = 3$ and the Poisson's ratio is varied from $\nu = 0.49$ up to $\nu = 0.49999$.

\begin{figure}[!tb]
	\centering
	\subfigure[Mesh 1]{\includegraphics[width=0.24\textwidth]{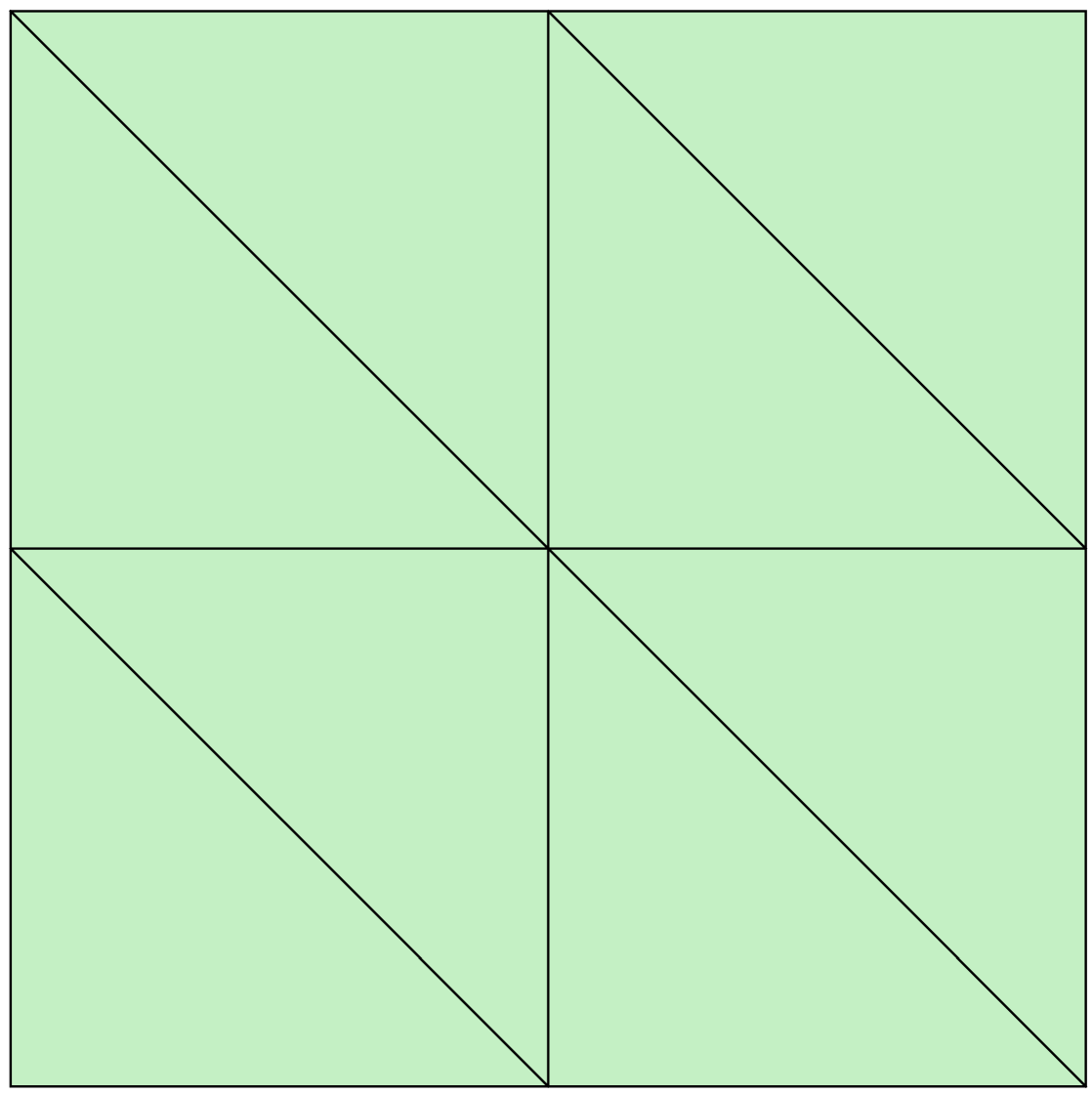}}
	\subfigure[Mesh 2]{\includegraphics[width=0.24\textwidth]{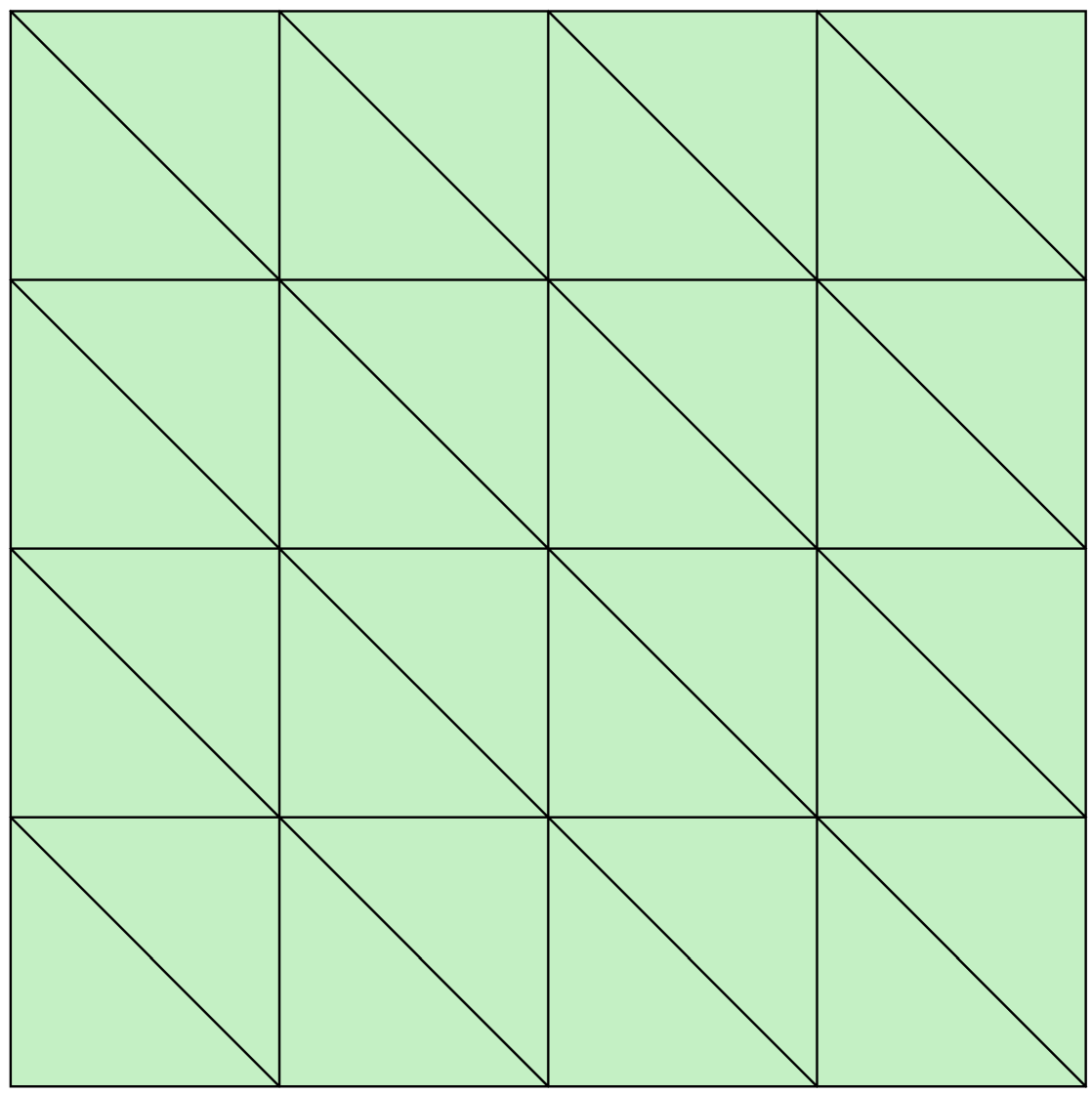}}
	\subfigure[Mesh 3]{\includegraphics[width=0.24\textwidth]{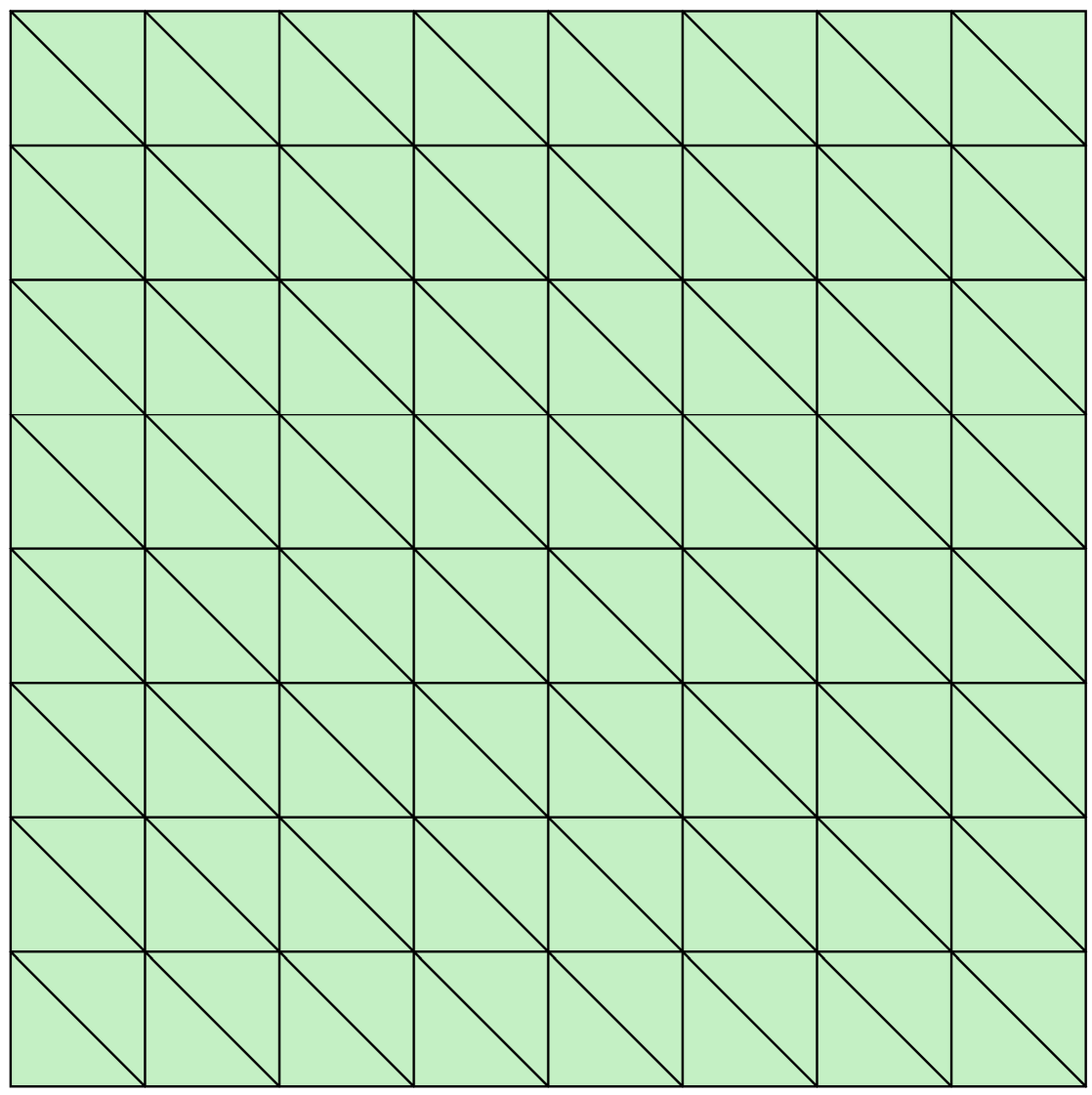}}
	\caption{Three two dimensional triangular meshes of $\Omega=[0,1]^2$ for the mesh convergence study with a nearly incompressible material.}
	\label{fig:TRImeshesOneDiag}
\end{figure}

Only triangular meshes with the arrangement represented in Figure~\ref{fig:TRImeshesOneDiag} are shown as this particular arrangement is known to exhibit a volumetric locking effect when considered with a traditional continuous Galerkin finite element formulation. Note that other the arrangements depicted in Figure~\ref{fig:2Dmeshes} also produce optimal rates of convergence.

\begin{figure}[!tb]
	\centering
	\subfigure[$u_1$]{\includegraphics[width=0.32\textwidth]{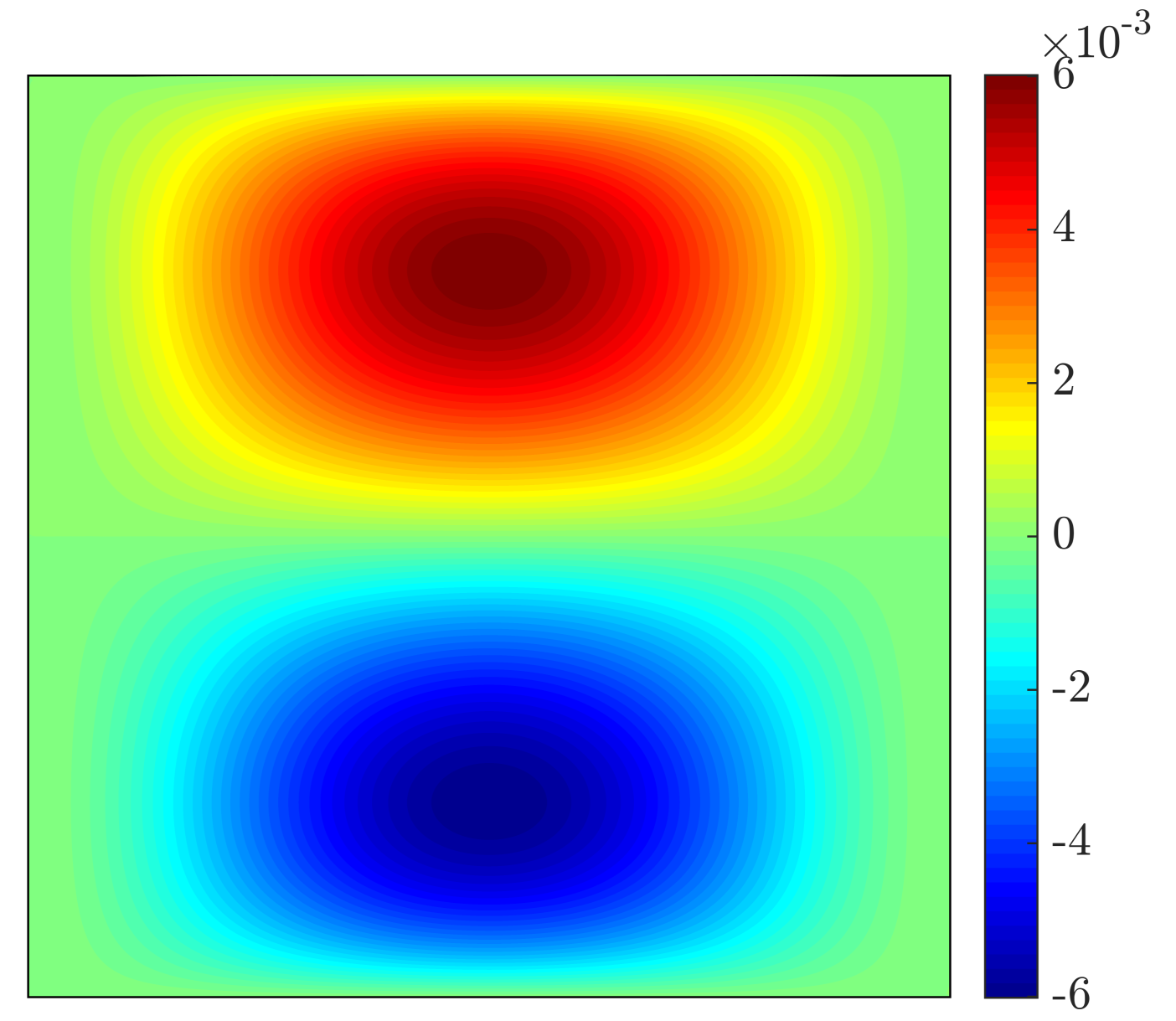}}
	\subfigure[$u_2$]{\includegraphics[width=0.32\textwidth]{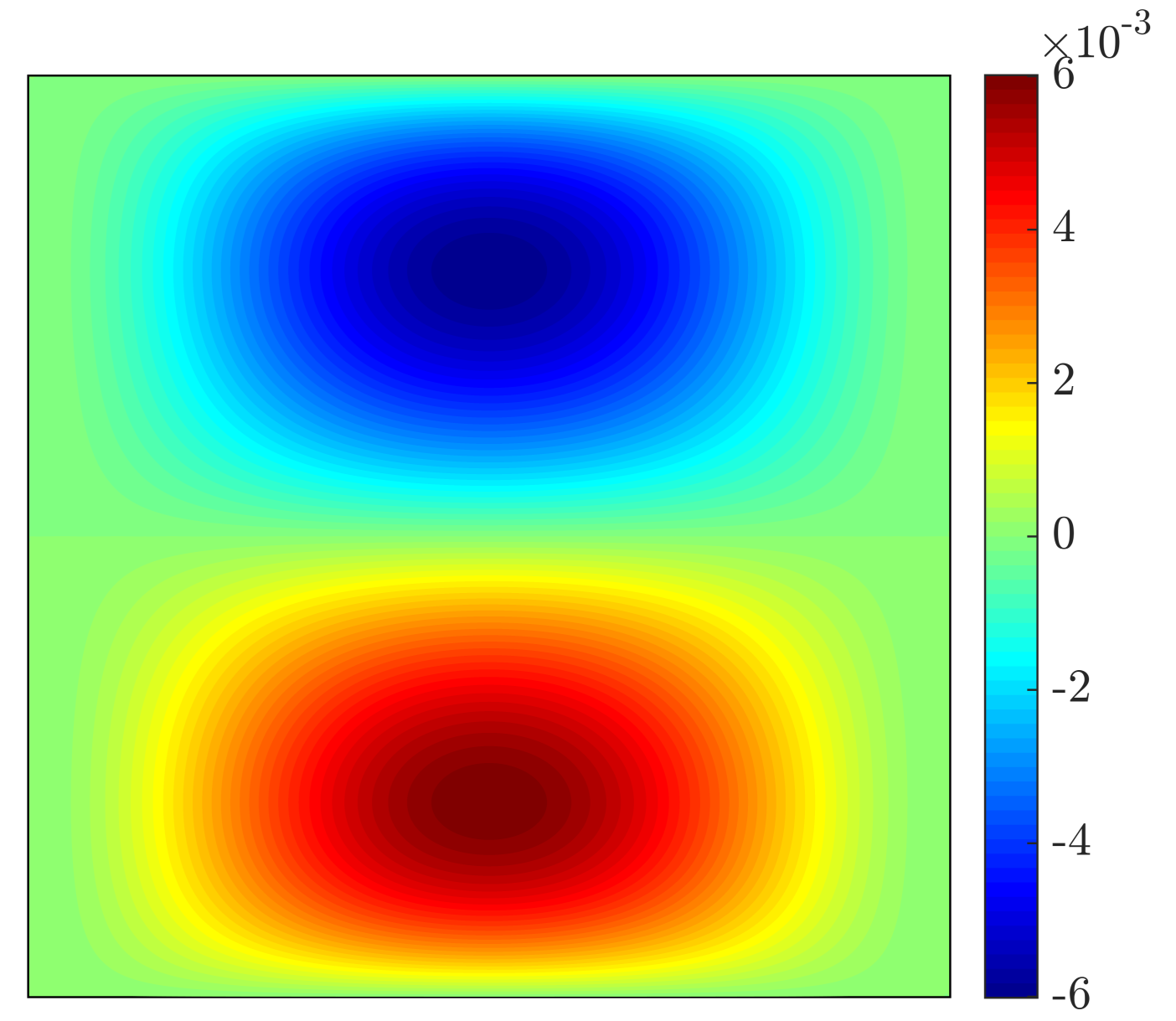}}
	\subfigure[$\sigma_{\texttt{VM}}$]{\includegraphics[width=0.32\textwidth]{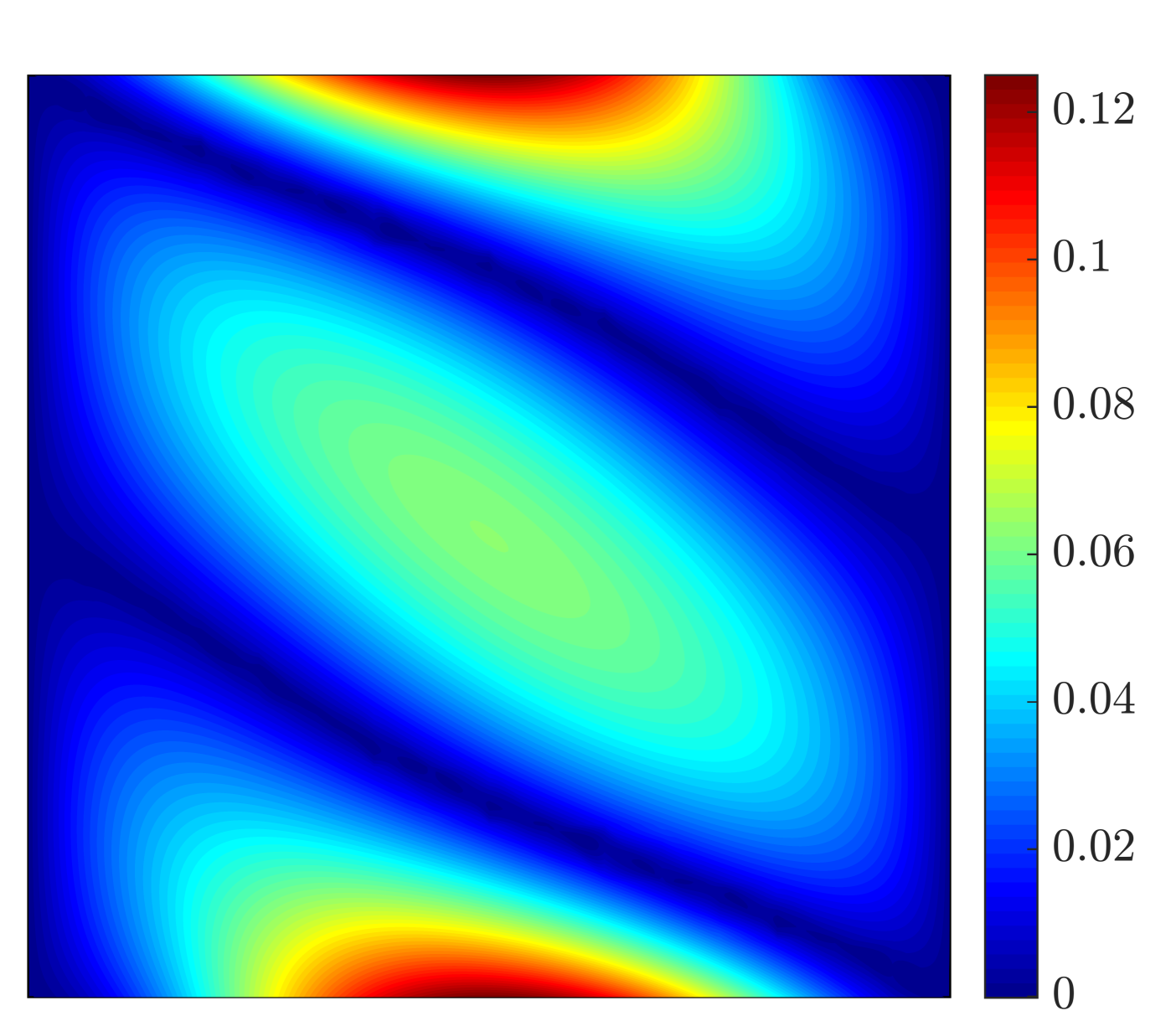}}
	\caption{HDG approximation of the displacement field and the Von Mises stress using the fourth triangular mesh and $k=3$ for a material with $\nu=0.49999$.}
	\label{fig:2DsolIncomp}
\end{figure}

The displacement field and the Von Mises stress computed on the fourth triangular mesh and using a cubic degree of approximation are depicted in Figure~\ref{fig:2DsolIncomp}.

\begin{figure}[!tb]
	\centering
	\subfigure[$\nu=0.49$]{\includegraphics[width=0.4\textwidth]{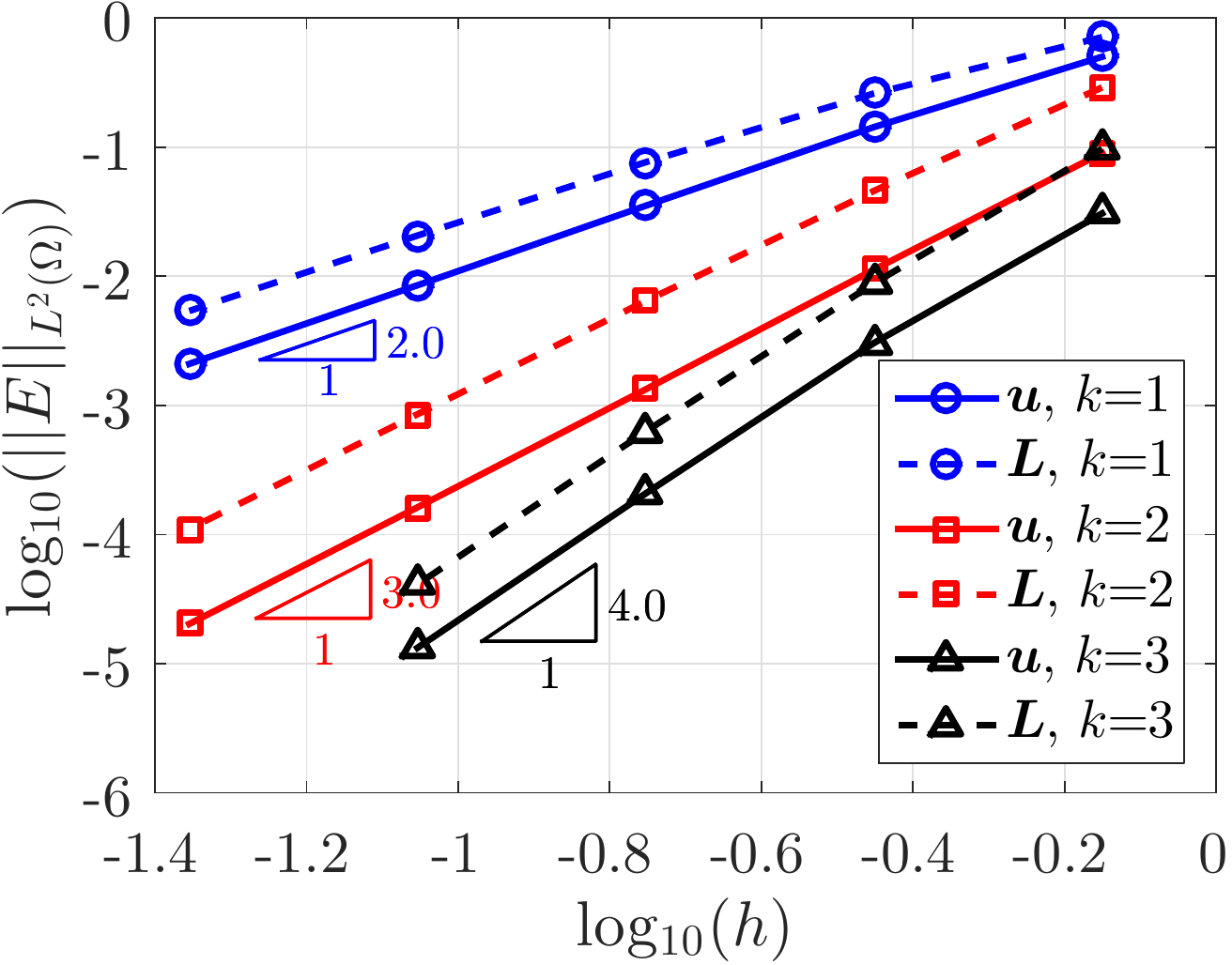}}
	\subfigure[$\nu=0.499$]{\includegraphics[width=0.4\textwidth]{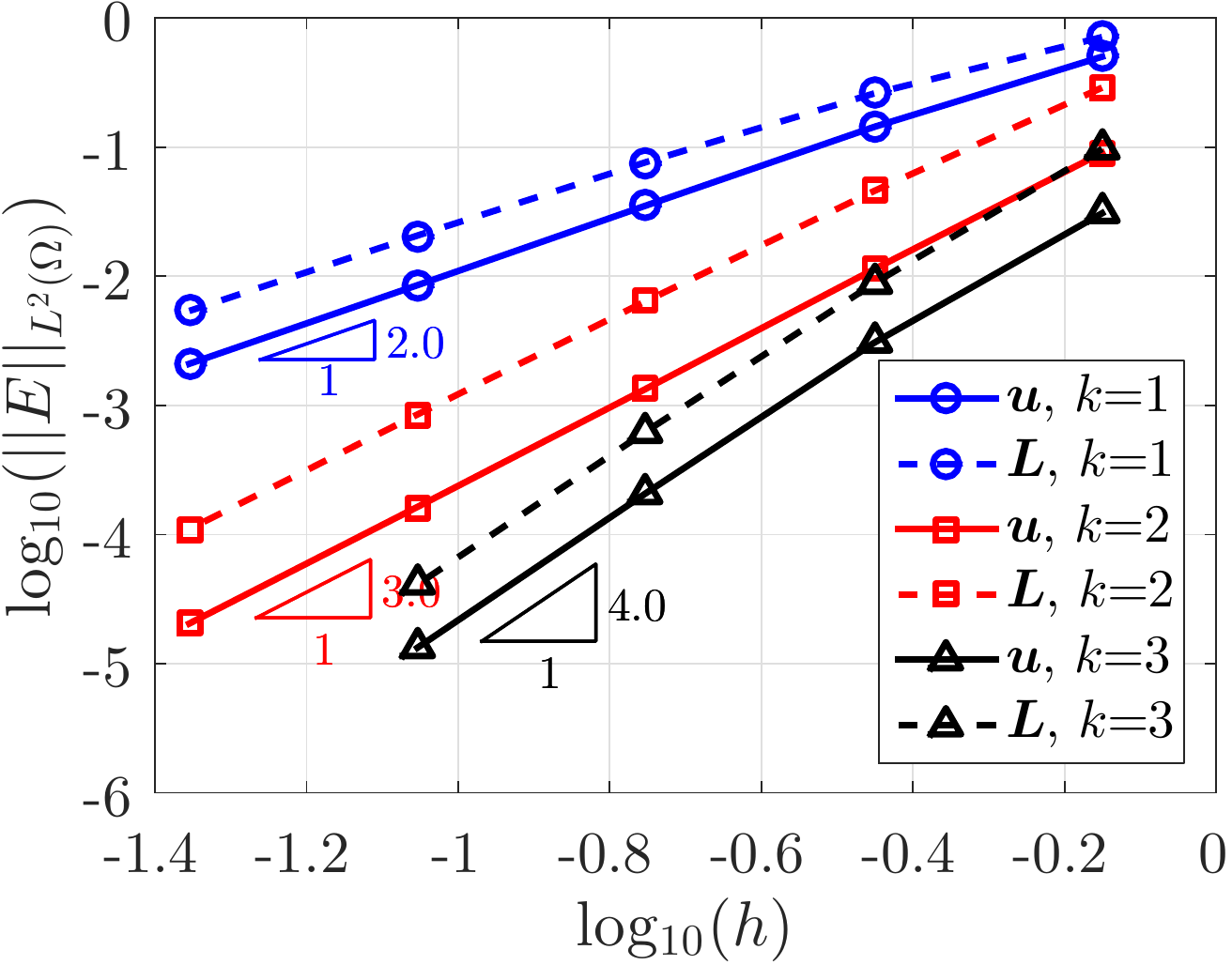}}
	\subfigure[$\nu=0.4999$]{\includegraphics[width=0.4\textwidth]{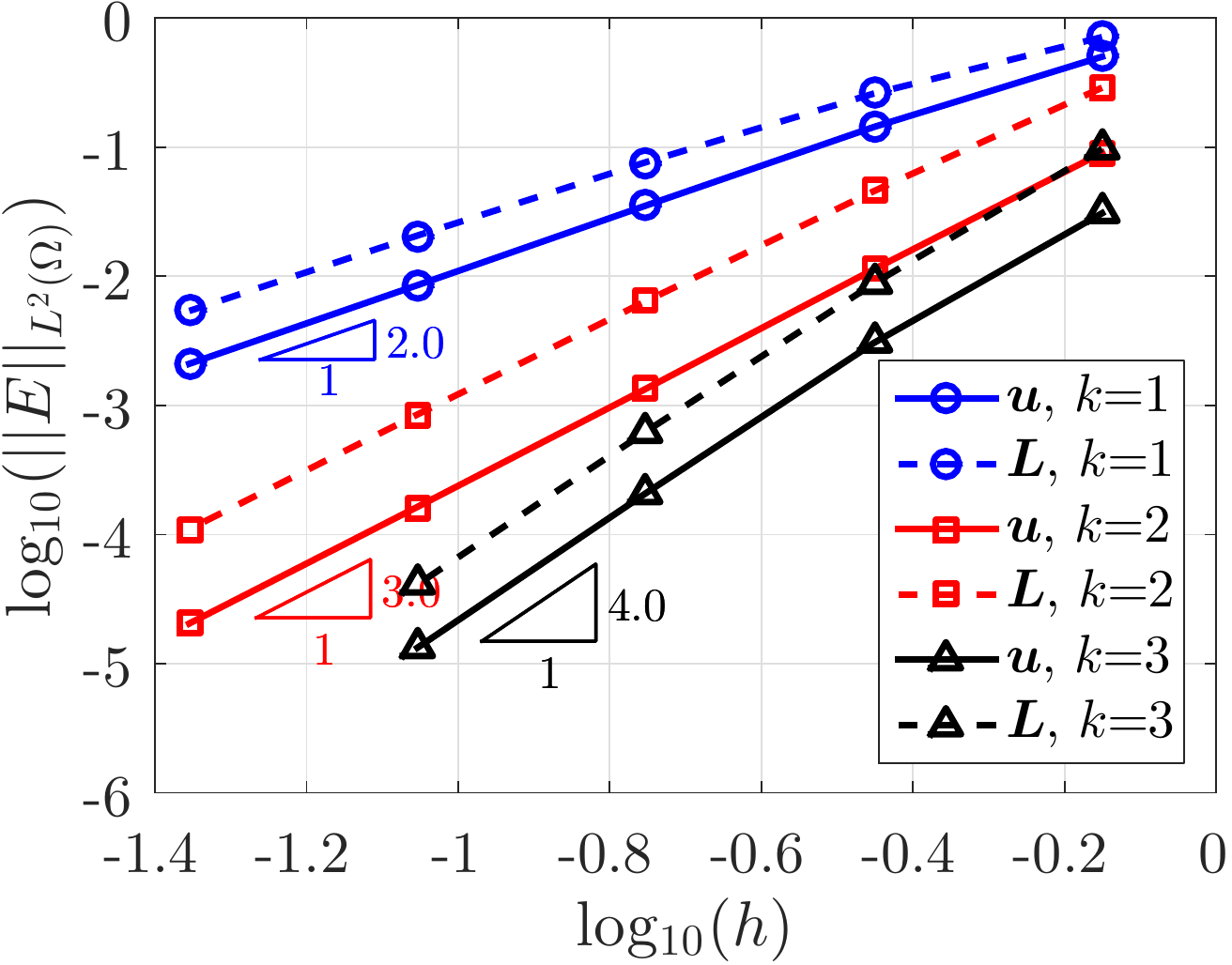}}
	\subfigure[$\nu=0.49999$]{\includegraphics[width=0.4\textwidth]{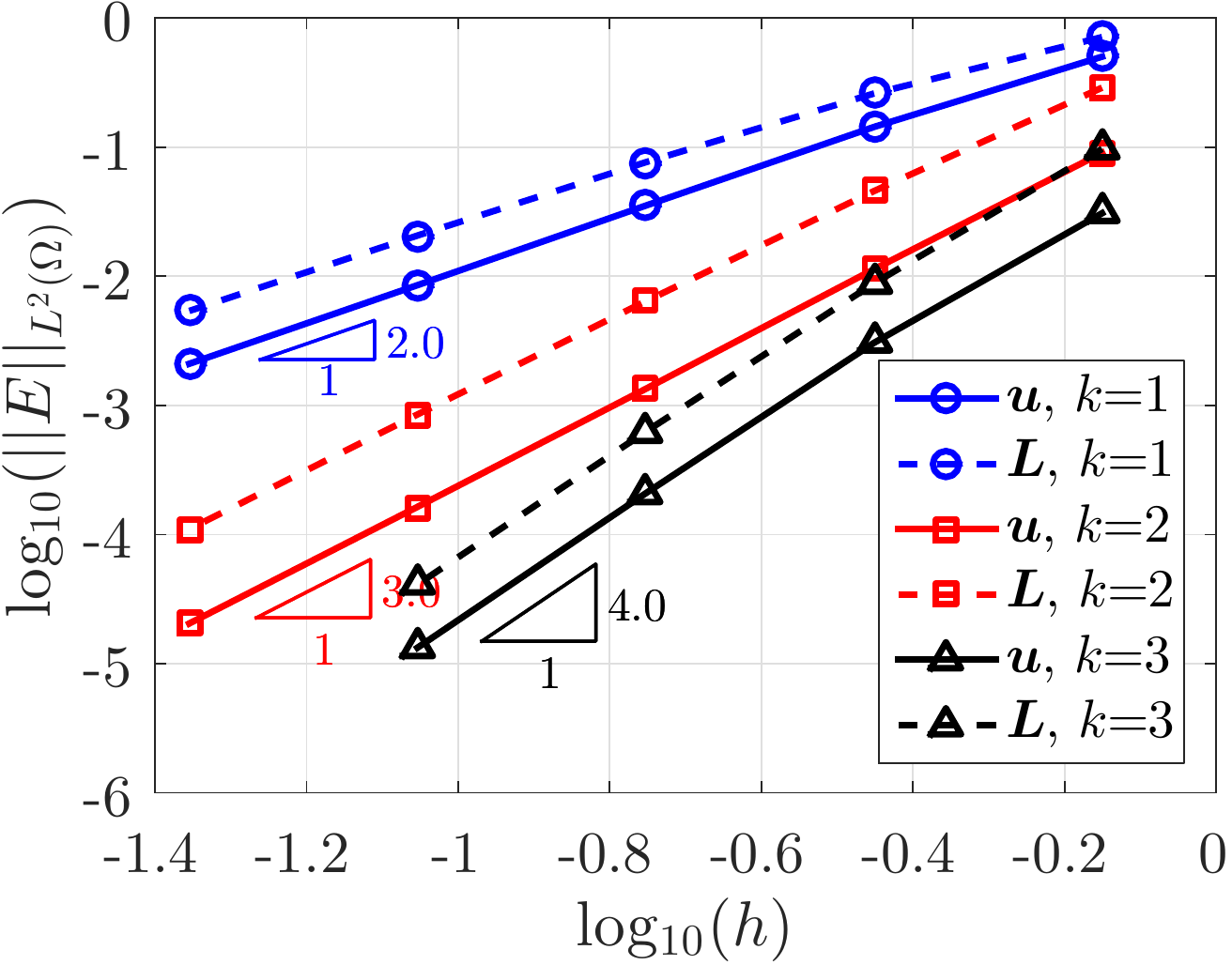}}
	\caption{$h$--convergence of the error of the primal and mixed variables, $\bu$ and $\bL$, in the $\eltwo(\Omega)$ norm for different orders of approximation and for an increasing value of the Poisson's ratio.}
	\label{fig:incompressibleUL}
\end{figure}

The mesh convergence results for the primal and mixed variables, $\bu$ and $\bL$, are represented in Figure~\ref{fig:incompressibleUL} for the triangular meshes shown in Figure~\ref{fig:TRImeshesOneDiag}, for different orders of approximation and for increasing value of the Poisson's ratio.
The results show that the proposed HDG formulation is volumetric locking--free. In addition, it is worth noting that the accuracy of the displacement and the stress is almost independent on the Poisson's ratio. This behaviour has also been observed when using a different HDG formulation~\cite{soon2009hybridizable}. However, contrary to the results reported in~\cite{soon2009hybridizable}, the proposed formulation shows the optimal rate of convergence, whereas the formulation in~\cite{soon2009hybridizable} exhibits a slight degradation of the rate of convergence for nearly incompressible materials. This degradation of the rate of convergence in the HDG formulation of~\cite{soon2009hybridizable} is sizeable when the error on the mixed variable is considered, even for high--order approximations.

\begin{figure}[!tb]
	\centering
	\subfigure[$\nu=0.49$]{\includegraphics[width=0.4\textwidth]{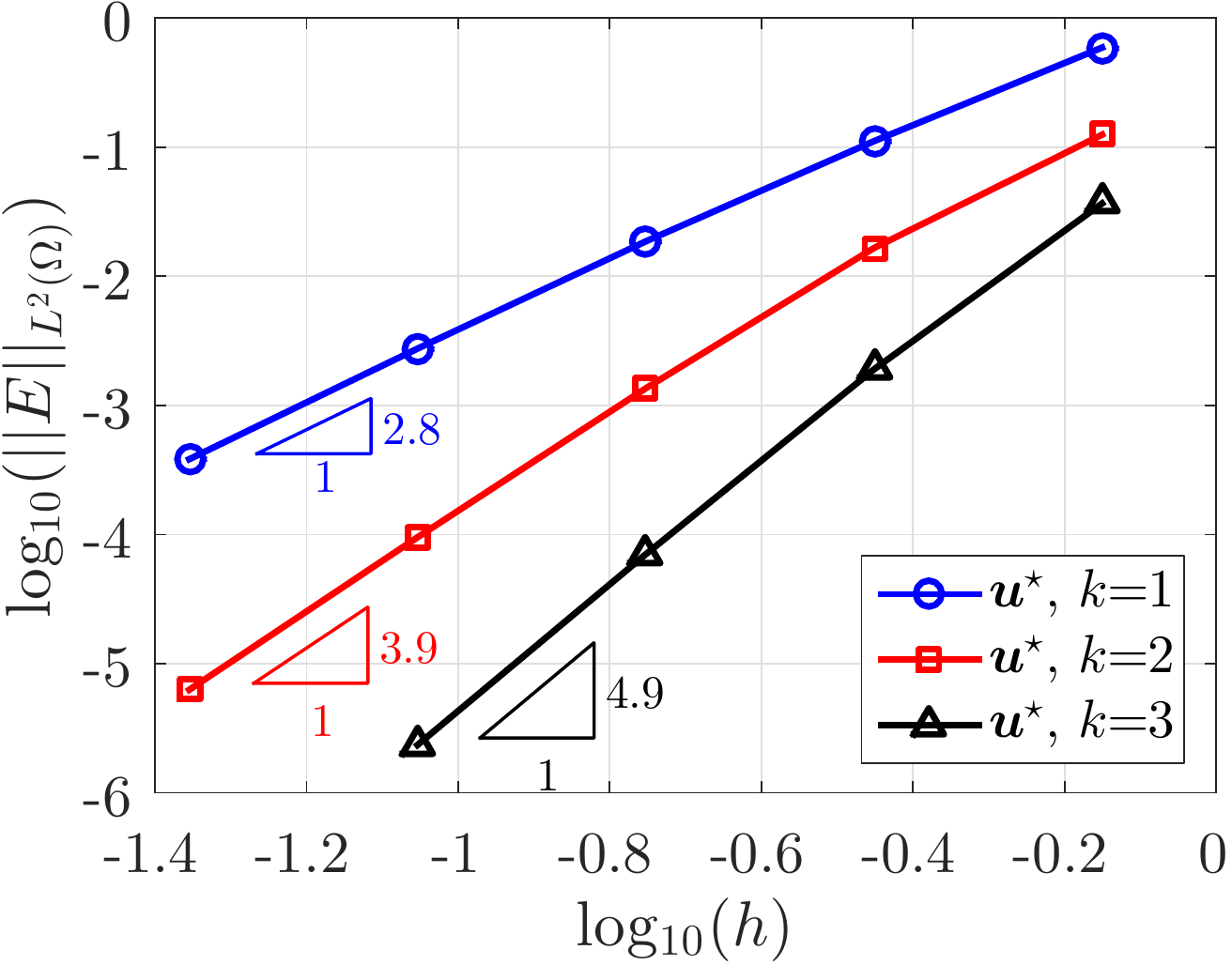}}
	\subfigure[$\nu=0.499$]{\includegraphics[width=0.4\textwidth]{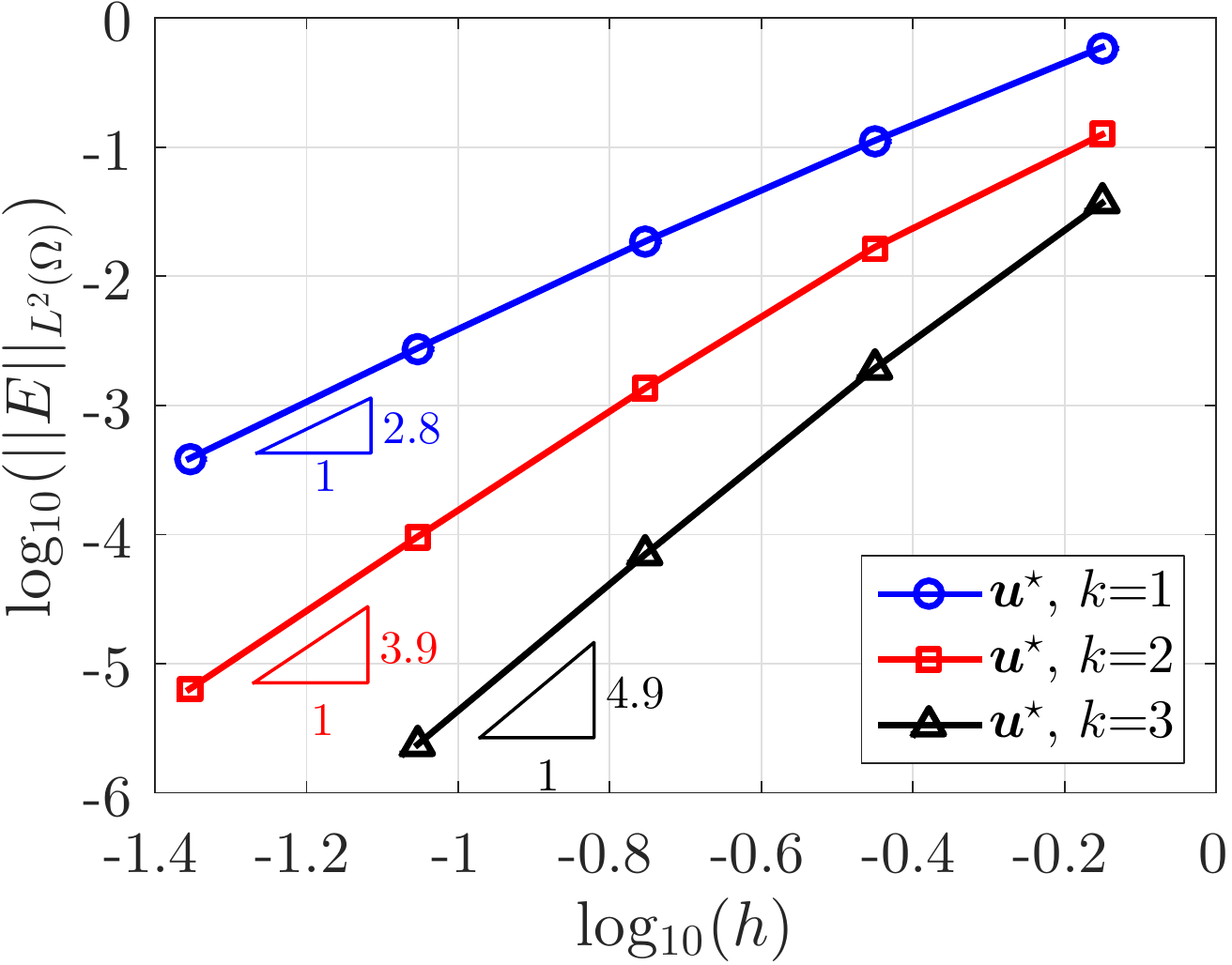}}
	\subfigure[$\nu=0.4999$]{\includegraphics[width=0.4\textwidth]{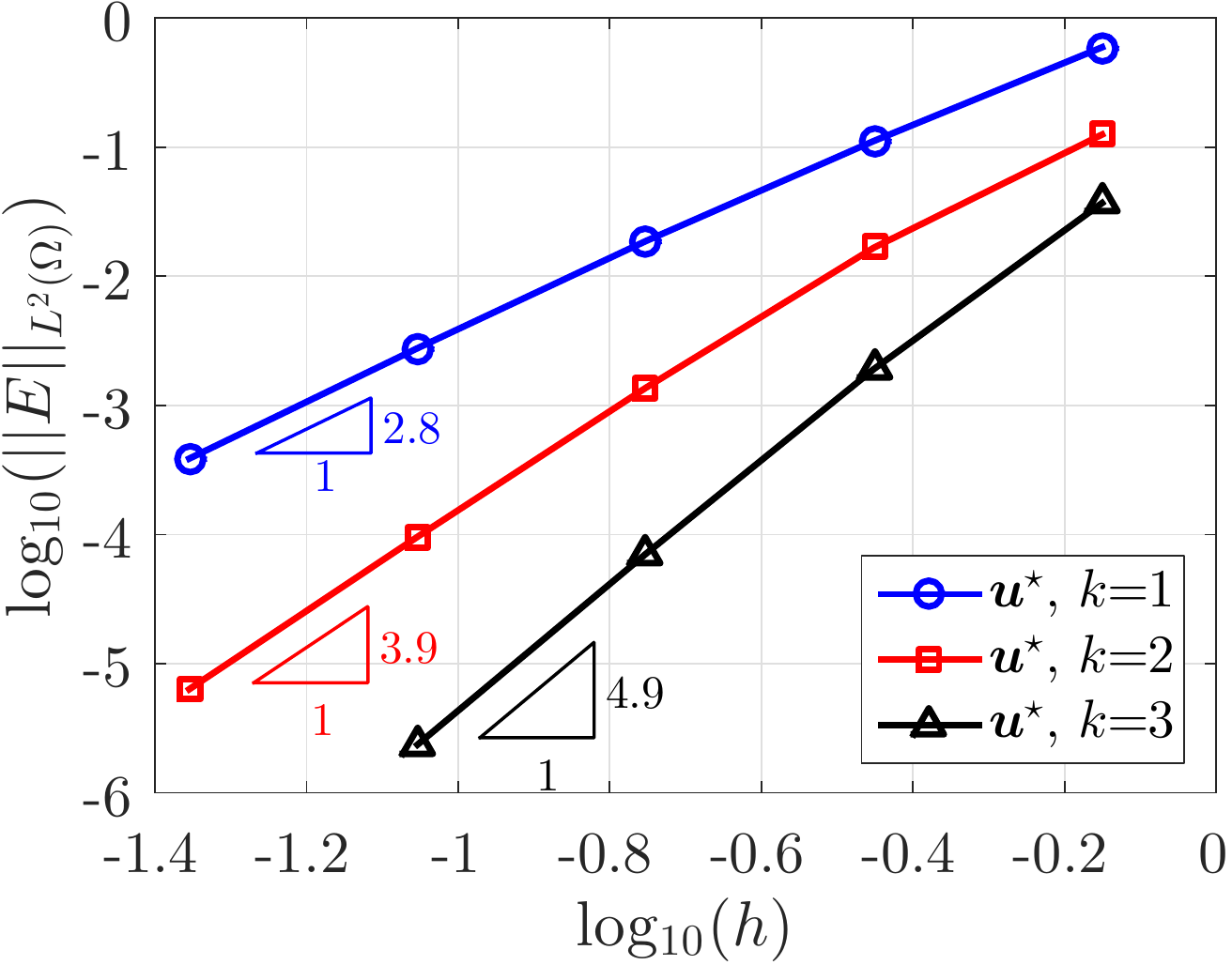}}
	\subfigure[$\nu=0.49999$]{\includegraphics[width=0.4\textwidth]{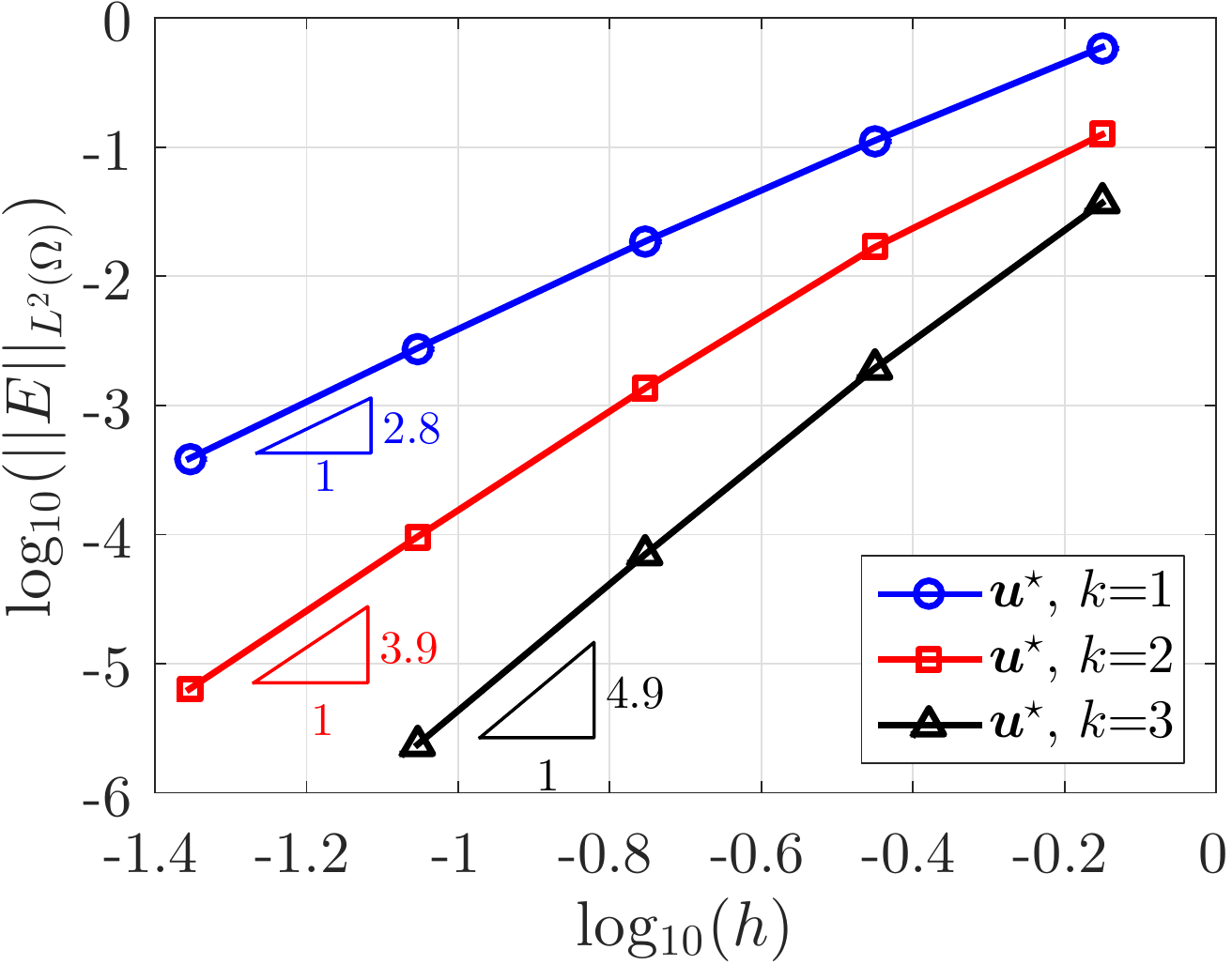}}
	\caption{$h$--convergence of the error of the post--processed variable, $\bu^\star$, in the $\eltwo(\Omega)$ norm for different orders of approximation and for an increasing value of the Poisson's ratio.}
	\label{fig:incompressibleUstar}
\end{figure}

Next, the mesh convergence study is performed for the post--processed variable using the technique proposed in this paper that resulted in optimal convergence in the numerical example of Section~\ref{sc:superconvergenceStudy}.  The mesh convergence results for the post-processed displacement field $\bu^\star$, are represented in Figure~\ref{fig:incompressibleUstar} for different orders of approximation and for increasing value of the Poisson's ratio. 
The results show again that the accuracy is independent on the Poisson's ratio. More important, the mesh convergence study demonstrates that the proposed formulation together with the proposed post--process is able to provide super--convergent solutions for all degrees of approximation, even for linear triangular elements in the particular arrangement that causes volumetric locking in a continuous Galerkin formulation. 

\section{Concluding remarks}
\label{sc:Conclusion}

This paper proposes a novel HDG formulation for the linear elastic problem strongly enforcing the symmetry of the stress tensor. 
Owing to the Voigt notation, the second--order tensors appearing in the linear elasticity equation are expressed as vectors featuring the diagonal and half of the off-diagonal terms.
Thus, the resulting method does not introduce any extra cost to guarantee the symmetry of the stress tensor and, in fact, it is more computationally efficient than other HDG formulations due to the reduced number of degrees of freedom of the mixed variable.

As all existing HDG formulations for linear elasticity, the resulting method provides optimal convergence rate of order $k+1$ for the displacement field. 
The optimal order $k+1$ is also obtained for the stress tensor which usually experiences sub--optimal behaviour using low--order approximations in the original HDG formulation by Cockburn and co--workers.
Furthermore, contrary to other proposed variants of HDG for linear elasticity, the optimality is achieved using equal order approximation spaces for the primal and mixed variables and no special enrichment is required. 

The optimally convergent stress tensor is thus utilised to locally construct a post--processed displacement field.
The element--by--element procedure uses the equilibrated stresses as boundary conditions of the local problems and exploits the optimal convergence of the trace of the displacements to remove the under--determination associated with the rigid rotational modes.
The post--processed displacement field belongs to the richer space of polynomials of degree at most $k+1$ in each element and super--converges with order $k+2$.
Therefore, the current formulation provides a workaround to avoid the construction of discrete spaces fulfilling the $\bm{M}$-decomposition property to guarantee the super--convergence of the post--processed solution.

An extensive set of numerical simulations has been presented to verify the optimal approximation properties of the method in 2D and 3D, to show the robustness of the formulation using meshes of different element types and to study the influence of the HDG stabilisation parameter.
Special attention has been dedicated to the analysis of the limit case of nearly incompressible materials: the method is locking--free and the optimal convergence and super--convergence rates of the primal, mixed and post--processed variables are preserved.

\section*{Acknowledgements}

This work was partially supported by the European Union's Horizon 2020 research and innovation programme under the Marie Sk\l odowska--Curie grant agreement No. 675919 and the Spanish Ministry of Economy and Competitiveness (Grant number: DPI2017-85139-C2-2-R).
The support of the Generalitat de Catalunya (Grant number: 2017SGR1278) is also gratefully acknowledged.
Finally, Alexandros Karkoulias was supported by the European Education, Audiovisual and Culture Executive Agency (EACEA) under the Erasmus Mundus Joint Doctorate Simulation in Engineering and Entrepreneurship Development (SEED), FPA 2013-0043.

\bibliographystyle{abbrv}
\bibliography{Ref-HDG}

\appendix

\section{Implementation details}\label{sc:implementation}

A standard isoparametric formulation is considered, where the approximation of the primal and mixed variables, $\bu^h$ and $\bL^h$, is defined in a reference element $\widetilde{\Omega}$, with local coordinates $\bxi = (\xi_1, \ldots, \xi_{\nsd})$, and the approximation of the hybrid variable, $\bhu^h$, is defined in a reference face $\widetilde{\Gamma}$, with local coordinates $\bet = (\eta_1, \ldots, \eta_{\nsd-1})$, as
\begin{equation*} 
\bu^h(\bxi) = \sum_{j=1}^{\nen} \nodaluV_j N_j(\bxi), \qquad
\bL^h(\bxi) = \sum_{j=1}^{\nen} \nodalLV_j N_j(\bxi), \qquad
\hu^h(\bet) = \sum_{j=1}^{\nfn} \nodaluhV_j \hat{N}_j(\bet),
\end{equation*} 
where $\nen$ and $\nfn$ denote the number of element and face nodes respectively and $N_j$ and $\hat{N}_j$ are the shape functions used to define the approximation within the reference element and face respectively.

The isoparametric transformation is used to relate local and Cartesian coordinates, namely
\begin{equation*}
\bx(\bxi) = \sum_{k=1}^{\nen} \bx_k N_k(\bxi),
\end{equation*}
where $\{\bx_k\}_{k=1,\ldots,\nen}$ denote the elemental nodal coordinates. 

The following matrices are introduced in two dimensions
\begin{equation*}
\mat{E}_1 = 
\begin{bmatrix}
1 & 0 & 0 \\
0 & 0 & 1
\end{bmatrix}^T
\qquad
\mat{E}_2 = 
\begin{bmatrix}
0 & 0 & 1 \\
0 & 1 & 0
\end{bmatrix}^T.
\end{equation*}
Similarly, in three dimensions, the following matrices are defined
\begin{equation*}
\mat{E}_1 = 
\begin{bmatrix}
1 & 0 & 0 & 0 & 0 & 0 \\
0 & 0 & 0 & 1 & 0 & 0 \\
0 & 0 & 0 & 0 & 1 & 0
\end{bmatrix}^T
\qquad
\mat{E}_2 = 
\begin{bmatrix}
0 & 0 & 0 & 1 & 0 & 0 \\
0 & 1 & 0 & 0 & 0 & 0 \\
0 & 0 & 0 & 0 & 0 & 1
\end{bmatrix}^T
\qquad
\mat{E}_3 = 
\begin{bmatrix}
0 & 0 & 0 & 0 & 1 & 0 \\
0 & 0 & 0 & 0 & 0 & 1 \\
0 & 0 & 1 & 0 & 0 & 0
\end{bmatrix}^T.
\end{equation*}
These matrices are used to express, in compact form, the matrices $\gradS$ and  $\bN$, defined in Equations~\eqref{eq:symmGrad} and~\eqref{eq:normalVoigt} respectively, as
\begin{equation*}
\gradS = \sum_{k=1}^{\nsd} \mat{E}_k \frac{\partial }{\partial x_k}, \qquad  \bN = \sum_{k=1}^{\nsd} \mat{E}_k n_k
\end{equation*}

In addition, the following compact form of the shape functions is introduced
\begin{equation*}
\Nmat = \begin{bmatrix} N_1\Insd & N_2\Insd & \dots & N_{\nen}\Insd \end{bmatrix}^T, \qquad
\Mmat = \begin{bmatrix} N_1\Imsd & N_2\Imsd & \dots & N_{\nen}\Imsd \end{bmatrix}^T,
\end{equation*}
\begin{equation*}
\Nmat_k = \begin{bmatrix} \frac{\partial N_1}{\partial x_k} \bm{E}_k \bDHalf & \frac{\partial N_2}{\partial x_k} \bm{E}_k \bDHalf & \dots & \frac{\partial N_{\nen}}{\partial x_k} \bm{E}_k \bDHalf \end{bmatrix}^T, \quad \text{for } k=1,\ldots,\nsd,
\end{equation*}
\begin{equation*}
\Nmat^n_k = \begin{bmatrix} N_1 n_k \bm{E}_k \bDHalf & N_2 n_k \bm{E}_k \bDHalf & \dots &  N_{\nfn} n_k \bm{E}_k \bDHalf \end{bmatrix}^T, \quad \text{for } k=1,\ldots,\nsd.
\end{equation*}
\begin{equation*}
\NmatHat = \begin{bmatrix} \hat{N}_1\Insd & \hat{N}_2\Insd & \dots & \hat{N}_{\nfn}\Insd \end{bmatrix}^T, \qquad 
\NmatHat_{\tau} = \begin{bmatrix} \hat{N}_1\btau & \hat{N}_2\btau & \dots & \hat{N}_{\nfn}\btau \end{bmatrix}^T,
\end{equation*}
\begin{equation*}
\NmatHat^n_k = \begin{bmatrix} \hat{N}_1 n_k \bm{E}_k \bDHalf & \hat{N}_2 n_k \bm{E}_k \bDHalf & \dots &  \hat{N}_{\nfn} n_k \bm{E}_k \bDHalf \end{bmatrix}^T, \quad \text{for } k=1,\ldots,\nsd.
\end{equation*}

The matrices and vectors resulting from the discretisation of Equation~\eqref{eq:HDGElasticityWeakLocalL} of the local problem are
\begin{equation*}
[\mat{A}_{LL}]_e = -\sumge \Mmat(\bxige) \Mmat^T(\bxige) |\bJ(\bxige) | \wge ,
\end{equation*}
\begin{equation*}
[\mat{A}_{Lu}]_e = \sum_{k=1}^{\nsd} \sumge \Nmat_k(\bxige) \Nmat^T(\bxige) |\bJ(\bxige) | \wge ,
\end{equation*}
\begin{equation*}
[\mat{A}_{L \hat{u}}]_e = \sum_{f=1}^{\numfa^e} \left( \sum_{k=1}^{\nsd} \sumgf \Nmat^n_k (\bxigf) \NmatHat^T(\bxigf) |\bJ(\bxigf) | \wgf \right) \left(1 - \chi_{\Gamma_D}(f) \right),
\end{equation*}
\begin{equation*}
[\vect{f}_{L}]_e = \sum_{f=1}^{\numfa^e} \left( \sum_{k=1}^{\nsd} \sumgf \Nmat^n_k(\bxigf) \bu_D \left(\bx(\bxigf) \right) |\bJ(\bxigf) | \wgf \right) \chi_{\Gamma_D}(f) ,
\end{equation*}
where $\numfa^e$ is the number of faces, $\Gamma_{e,j}$ for $j=1,\ldots,\numfa^e$ of the element $\Omega_e$ and $\chi_{\Gamma_D}$ is the indicator function of $\Gamma_D$, i.e.
\begin{equation*}
\chi_{\Gamma_D}(f) = 
\Bigg\{
\begin{array}{ll}
1 & \text{ if } \Gamma_{e,j} \cap \Gamma_D \neq \emptyset  \\
0 & \text{ otherwise}
\end{array}.
\end{equation*}

In the above expressions, $\bxige$ and $\wge$ are the $\nipe$ integration points and weights defined on the reference element and $\bxigf$ and $\wgf$ are the $\nipf$ integration points and weights defined on the reference face.

Similarly, the matrices and vectors resulting from the discretisation of Equation~\eqref{eq:HDGElasticityWeakLocalU} of the local problem are
\begin{equation*}
[\mat{A}_{uu}]_e = \sum_{f=1}^{\numfa^e}  \sumgf \NmatHat(\bxigf) \NmatHat_{\tau}^T(\bxigf) |\bJ(\bxigf) | \wgf,
\end{equation*}
\begin{equation*}
[\mat{A}_{u \hat{u}}]_e = \sum_{f=1}^{\numfa^e} \left(  \sumgf \Nmat(\bxigf) \NmatHat_{\tau}^T(\bxigf) |\bJ(\bxigf) | \wgf \right)  \left(1 - \chi_{\Gamma_D}(f) \right),
\end{equation*}
\begin{equation*}
[\vect{f}_{u}]_e = \sumge \Nmat(\bxige) \bm{f} \left(\bx(\bxige) \right) |\bJ(\bxige) | \wge + \sum_{f=1}^{\numfa^e} \left(  \sumgf \Nmat(\bxigf) \btau \bu_D \left(\bx(\bxigf) \right) |\bJ(\bxigf) | \wgf \right) \chi_{\Gamma_D}(f).
\end{equation*}

Finally, the matrices and vectors resulting from the discretisation of Equation~\eqref{eq:HDGElasticityWeakGlobal} of the local problem are
\begin{equation*}
[\mat{A}_{\hat{u} \hat{u}}]_e = -\sum_{f=1}^{\numfa^e} \left(  \sumgf \NmatHat(\bxigf) \NmatHat_{\tau}^T(\bxigf) |\bJ(\bxigf) | \wgf \right)  \left(1 - \chi_{\Gamma_D}(f) \right),
\end{equation*}
\begin{equation*}
[\vect{f}_{\hat{u}}]_e = - \sum_{f=1}^{\numfa^e} \left(  \sumgf \NmatHat(\bxigf)  \bg \left(\bx(\bxigf) \right) |\bJ(\bxigf) | \wgf \right) \chi_{\Gamma_N}(f),
\end{equation*}
where $\chi_{\Gamma_N}$ is the indicator function of $\Gamma_N$.

\end{document}